\documentclass[11pt]{book}
\usepackage{amsmath,amsfonts,latexsym,amssymb,amscd, graphicx,pictexwd}

\renewcommand{\P}{{\mathbb P}}
\newcommand{\E}{{\mathbb E}}
\newcommand{\R}{{\mathbb R}}
\newcommand{\LL}{{\mathbb L}}

\newcommand{\ZZ}{{\mathbb Z}}
\newcommand{\PP}{{\mathbf P}}

\newcommand{\0}{{\mathbf 0}}
\newcommand{\pp}{{\mathbf p}}
\newcommand{\qq}{{\mathbf q}}
\newcommand{\cc}{{\mathbf c}}
\newcommand{\zz}{{\mathbf z}}
\newcommand{\yy}{{\mathbf y}}
\newcommand{\mm}{{\mathbf m}}
\newcommand{\xx}{{\mathbf x}}
\newcommand{\nn}{{\mathbf n}}

\newcommand{\bft}{{\mathbf t}}
\renewcommand{\cal}{\mathcal}

\newcommand{\Co}{{\rm Co}}
\newcommand{\Cy}{{\rm Cyl}}

\newcommand{\calN}{{\mathcal N}}

\newcommand{\calF}{{\mathcal F}}

\newcommand{\eps}{{\varepsilon }}

\newcommand{\Var}{{\rm Var}}
\newcommand{\Cov}{{\rm Cov}}
\newcommand{\la}{{\lambda}}

\newtheorem{thm}{Theorem}[section]
\newtheorem{coro}[thm]{Corollary}
\newtheorem{lem}[thm]{Lemma}
\newtheorem{prop}[thm]{Proposition}
\numberwithin{equation}{section}

\begin{document}



\title{Hydrodynamical Methods in \\ Last Passage Percolation Models}

\author{Eric A. Cator and Leandro P. R. Pimentel}




\maketitle



\chapter*{Preface}

These lecture notes are written as reference material for the Advanced Course ``Hydrodynamical Methods in Last Passage Percolation Models'', given at the $28^{\rm th}$ \emph{Col\'oquio Brasileiro de Matem\'atica} at IMPA, Rio de Janeiro, July 2011.\\

The notes focus on recent progress made in two dimensional Poisson Last Passage Percolation (LPP), using what is generally known as ``hydrodynamical methods''\cite{CaGr1,CaGr2,CPshape,CaPi,CaPi2,CaPi3}. These methods have been developed in recent years by the authors, but also by many others, such as Groeneboom, Sepp\"al\"ainen, Bal\'azs, Aldous, Diaconis, Martin and Ferrari. The methods try to use probabilistic tools to analyze last passage percolation, such as Markov processes, interacting particle processes, concentration inequalities etc. This is somewhat in contrast to the analytic methods used by Baik, Deift and Johansson in their ground breaking paper \cite{BaDeJo}, which were further developed by many other authors as well, such as for example Rains, Pr\"aehofer, Spohn and O'Connell. The analytic methods (which we will not discuss in any detail) are able to get strong results on limiting distributions and fluctuation probabilities of the length of a longest path in the characteristic direction for some specific LPP models, using exact formulas for certain distributions and analyzing the limiting behavior of rescaled versions of these formulas. The hydrodynamical methods are not able to find these limiting distributions (so far), but they were able to find the correct scaling for these specific models, and extend this scaling to the transversal fluctuations of the paths, giving information on second class particles, for example. Furthermore, it turns out that these methods can be (partially) extended to models for which there are no analytical methods available yet. Also, certain tightness properties of rescaled functionals of the LPP model can be established using the hydrodynamical methods, whereas limiting distributions for these functionals can be determined using the analytic methods.\\

 The main object of study is LPP using iid weights on a homogeneous Poisson process in $\R^2$. The strongest results will be valid for the case where all weights are equal to 1 (called the classical Hammersley process). Almost all results can (and have been) extended to LPP on ${\mathbb Z}^2$ with iid exponential weights, which corresponds to the Totally Asymmetric Exclusion process. However, for the purpose of this one week course, we will not go into this specific model.\\
 
Chapter \ref{ch:poissonLPP} introduces the model and shows the existence of Busemann functions, a concept originally developed for Riemannian geometry, but which proves to be very useful in our context. In Chapter \ref{ch:fluid} we introduce the Hammersley Interacting Fluid System, a generalization of the Hammersley Interacting Particle System, which turns out to be an important tool in studying the Poisson LPP. Chapter \ref{ch:Busemannequi} connects the Busemann functions and the fluid process in the sense that the equilibrium measures for the fluid process can be described by the Busemann functions. This leads to a description of the Busemann function in the classical Hammersley process and to a new description of the multi-class particle (or fluid) system. Chapter \ref{ch:2ndclass} introduces second-class particles in the fluid process, and shows they satisfy the usual strong law when in equilibrium, or close to equilibrium. It also shows how the Busemann functions can be used to describe the asymptotic speed of second class particles in a rarefaction fan. Finally, Chapter \ref{ch:CLT} introduces some remarkable identities that link the variance of the longest path to the variance of the starting configuration. In the classical Hammersley process, this leads to normal central limit theorems for longest paths in non-characteristic directions, and to cube-root fluctuations of the longest path along the characteristic direction. \\


\tableofcontents

\chapter{The Poisson Last Passage Percolation Model}\label{ch:poissonLPP}
\section{The Graphical Construction}
Let $\PP\subseteq\R^2$ be a two-dimensional Poisson random set of intensity one. On each point $\pp\in\PP$ we put a random positive weight $\omega_\pp$ and we assume that $\{\omega_\pp\,:\,\pp\in\PP\}$ is a collection of i.i.d. random variables, distributed according to a distribution function $F$, which are also independent of $\PP$. When $F$ is the Dirac distribution concentrated on $1$ (each point has weight $1$; we will denote this $F$ by $\delta_1$), then we refer to this model as the classical Hammersley model \cite{AlDi}. \\

For each $\pp,\qq\in\R^2$, with $\pp<\qq$ (inequality in each coordinate, $\pp\neq \qq$), let $\Pi(\pp,\qq)$ denote the set of all increasing (or up-right) paths, consisting of points in $\PP$, from $\pp$ to $\qq$, where we exclude all points that share (at least) one coordinate with $\pp$. So we consider the points in the rectangle $(\pp,\qq]$, where we leave out the south and the west side of the rectangle. In this probabilistic model, the ``metric'' (or last passage time) $L$ between  $\pp\leq\qq$ is defined by
$$L(\pp,\qq):=\max_{\varpi\in\Pi(\pp,\qq)}\big\{\sum_{\pp'\in\varpi}\omega_{\pp'}\big\}\,.$$
Then $L$ is superadditive,
$$L(\pp,\qq)\geq L(\pp,\zz)+L(\zz,\qq)\,.\\$$

When we consider a path $\varpi$ from $\pp$ to $\qq$ consisting of increasing points $(\pp_1,\ldots,\pp_n)$, we will view $\varpi$ as the lowest increasing continuous path connecting all the points, starting at $\pp$ and ending at $\qq$, and then excluding $\pp$. This way we can talk about crossings with other paths or with lines. Suppose $\varpi_1$ and $\varpi_2$ are two different paths between $\pp$ and $\qq$ that both have the maximal weight $L(\pp,\qq)$. Now suppose that $(x_1,t_1)\leq (x_2,t_2)$ are two points on the intersection of the two paths, such that the paths do not intersect in the open strip $(x_1,x_2)\times \R$. Then one path must lie entirely below the other path in this strip, and the total weight of the two paths in this strip must be equal (otherwise we would be able to construct a ``longer'' path). This means that we can define the minimum of $\varpi_1$ and $\varpi_2$ by choosing the lowest path in each such strip.  The finite geodesic between $\pp$ and $\qq$ is given by the lowest path (in the sense we just described) that attains the maximum in the definition of $L(\pp,\qq)$. We will denote this geodesic by $\varpi(\pp,\qq)$ (this is well defined for any ordered pair $(\pp,\qq)$, even if we do not specify the order).\\

\section{The shape Theorem}
A crucial result for LPP models is the following shape theorem: set $\0=(0,0)$, $\nn=(n,n)$,
\[F(x) = \P(\omega_\pp \leq x)\,\mbox{ and }\,\gamma=\gamma(F)=\sup_{n\geq 1}\frac{\E L(\0,\nn)}{n}> 0 \,.\]
\begin{thm}\label{thm:shape}
Suppose that
\begin{equation}\label{eq:a1}
\int_0^\infty \sqrt{1-F(x)}\,dx <+\infty\,.
\end{equation}
Then $\gamma(F)<\infty$ and for all $x,t>0$, as $r \to \infty$,
$$\frac{\E L\left(\0,(rx,rt)\right)}{r} \to \gamma \sqrt{xt}\,,$$
and, a.s., 
\begin{equation}\label{eq:shape}
\frac{L\left(\0,(rx,rt)\right)}{r}\to \gamma \sqrt{xt}\,.
\end{equation}
\end{thm}

Theorem \ref{thm:shape} shows that $L$ has a curved limiting shape, mainly due to the invariance of the Poisson process under volume preserving maps: if $x,t,\lambda >0$ and $\pp\in \R^2$, then
\begin{equation}\label{eq:sym}
L(\0,(x,t)) \stackrel{\cal D}{=} L\left(\pp,\pp + (\lambda x,t/\lambda)\right)\,
\end{equation}
(choose $\lambda=\sqrt{t/x}$ and $\pp=\0$). \\

\noindent{\bf Proof:} Equation \eqref{eq:sym} shows that it is enough to prove the theorem for $(x,t)=(1,1)$. When considering only one ray, the convergence of $L$ is a standard consequence of Liggett's version of the Subadditive Ergodic Theorem \cite{Li}, as soon as we can show that
\[ \limsup_{r \to \infty} \frac{\E L\left(\0,(r,r)\right)}{r} <\infty\,.\]
Recall that $\delta_1$ denotes the Dirac distribution concentrated on $1$. For each $p\in[0,1]$ denote by $\E_p$ expectation for the Hammersley last passage model induced by Bernoulli weights $w'_{\pp}$, where $\P(w'_\pp = 1) = p$. This coincides with the classical Hammersley model, but with Poisson intensity $p$ (instead of $1$). Thus, it is not hard to see that 
\[ \lim_{r \to \infty} \frac{\E_pL\left(\0,(r,r)\right)}{r} = \gamma(\delta_1)\sqrt{p}\,,\]
for some $\gamma(\delta_1)<\infty$\footnote{It is actually known that $\gamma(\delta_1)=2$ \cite{AlDi}.}. Now we use an idea introduced in \cite{Ma}:
\begin{eqnarray*}
L(\0,\pp) & = & \max_{\varpi\in\Pi(\0,\pp)}\big\{\sum_{\pp'\in\varpi}w_{\pp'}\big\}\\
& = & \max_{\varpi\in\Pi(\0,\pp)}\big\{ \int_0^\infty \sum_{\pp'\in\varpi}1_{\{w_{\pp'}>x\}} \,dx\big\}\\
& \leq & \int_0^\infty \max_{\varpi\in\Pi(\0,\pp)}\big\{\sum_{\pp'\in\varpi}1_{\{w_{\pp'}>x\}} \big\}\,dx\,.
\end{eqnarray*}
The integrand in the last line corresponds to the Bernoulli model with $p=1-F(x)$. This means that
\begin{eqnarray*}
\limsup_{r \to \infty} \frac{\E L\left(\0,(r,r)\right)}{r} & \leq & \limsup_{r \to \infty} \int_0^\infty \frac{\E_{1-F(x)}L\left(\0,(r,r)\right)}{r}\, dx \\
& = & \gamma(\delta_1) \int_0^\infty \sqrt{1-F(x)}\,dx\,.
\end{eqnarray*}

\hfill $\Box$\\

\begin{thm}\label{thm:fluctuation}
If \eqref{eq:a1} is strengthened to: there exists $a>0$ such that
\begin{equation}\label{eq:a2}
 \int_0^\infty \exp(a x)\,dF(x) <+\infty\,,
\end{equation}
then there exist constants $a_0,a_1,a_2,a_3,a_4>0$ such that for all $r\geq a_0$
\begin{equation}\label{eq:kesten}
\P\big(|L(\0,(r,r))-\gamma r |\geq u\big)\leq a_1\exp\Big(-a_2\frac{u}{\sqrt{r}\log r}\Big)\,
\end{equation}
for  $u\in \left[\,a_3 \sqrt{r}\log^2 r\,,\,a_4 r^{3/2}\log r\,\right]$.
\end{thm}

The proof of \eqref{eq:kesten}, which was first stated in \cite{CPshape}, follows Kesten's approach developed for lattice First Passage Percolation models \cite{Ke}, which is based on the Method of Bounded Increments applied to $L$ (see Lemma \ref{Kesten} below). This gives a bound on the fluctuation of $L$ around its expectation. Then adapting a clever argument used by Howard and Newman in \cite{HoNe} shows that one can replace $\E L(\0,(r,r))$ by the shape function.
\begin{lem}\label{Kesten}
Let $\{\calF_k\}_{0\leq k\leq N}$ be a filtration and let $\{U_k\}_{0\leq k\leq N}$ be a family of positive random variables that are $\calF_N$ measurable. Let $\{M_k\}_{0\leq k\leq N}$ be a martingale with respect to $\{\calF_k\}_{0\leq k\leq N}$. Assume that the increments $\Delta_k:=M_k-M_{k-1}$ satisfy
\begin{equation}\label{hip1}
|\Delta_k|\leq c_0\mbox { for some }c_0>0
\end{equation}
and
\begin{equation}\label{hip2}
\E \left(\Delta_k^2\mid\calF_{k-1}\right)\leq \E \left(U_k\mid\calF_{k-1}\right)\,.
\end{equation}
Assume further that for some constants $0<c_1,c_2<\infty$ and $x_0\geq e^2c^2$ we have
\begin{equation}\label{hip3}
\P\left(\sum_{k=1}^N U_k>x\right)\leq c_1 \exp(-c_2 x)\mbox{ when }x\geq x_0\,.
\end{equation}
Then irrespective of the value of $N$, there exists universal constants $0<c_3,c_4<\infty$ that do not depend on $N,c_0,c_1,c_2$ and $x_0$, nor on the distribution of $\{M_k\}_{0\leq k\leq N}$ and $\{U_k\}_{0\leq k\leq N}$, such that
\begin{equation}\label{bound1}
\P(M_N-M_0\geq x)\leq c_3\left\{1+c_1+\frac{c_1}{c_2x_0}\right\}\exp\left(-c_4\frac{x}{\sqrt{x_0}}\right)\,
\end{equation}
whenever $x\leq c_2 \,x_0^{3/2}$.
\end{lem}
\noindent{\bf Proof:} See Theorem 3 in \cite{Ke}.

\hfill $\Box$\\

We decompose $L$ as a sum of martingales increments as follows. For each integer $r\geq 1$, let $L_r:=L(\0,(r,r))$ and consider a partition of the two dimensional square 
$$[0,r]^2=\cup_{l=1}^{N}B_l\,$$
 into $N=r^2$ disjoint squares of size one. Let $\calF_0=\{\emptyset, \Omega\}$ and for each $k=1,\dots,N$ consider the $\sigma$-algebra
$$\calF_k=\sigma\left(\{\,(\pp,\omega_\pp)\,:\,\pp\in\cup_{l=1}^{k}B_l\cap\PP\}\right)\,,$$
and the Doob martingale
$$M_k:=\E\left(L_r\mid\calF_k\right)\,.$$
Denote by $\P_l$ the probability law induced by 
$$\{\,(\pp,\omega_\pp)\,:\,\pp\in B_l\cap\PP\}\,,$$
and by $\Omega_l$ the underlying sample space. For $\omega,\sigma\in\prod_{l=1}^{N}\Omega_l$ let 
$$[\omega,\sigma]_k:=(\omega_1,\dots,\omega_k,\sigma_{k+1},\dots,\sigma_{r^2})\in\prod_{l=1}^{N}\Omega_l\,.$$ 
Then
\begin{eqnarray}
\nonumber\Delta_k(\omega_1,\dots,\omega_k)&:=&M_k-M_{k-1}\\
\nonumber&=&\int L_r[\omega,\sigma]_{k}\prod_{l=k+1}^{r^2}d\P_l(\sigma_l)\\
\nonumber&&\qquad\qquad-\int L_r[\omega,\sigma]_{k-1}\prod_{l=k}^{N}d\P_l(\sigma_l)\\
\nonumber&=&\int L_r[\omega,\sigma]_{k}-L_r[\omega,\sigma]_{k-1}\prod_{l=k}^{N}d\P_l(\sigma_l)\,.
\end{eqnarray}
(To integrate $M_k$ over $\sigma_k$ does not change it, and it allows us to put $M_k$ and $M_{k-1}$ under the same integral.).\\

For each $1\leq k\leq N$ let
\begin{equation}
\nonumber Z_k:=\sum_{\pp\in B_k\cap\PP}w_\pp \,.
\end{equation}
Then $\{Z_k\}_{1\leq k\leq N}$ is an i.i.d. collection of random variables such that
$$\E e^{aZ_1}=\exp\left(\E e^{aw}-1\right)<\infty\,.$$
\begin{lem}\label{lem:increm}
Let $I_k$ denote the indicator function of the event that the geodesic $\varpi_r:=\varpi(\0,(r,r))$ has a point $\pp\in B_k\cap\PP$. Then
\begin{eqnarray}
\nonumber |L_r[\omega,\sigma]_{k}-L_r[\omega,\sigma]_{k-1}|&\leq&\max\left\{I_k[\omega,\sigma]_k,I_k[\omega,\sigma]_{k-1}\right\}\times\\
\nonumber&&\qquad\qquad\max\left\{ Z_k(\omega_k),Z_k(\sigma_k)\right\}\,.\\
\label{increm}
 \end{eqnarray}
 \end{lem}
\noindent{\bf Proof:}
We note that there will be no difference between $L_r[\omega,\sigma]_{k}$ and $L_r[\omega,\sigma]_{k-1}$, if no geodesic has a point $\pp\in B_k\cap\PP$ (recall that $[\omega,\sigma]_{k}$ and $[\omega,\sigma]_{k-1}$ only differ inside $B_k$), and this corresponds to the first factor in the right hand side of  \eqref{increm}. And, if one of them does intersect, then the increment can not be greater then the total weight inside the box $B_k$. (Compare with (2.12) in \cite{Ke}.)

\hfill $\Box$\\

The next step is to construct, from $\{(\pp,w_\pp):\pp\in\PP\}$, a new process $\{(\pp,\bar{w}_\pp):\pp\in\bar\PP\}$, by truncating the original model inside each box $B_k$, if $Z_k> b\log r$, in order to have
$$\sum_{\pp\in B_k\cap\bar\PP}\bar{w}_\pp \leq b\log r\,.$$
We do this truncating by putting all the weights in box $B_k$ equal to $0$ if the original weight is greater than $b\log r$. Note that the truncated process in each box is still independent of all the other boxes. Let us denote the configuration induced by the truncated model by $\bar{\omega}$ and for each random variable $Y$ that depends on $\omega$, let us write $\bar{Y}(\omega):=Y(\bar\omega)$. Thus,
$$|\bar L_r[\omega,\sigma]_{k}-\bar L_r[\omega,\sigma]_{k-1}|\leq (b\log r)\max\left\{\bar I_k[\omega,\sigma]_k,\bar I_k[\omega,\sigma]_{k-1}\right\}\,,$$
and hence
\begin{equation}\label{truncincrem2}
|\bar\Delta_{k}(\omega_1,\dots,\omega_{k})|\leq b\log r\,.
\end{equation}
The upper bound \eqref{truncincrem2} allows us to apply Lemma \ref{Kesten} to get concentration inequalities for 
$$\bar M_{N}-\bar M_0=\bar L_r(\omega)-\E \bar L_r(\omega)\,.$$
\begin{lem}\label{conditions}
Let $U_k:=2(b\log r)^2 I_k$. Then $\bar\Delta_k\leq b\log r$ and
$$\E(\bar\Delta_k^2\mid\calF_{k-1})\leq \E(\bar U_k\mid\calF_{k-1})\,.$$
\end{lem}
\noindent{\bf Proof:} By Schwarz's inequality and Lemma \ref{lem:increm},
\begin{eqnarray*}
\E(\bar\Delta_k^2\mid\calF_{k-1})&\leq &\int\int \max\left\{\bar I_k[\omega,\sigma]_k,\bar I_k[\omega,\sigma]_{k-1}\right\}\times\\
&&\qquad\max\left\{ \bar Z_k(\omega_k),\bar Z_k(\sigma_k)\right\}^2\prod_{l=k}^{N}d\bar{\P}_l(\sigma_l)d\bar{\P}_k(\omega_k)\\
&\leq&\int\int \big(\bar I_k[\omega,\sigma]_k+ \bar I_k[\omega,\sigma]_{k-1}\big)\\
&&\qquad\qquad\qquad\qquad(b\log r)^2\prod_{l=k}^{N}d\bar{\P}_l(\sigma_l)d\bar{\P}_k(\omega_k)\\
&= & \E(\bar U_k\mid\calF_{k-1})\,.
\end{eqnarray*}

 \hfill $\Box$\\

Since the number of boxes intersected the up-right path $\varpi_r$ is $2r-1$ (with probability one), we conclude that
\begin{equation}\label{eq:sumUk}
\sum_{k=1}^N\bar U_k \,=\,2(b\log r)^2\sum_{k=1}^N\bar I_k\,\leq\, 2(b\log r)^2 (2r-1)\,.
\end{equation}

By choosing $b>0$ large enough, one shows that the truncated model is a good approximation of the original model, in the sense that the probability that they will differ by $x$ goes exponentially fast to zero in $x$ (Compare with (2.30) and (2.34) in \cite{Ke}).
\begin{lem}\label{lem:trunc}
Let $b=6/a$ and $r\geq \E e^{a Z_1}/\log 2$. Then
$$\P\left(L_r(\omega)-\bar L_r(\omega)>x\right)\leq 2 \exp\left(-\frac{a}{2}x\right)\,.$$
\end{lem}

\noindent{\bf Proof:} Fix $b>0$ and a positive integer $r\geq 1$ ($N=r^2$). Notice that
\begin{equation}\label{comparing}
0\leq L_r(\omega)-\bar L_r(\omega)\leq \sum_{l=1}^{N} Z_lI\left\{Z_l>b\log r\right\}\,.
\end{equation}
By Markov's inequality,
\begin{equation*}
\P\left(\sum_{l=1}^{N}Z_l I\{Z_l> b\log r\} > x\right)\leq e^{-\frac{a}{2}x}\Big[\E\left(e^{\frac{a}{2}Z_1I\{\frac{a}{2}Z_1>\frac{a b}{2}\log r\}}\right)\Big]^{N}\,.
\end{equation*}
On the other hand,
\begin{eqnarray*}
e^{\frac{a}{2}Z_1I\{\frac{a}{2}Z_1>\frac{a b}{2}\log r\}}&=&e^{\frac{a}{2}Z_1I\{e^{aZ_1}>r^{\frac{ab}{2}}e^{\frac{aZ_1}{2}}\}}\\
&\leq& 1+ e^{\frac{a}{2}Z_1}I\{e^{aZ_1}>r^{\frac{ab}{2}}e^{\frac{aZ_1}{2}}\}\\
&\leq& 1 +\frac{e^{aZ_1}}{r^{\frac{ab}{2}}}\,,
\end{eqnarray*}
and hence,
$$ \P\left(\sum_{l=1}^{N}Z_l I\{Z_l> b\log r\} > x\right)\leq e^{-\frac{a}{2}x}\Big[1+\frac{\E e^{aZ_1}}{r^{\frac{ab}{2}}}\Big]^{N}\,.$$
Now,
\begin{eqnarray}
\nonumber \log\left(\left[1+\frac{\E e^{aZ_1}}{r^{\frac{ab}{2}}}\right]^{N}\right)&=&r^2\log\left(1+\frac{\E e^{aZ_1}}{r^{\frac{ab}{2}}}\right)\\
\nonumber&\leq&r^2\frac{\E e^{aZ_1}}{r^{\frac{ab}{2}}}\\
\nonumber&=&r^{\frac{4-ab}{2}}\E e^{aZ_1}\\
\nonumber&\leq&\log 2\,.
\end{eqnarray}
if we take $b=6/a$ and $r\geq \E e^{a Z_1}/\log 2$. Together with \eqref{comparing}, this proves Lemma \ref{lem:trunc}.

\hfill $\Box$\\

\begin{lem}\label{concentration}
Under \eqref{eq:a2}, there exist constants $b_i>0$, such that for all $r\geq b_0$
\begin{equation}\label{e-talagrand}
\P\left(|L\left(\0,(r,r)\right)-\E L\left(\0,(r,r)\right)|\geq u\right)\leq b_1\exp\left(-b_2\frac{u}{(\log r)\sqrt{r}}\right)\,,
\end{equation}
 for $u\in(0,b_3r^{3/2}\log r]$.
\end{lem}
\noindent{\bf Proof:} We have checked all the conditions of Lemma \ref{Kesten} applied to the truncated process, where we take $x_0 = 4(b\log r)^2r$, $c_0=b\log r$, $c_1=1$ and $c_2=1/4(b\log r)^2$ (this follows from Lemma \ref{conditions} and \eqref{eq:sumUk}). So there exist $c_3, c_4>0$ such that
\[ \P\left(\bar M_N - \bar M_0 \geq u\right) \leq c_3\left\{2+\frac{1}{r}\right\}\exp\left(-c_4\frac{u}{(2b\log r)\sqrt{r}}\right),\]
for all $u\leq 2(b\log r)r^{3/2}$. \\

Using Lemma \ref{lem:trunc} and the fact that $L_r\geq \bar L_r$,we can see that there exists $M>0$, such that for all $r\geq 1$, $|\E(L_r) - \E(\bar L_r)|\leq M$. Therefore, for $u\geq M$,
\begin{eqnarray*}
\P(|L_r - \E(L_r)|\geq 2u) & \leq & \P(|L_r - \E(\bar L_r)|\geq 2u - M) \\
& \leq & \P(|\bar L_r - \E(\bar L_r)|\geq u)\\
&&\qquad\qquad +\,\, \P(|L_r - \bar L_r|\geq u - M)\\
& = & \P\left(\bar M_N - \bar M_0 \geq u\right) \\
&&\qquad\qquad +\,\, \P(L_r - \bar L_r\geq u - M).
\end{eqnarray*}

Again using Lemma \ref{lem:trunc}, we can choose $b_0 = \E e^{a Z_1}/\log 2$, $b_1>0$ large enough and $b_2,b_3>0$ small enough such that \eqref{e-talagrand} holds not only for $u\in [2M,b_3(\log r)r^{3/2}]$, but also for $0\leq u\leq 2M$.

\hfill $\Box$\\

\begin{lem}\label{HoNe}
There exists a constant $b_4>0$ such that
\begin{equation}\label{eqHoNe}
\gamma r - b_4(\log^2 r)\sqrt{r}\leq \E L\left(\0,(r,r)\right) \leq \gamma r\,.
\end{equation}
\end{lem}
\noindent{\bf Proof:} We note that the right hand side of \eqref{eqHoNe} follows from the definition of $\gamma$. To prove that the left hand side of \eqref{HoNe} also holds, we parallel the arguments developed by Howard and Newmann in \cite{HoNe}.\\

We start by noting that it is enough to prove \eqref{HoNe} for integer values of $r$. Denote by $H_r$ the set of points $(x,t)$ such that $x,t\geq 0$ and $x+t=r$. Then
$$L(\0,(2r,2r))\leq \max_{\xx\in H_{2r}}L(\0,\xx) + \max_{\xx\in H_{2r}}L(\xx,(r,r))\,,$$
and hence, by symmetry with respect to $H_r$,
$$\E L\left(\0,(2r,2r)\right)\leq 2\E \max_{\xx\in H_{2r}} L(\0,\xx)\,.$$
Define, for $k\in \{0,1,\ldots, 2r\}$, $\xx_k=(k,2r-k)$, and for $k\in \{1,\ldots,2r\}$, $\zz_k=(k-1,2r-k)$. For any $\xx\in H_r$, there exists a $k\in \{1,\ldots,2r\}$ such that $\zz_k\leq \xx$. We have
\[ L\left(0,\xx\right) \leq \max\left\{L(\0,\xx_{k-1})\,,\,L(\0,\xx_k)\right) + L(\zz_k,\zz_k+(1,1)\}.\]
This implies that
\[ \E L\left(\0,(2r,2r)\right)\leq 2\E\max_{0\leq k\leq 2r}L(\0,\xx_k) + 2\E\max_{1\leq k\leq 2r}L(\zz_k,\zz_k+(1,1)).\]

The second term on the righthand side is bounded by the expectation of the maximum of the total weight in the $2r$ squares $[\zz_k,\zz_k+(1,1)]$, which is clearly bounded by $c\log r$, for some constant $c>0$ (here we use \eqref{eq:a2}). So we get
\[ \E L\left(\0,(2r,2r)\right)\leq 2\E\max_{0\leq k\leq 2r}L(\0,\xx_k) + c\log r.\]
By \eqref{eq:sym},
\begin{eqnarray*}
\max_{\xx\in H_{2r}}\E L(\0,\xx)&=&\max_{x\in[0,2r]}\E L\left(\0,(\sqrt{2r-x}\sqrt{x},\sqrt{2r-x}\sqrt{x})\right)\\
&=&\E L\left(\0,(r,r)\right)\,,
\end{eqnarray*}
and thus
\begin{eqnarray}
\nonumber\E L\left(\0,(2r,2r)\right)&\leq& 2\E \max_{0\leq k\leq 2r} L(\0,\xx_k) + c\log r\\
\nonumber&\leq& 2\E \max_{0\leq k\leq 2r} \left\{L(\0,\xx_k)-\E L(\0,\xx_k)\right\}\\
\nonumber &&\qquad + \,2\max_{0\leq k\leq 2r} \E L(\0,\xx_k) + c\log r\\
\nonumber&\leq& 2\E \max_{0\leq k\leq 2r} \left\{L(\0,\xx_k)-\E L(\0,\xx_k)\right\} \\
\nonumber&&\qquad+\,2\E L(\0,(r,r)) + c\log r\,.\\
\label{eq:sub1}
\end{eqnarray}

Let
\[ M_r = \max_{0\leq k\leq 2r} \left\{L(\0,\xx_k)-\E L(\0,\xx_k)\right\}.\]
Define for a large constant $C>0$, the event
\[ A = \bigcap_{k=0}^{2r}\{ L(\0,\xx_k)-\E L(\0,\xx_k)\leq C(\log^2 r)\sqrt{r}\}.\]
Then
\[ M_r \leq C(\log^2 r)\sqrt{r}1_A + L(\0,(2r,2r))1_{A^c}.\]
Therefore,
\[ \E M_r \leq C(\log^2 r)\sqrt{r} + \sqrt{\E \left[L(\0,(2r,2r))^2\right]\cdot \P(A^c)}.\]

We crudely bound $L(\0,(2r,2r))$ by the total weight in the square $[\0,(2r,2r)]$, and see that there exists a constant $c'>0$ such that
\[ \E \left[L(\0,(2r,2r))^2\right] \leq c'r^4.\]
We can use Lemma \ref{concentration} to conclude that
\begin{eqnarray*}
\P(A^c) & \leq & \sum_{k=0}^{2r}\P\left(L(\0,\xx_k)-\E L(\0,\xx_k)> C(\log^2 r)\sqrt{r}\right)\\
& \leq & (2r+1)b_1\exp\left(-b_2C\log r\right).
\end{eqnarray*}
By increasing $C$, this shows that there exists $c''>0$ such that for all $r\geq 1$,
\[ \E M_r\leq c''(\log^2 r)\sqrt{r}.\]
Together with \eqref{eq:sub1}, this proves that there exists $b_4>0$ such that for all $r\geq 1$
\begin{equation}\label{eq:sub2}
\E L\left(\0,(2r,2r)\right)-b_4(\log^2 r) \sqrt{r}\leq 2\E L\left(\0,(r,r)\right)\,.
\end{equation}
By Lemma 4.2 in \cite{HoNe}, \eqref{eq:sub2} implies that the left hand side of \eqref{eqHoNe} is true.

\hfill $\Box$\\

\noindent{\bf Proof of Theorem \ref{thm:fluctuation}:}
The results of Lemma \ref{concentration} and Lemma \ref{HoNe} now easily combine to Theorem \ref{thm:fluctuation}.

 \hfill $\Box$\\

\section{Geodesics and $\alpha$-Rays}
The second subject of interest to us are the geodesics. The existence of semi-infinite geodesics (or rays) for percolation like models has already been extensively study by Newman and coauthors \cite{Ne}. They developed a general approach, based on Theorem \ref{thm:fluctuation} and on the curvature of the limit shape, that leads us to what they called the $\delta$-straightness of geodesics. This property is the key for proving the existence of rays.\\

For each $\pp\in\R^2$ and $\theta\in(0,\pi/4)$, let $\Co(\pp,\theta)$ denote the cone through the axis from $\0$ to $\pp$ and of angle $\theta$. Let $R_\0^{out}(\pp)$ be the set of points $\qq\geq \pp$ such that $\pp\in\varpi(\0,\qq)$. For fixed $\delta\in(0,1)$, we say that the geodesics starting at $\0$ in $\Co((1,1),\theta)$ are $\delta$-straight if there exist constants $M,c>0$ such that for all $\pp\in \Co((1,1),\theta)$ and $\|\pp\|\geq M$,
\begin{equation}\label{eq:delta}
R_\0^{out}(\pp)\subseteq \Co(\pp,c|\pp|^{-\delta})\,.
\end{equation}

Our argument on how to control the fluctuations of the geodesics will very closely follow the proofs given in \cite{Wu} for the classical Hammersley process. In the classical case, the control of $L(\0,(r,r))$ around its asymptotic value is stronger than our Theorem \ref{thm:fluctuation}, but our result is strong enough to extend the method to the more general Hammersley process. For details of the proof, we refer to \cite{Wu}.\\

For each $L>0$, let $\partial\Cy(\pp,L)$ denote the side-edge of the truncated cylinder of width $L$, that is composed of points $\qq\in\R^2$ with $\qq\geq \pp$ and $|\qq|\leq 2|\pp|$, and such that the Euclidean distance between $\qq$ and the line through $\0$ and $\pp$ equals $L$. Assume that $\pp=(x,t)\in \Co((1,1),\theta)$, that $\qq\in R^{out}_{\0}(\pp)$ and that $\qq\in \partial\Cy(\pp,|\pp|^{1-\delta})$. Then
$$L(\0,\qq)=L(\0,\pp)+L(\pp,\qq)\,$$
or, equivalently,
$$f(\qq)-f(\qq-\pp)-f(\pp)=\Delta(\0,\pp)+\Delta(\pp,\qq)-\Delta(\0,\qq)\,,$$
where $f(\pp)=f(x,t)=\gamma\sqrt{xt}$ is the shape function and
$$\Delta(\pp,\qq)=L(\pp,\qq)-f(\qq-\pp)\,.$$
Since $\qq\in \partial\Cy(\pp,|\pp|^{1-\delta})$,
$$f(\qq)-f(\qq-\pp)-f(\pp)\geq c|\pp|^{1-2\delta}\,,$$
for a finite constant $c>0$, depending on $\theta\in (0,\pi/4)$ (here we use the curvature of the shape function; see Lemma 2.1 in \cite{Wu}). Notice that if $\delta\in(0,1/4)$ then $1-2\delta\in(1/2,1)$ and so
$$|\pp|^{1-2\delta}>> |\pp|^{1/2}\log|\pp|\,.$$
Hence, by Theorem \ref{thm:fluctuation}, for $\delta\in(0,1/4)$, we must have that if $\qq\in \partial\Cy(\pp,|\pp|^{1-\delta})$, then with very high probability $\qq\not\in R^{out}$. This can be formalized to prove the following lemma:
\begin{lem}\label{delta1}
Fix $\delta\in(0,1/4)$ and $\theta\in(0,\pi/4)$. For each $\pp=(x,t)\in \Co((1,1),\theta)$ and $\qq\in \partial \Cy(\pp,|\pp|^{1-\delta})$, let $G_\delta(\pp,\qq)$ be the event that that there exists $\pp'\in \pp+[0,1]^2$ such that $\pp'\in\varpi(\0,\qq)$. Then there exist constants $\kappa,c_i>0$ such that, if $|\pp|\geq c_0$ then
$$\P\left(G_\delta(\pp,\qq)\right)\leq c_1e^{-c_2|\pp|^{\kappa}}\,.$$
\end{lem}

We extend this Lemma to hold uniformly for $\qq$ and $\pp$ in a fixed-size finite box, then use the boxes around $\qq$ to cover the side-edge of the cylinder to get:
\begin{lem}\label{delta2}
Fix $\delta\in(0,1/4)$ and $\theta\in(0,\pi/4)$. For each $\pp=(x,t)\in \Co((1,1),\theta)$, let $G_\delta(\pp)$ be the event that there exists $\pp'\in \pp+[0,1]^2$ and $\qq\in \partial \Cy(\pp,|\pp|^{1-\delta})$ such that $\pp'\in\varpi(\0,\qq)$. Then there exist constants $\kappa,c_i>0$ such that, if $|\pp|\geq c_3$ then
$$\P\left(G_\delta(\pp)\right)\leq c_4e^{-c_5|\pp|^{\kappa}}\,.$$
\end{lem}

Now we can show $\delta$-straightness by ``gluing'' together these cylinders: if a geodesic starts close to $\pp$, with high probability it will exit the bottom edge of the cylinder $\Cy(\pp,|\pp|^{1-\delta})$. Then we cover this bottom edge with boxes $\pp_2+[0,1]^2$, where $|\pp_2|\geq 2|\pp|$, and for each of these $\pp_2$ we consider the cylinder $\Cy(\pp_2,|\pp_2|^{1-\delta})$, and so on. With Borel-Cantelli we can make the probability that a geodesic starting close to $\pp$ will ever leave through the outer edges of the boundary cylinders very small. The cylinders at the next step of the procedure have a slightly different angle than the cylinders in the previous step, but the changes in these angles are bounded by a geometric series, which means that all cylinders are contained in a cone starting at $0$, of angle $|\pp|^{-\delta}$. This is basically the same argument used for the proof of Lemma 2.4 in \cite{Wu} in the classical set-up.  This leads us to:
\begin{lem}\label{delta3}
Fix $\delta\in(0,1/4)$ and $\theta\in(0,\pi/4)$. There exist constants $\kappa,c_i>0$ such that for all $\pp\in \Co((1,1),\theta)$ with $|\pp|>c_6$, we have
\[\P\left(\left(\bigcup_{\pp'\in \pp+[0,1]^2} R_\0^{out}(\pp')\right)\subset \Co(\pp,|\pp|^{-\delta})\right)\geq 1-c_7e^{-c_8|\pp|^\kappa}\,.\]
Furthermore, with probability one, there exists $M>0$ such that for all $\pp\in \Co((1,1),\theta)$ with $|\pp|\geq M$,
\[R_\0^{out}(\pp)\subset \Co(\pp,|\pp|^{-\delta}).\]
\end{lem}

An up-right (down-left, resp.) semi-infinite path starting at $\xx$ is the lowest continuous increasing path through an ordered sequence $(\pp_n)_{n\geq 0}$ in $\R^2$, with $\pp_0=\xx$, $\pp_n\in \PP$ and $\pp_n\leq \pp_m$ ($\pp_n\geq \pp_m$, resp.) whenever $ n\leq m$. We call $(\pp_n)_{n\geq 0}$ a semi-infinite geodesic if  $\varpi(\pp_n,\pp_m) \subset (\pp_n)_{n\geq 0} $  (every part of the path is a geodesic). Finally, for each angle $\alpha\in (0,\pi/2)\cup(\pi, 3\pi/2)$, we will call a semi-infinte geodesic $(\pp_n)_{n\geq 0}$ an $\alpha$-ray if
\begin{equation}\label{eq:ray}
 \lim_{n\to \infty} \frac{\pp_n}{\|\pp_n\|} = \vec{\alpha}:=(\cos\alpha, \sin\alpha)\,.
\end{equation}
\begin{thm}\label{thm:existence}
With probability one, every semi-infinite geodesic is an $\alpha$-ray for some $\alpha\in(0,\pi/2)\cup(\pi, 3\pi/2)$, and for every $\xx\in\R^2$ and $\alpha\in(0,\pi/2)\cup(\pi, 3\pi/2)$ there exists at least one $\alpha$-ray starting at $\xx$.
\end{thm}

Uniqueness and coalescence of $\alpha$-rays do not depend upon $\delta$-straightness. The proof of these can be done by using a method  introduced in \cite{Ne}, that would work in a wide context. In \cite{Wu}, W\"uthrich applied this method to the classical Hammersley mode to get uniqueness and coalescence for fixed directions. Here we state without proof the analogous result for the Hammersley model with random weights. The reader can convince her- or himself of the validity of the theorem by checking that the proof given by W\"uthrich can be adapted mutatis mutandis to our set-up. \\

Before we state the next theorem, we shall define what we mean by convergence of paths: we say that a sequence of paths $\varpi^n$ converges to $\varpi$, and denote
$\lim_{n\to\infty}\varpi^n= \varpi$, if for all bounded subsets  $B\subset \R^2$ there exists $n_0$ such that  $\varpi^n\cap B=\varpi\cap B$ for all $n\geq n_0$.
\begin{thm}\label{thm:unicoa}
For fixed $\alpha\in(0,\pi/2)\cup(\pi, 3\pi/2)$, with probability one, for each $\xx\in\R^2$ there exists a unique $\alpha$-ray starting at $\xx$, which we denote by $\varpi_\alpha(\xx)$, and if
 $$\lim_{n\to \infty} \frac{\zz_n}{\|\zz_n\|} = \vec{\alpha}\,,$$
then
$$\lim_{n\to\infty}\varpi(\zz_0,\zz_n)= \varpi_\alpha(\zz_0)\,.$$
Furthermore, for any $\xx,\yy\in\R^2$ there exists $\cc_\alpha(\xx,\yy)\in\R^2$ such that $\varpi_\alpha(\xx)$ and $\varpi_\alpha(\yy)$ coalesce at $\cc_\alpha(\xx,\yy)$.
\end{thm}

\section{Busemann Functions}
In the middle of the fifties Busemann \cite{Bu} introduced a collection of functions to study geometrical aspects  of metric spaces. These functions are induced by a metric $d$, and by a collection of rays (semi-infinite geodesics) as follows: the Busemann function $b_\varpi(\cdot)$, with respect to a ray $(\varpi(r))_{r\geq 0}$, is the limit of 
$$b_\varpi(r,\cdot):=d(\varpi(r),\varpi(0))-d(\varpi(r),\cdot)$$
as $r$ goes to infinity. Along a ray $\varpi$ the metric $d$ becomes additive. By using the triangle inequality, this implies that the defining sequence is nondecreasing and bounded from above, and so it always converges. Using analogous considerations, one can construct Busemann functions over spaces equipped with a superadditve ``metric''  $L$ (one needs the reversed triangle inequality).\\

Using the concept of $\alpha$-rays, we will study the function $B_\alpha(\xx,\yy)$, which is defined by taking the first coalescence point $\cc=\cc(\alpha,\xx,\yy)$ between the $\alpha$-ray that starts from $\xx$ and the one that starts from $\yy$ (remember that these two rays coalesce), and setting
\begin{equation}\label{eq:defLa}
B_{\alpha}(\xx,\yy)=L(\cc,\yy)-L(\cc,\xx)\,.
\end{equation}

Note that if we take a different coalescence point $\cc'$, then $\cc\geq\cc'$ and they both lie on a geodesic. Since $L$ is additive on a geodesic, we get 
$$L(\cc',\xx) = L(\cc',\cc) + L(\cc,\xx)\,,$$
which shows that the definition of $B_\alpha(\xx,\yy)$ does not depend on the choice of the coalescence point.

Let $(\zz_n)_{n\geq 1}$ be any unbounded decreasing sequence that follows direction $(\cos\alpha,\sin\alpha)$, and let $\cc(\zz_n,\xx,\yy)$ denote the most up-right coalescence point between $\varpi\left(\zz_n,\xx\right)$ and $\varpi\left(\zz_n,\yy\right)$. By Theorem \ref{thm:unicoa}, with probability one, there exists $n_0>0$ such that
\begin{equation*}
\forall\,n\geq n_0\,\,\,\,\cc(\zz_n,\xx,\yy)=\cc(\alpha,\xx,\yy)\,
\end{equation*}
and so
\begin{equation}\label{eq:coal1}
\forall\,n\geq n_0\,\,\,\,L(\zz_n,\yy) - L(\zz_n,\xx)=B_\alpha(\xx,\yy)\,.
\end{equation}
Therefore, in geometrical terms, $B_\alpha(\xx,\cdot)$ can be seen as the \emph{Busemann function} along the ray $\varpi_\alpha(\xx)$.\\

Some properties of $B_\alpha$ are summarized in the following propositions:
\begin{prop}\label{prop:Busemann}
The distribution of the function $B_\alpha$ is translation invariant:
\begin{equation}\label{prop:Buse1}
B_{\alpha}(\cdot+\zz,\cdot+\zz)\stackrel{\cal D}{=}B_{\alpha}(\cdot,\cdot)\,.
\end{equation}
The Busemann function is anti-symmetric and additive:
\begin{equation}\label{prop:Buse2}
B_{\alpha}(\xx,\yy)=-B_{\alpha}(\yy,\xx)\ \ \mbox{and}\ \ \ B_\alpha(\xx,\zz) = B_\alpha(\xx,\yy) + B_\alpha(\yy,\zz)\,.
\end{equation}
Fix $\xx, \yy\in\R^2$ and $\pp,\qq\in \R^2$ such that $\pp,\qq\geq\0$. The function
\begin{equation}\label{prop:Buse4}
\lambda\mapsto B_\alpha(\xx+\lambda\pp,\yy+\lambda\qq)\,
\end{equation}
is c\`adl\`ag in $\lambda\in \R$.
\end{prop}

\noindent{\bf Proof:} The translation invariance of the underlying compound two-dimensional Poisson process,
$$\{\omega_\qq\,:\,\qq\in(\PP+\zz)\}\stackrel{\cal D}{=}\{\omega_\pp\,:\,\pp\in\PP\}\,,$$
implies \eqref{prop:Buse1}. Anti-symmetry follows directly from the definition of the Busemann function. Now, by taking a coalescence point  $\cc$ between $\varpi_\alpha(\xx)$, $\varpi_\alpha(\yy)$ and $\varpi_\alpha(\zz)$ we have that
\begin{eqnarray}
\nonumber B_\alpha(\xx,\zz)&=&L(\cc,\zz)-L(\cc,\xx)\\
\nonumber&=&L(\cc,\zz)-L(\cc,\yy)+L(\cc,\yy)- L(\cc,\xx)\\
\nonumber&=&B_\alpha(\yy,\zz)+B_\alpha(\xx,\yy)\,,
\end{eqnarray}
which clearly shows additivity, and finishes the proof of \eqref{prop:Buse2}.\\ 

By additivity and anti-symmetry, to prove cadlag, we can restrict our attention to  $\yy=(0,0)$, $\lambda=0$ and $\qq=(1,0)$ so that $\epsilon$ is varying on the horizontal direction close to the origin. For different values of $\yy$ the argument is similar. Choose $\cc$ as the coalescing point of $\varpi_\alpha(\0)$ and $\varpi_\alpha(\qq)$. Furthermore, define $y<0$ such that $(0,y)$ is the crossing point of  $\varpi_\alpha(\qq)$ with the $y$-axis. Then for every $\lambda \in [0,1]$, we have that
\[ B_\alpha(\0,\lambda \qq) = L(\cc,\lambda \qq) - L(\cc,\0).\]
This is because the $\alpha$-rays starting at $\lambda\qq$ are wedged in between $\varpi_\alpha(\0)$ and $\varpi_\alpha(\qq)$. Furthermore, for $\lambda$ small enough, there will be no Poisson points in the rectangle $[0,\lambda]\times [y,0]$, which means that for those $\lambda$, $L(\cc,\lambda \qq) = L(\cc,\0)$, and therefore
\[ B_\alpha(\0,\lambda \qq) = 0.\]
This proves right continuity. The existence of the left limit follows from monotonicity.

\hfill $\Box$\\

We also have the symmetries:
\begin{prop}\label{prop:sym}
For any $(x,t)\in \R^2$,
\begin{equation}\label{prop:Busesym}
B_{5\pi/2-\alpha}(\0,(x,t)) \stackrel{\cal D}{=} B_\alpha(\0,(t,x))\,
\end{equation}
and
\begin{equation}\label{prop:Busescale}
B_\alpha(\0,(x,t)) \stackrel{\cal D}{=} B_{\tilde\alpha}(\0,(r x, t/r))\,,
\end{equation}
where $\tan\tilde\alpha= r^{-2}\tan\alpha$.
\end{prop}

\noindent{\bf Proof:} Define $S$ as the reflection in the diagonal $x=t$. Then $S(\PP)$ has the same distribution as $\PP$, and $\xx\leq \yy \Leftrightarrow S(\xx)\leq S(\yy)$. This shows that for all $\xx,\yy\in\R^2$,
\[ L(\xx,\yy) \stackrel{\cal D}{=} L(S(\xx),S(\yy)).\]
Note that the lowest geodesic in the reflected case is not necessarily the reflection of the original lowest geodesic. However, we can now use \eqref{eq:coal1} to prove \eqref{prop:Busesym}. The scaling invariance \eqref{prop:Busescale} follows directly from symmetry \eqref{eq:sym}.

\hfill $\Box$\\

We now prove integrability of $B_\alpha$. Before that we require the following lemma:
\begin{lem}\label{lem:int}
Let $\yy=0$, $\alpha=5\pi/4$, and $\xx=(-1,-1)$. Let $\zz$ be the intersection point between $\varpi_{5\pi/4}(\0)$ and the one dimensional boundary of $\{\,\pp\,:\,\pp\leq\xx\,\}$. Then
$$\E L(\zz,\0)<\infty\,.$$
\end{lem}

\noindent{\bf Proof:}  Assume that $c>0$ is a small constant and that $|\zz|\leq c r$. Then
$$L(\zz,\0)\leq \min\left\{L\left((-1,-cr),\0\right),L\left((-cr,-1),\0\right)\right\}\leq P(cr)\,$$
where $P(r)$ denotes the total weight of the compound Poisson process in the set $[-1,0]\times [-cr,0]\cup [-cr,0]\times [-1,0]$. Hence
$$\P\left(L(\zz,\0)\geq r\right)\leq \P\left(|\zz|> c r\right)+\P\left(P(cr)\geq r\right)\,.$$

By choosing $c$ small enough, one can make  $\P\left(P(cr)\geq r\right)$ integrable over $r>0$. It is therefore enough to prove that $\P\left(|\zz|> c r\right)$ is integrable over $r>0$. We first consider the case where $\zz\leq (-1,-cr)$; the case $\zz\leq (-cr,-1)$ can be handled in exactly the same way. 

By Lemma \ref{delta3} ($\delta$-straightness of geodesics): fix $\delta\in(0,1/4)$ and $\theta\in(0,\pi/4)$, then for all $\pp\in \Co((1,1),\theta)$ with $|\pp|>c_6$, we have
\begin{equation}\label{eq:deltastraight}
\P\left(\left(\bigcup_{\pp'\in \pp+[0,1]^2} R_\0^{out}(\pp')\right)\subset \Co(\pp,|\pp|^{-\delta})\right)\geq 1-c_7e^{-c_8|\pp|^\kappa}\,.
\end{equation}

Let 
$$b:=(\tan(23\pi/16))^{-1},\,\,b':=(\tan(21\pi/16))^{-1}\mbox{  and }\pp_n:=(-bn,-n)\,.$$ 
Fix $n_0>0$ big enough so that ${\rm Co}(\pp_n,|\pp_n|^{-\delta})$ does not intersect $\{(-b't,-t)\,:\,t> 0\}$ for all $n\geq n_0$ (recall that $\delta\in(0,1/4)$ is fixed). Since the asymptotic direction of $\varpi_{5\pi/4}(\0)$ equals $(-1,-1)$ it will eventually be to the left of $\{(-b't,-t)\,:\,t> 0\}$. Therefore, if $\zz\leq (-1,-cr)$ and $cr\geq n_0$, there must be an $n\geq c r$ and $\pp'_n\in(\pp_n+[0,1]^2)\cap \varpi_{5\pi/4}(\0)$ such that
\[ R_\0^{out}(\pp'_n)\not \subset \Co(\pp_n,|\pp_n|^{-\delta}).\]
Using \eqref{eq:deltastraight}, we see that
$$\P\left(  \zz\leq (-1,-cr)\right)\leq \sum_{n\geq c r}c_7e^{-c_8|\pp_n|^\kappa}\,,$$
for $cr>\max\{n_0,c_6\}$, which shows integrability.

 \hfill $\Box$\\

\begin{prop}\label{prop:int}
If $\xx\leq\yy$ and $\xx\neq\yy$ then
\begin{equation}\label{eq:int}
B_\alpha(\xx,\yy)\geq 0\,\,\mbox{ and }\,\, 0<\E B_\alpha(\xx,\yy)<\infty\,.
\end{equation}
\end{prop}

\noindent{\bf Proof:} Since $L(\cc,\yy)\geq L(\cc,\xx)$ whenever $\cc\leq\xx\leq\yy$ we have that $B_\alpha(\xx,\yy)\geq 0$ whenever $\xx\leq\yy$. Note that
\[ B_\alpha(\xx,\xx+(1,1))\geq L(\xx,\xx+(1,1))\,.\]
Clearly, $\E L(\xx,\xx+(1,1))>0$ (there is a positive probability that a Poisson point will fall in between $\xx$ and $\xx+(1,1)$). Also,
\[\E B_\alpha(\xx,\xx+(1,1))=\E B_\alpha(\xx,\xx+(1,0)) + \E B_\alpha(\xx,\xx+(0,1))\,.\]
Property \eqref{prop:Busesym} then shows that
\[ \E B_\alpha(\xx,\xx+(1,0)) = \E B_\alpha(\xx,\xx+(0,1)) > 0\,.\]
Properties \eqref{prop:Buse1} and \eqref{prop:Buse2} show that for $h>0$, 
$$\E B_\alpha(\xx,(h,0))= h\E B_\alpha(\xx,\xx+(1,0))>0\,,$$ 
and the same for the vertical direction.\\

To prove that its expected value is finite, without loss of generality, assume that $\yy=0$, that $\alpha=5\pi/4$, and that $\xx=(-1,-1)$. Let $\zz$ be the intersection point between $\varpi_{5\pi/4}(\0)$ and the one dimensional boundary of $\{\,\pp\,:\,\pp\leq\xx\,\}$. By taking a coalescence point $\cc\leq \zz\leq\xx$, we have that
\begin{eqnarray}
\nonumber B_{5\pi/4}(\xx,\0)&=&L(\cc,\0)-L(\cc,\xx)\\
\nonumber &=& L(\cc,\zz)+L(\zz,\0)-L(\cc,\xx)\\
\nonumber &=& L(\zz,\0)- \{L(\cc,\xx)-L(\cc,\zz)\}\\
\nonumber&\leq&L(\zz,\0)\,.
\end{eqnarray}
By Lemma \ref{lem:int}, this proves integrability.

\hfill $\Box$\\

The most important aspect of $B_\alpha$ is a Markovian structure described below:
\begin{prop}\label{prop:BusMarkov}
For all $s\leq t$ and $x\in \R$ we have
\[ B_\alpha((0,s),(x,t)) = \sup_{z\leq x}\left\{ B_\alpha((0,s),(z,s)) + L((z,s),(x,t))\right\}.\]
\end{prop}
\noindent{\bf Proof:} Without loss of generality we can take $s=0$ (and therefore $t\geq 0$). Define $Z_\alpha=Z_\alpha(x,t)\in \R$ as the crossing-point of the $\alpha$-ray starting at $(x,t)$ with the $x$-axis. Clearly, $Z_\alpha\leq x$ and
\begin{eqnarray*}
B_\alpha(\0,(x,t)) & = & B_\alpha(\0,(Z_\alpha,0)) + B_\alpha((Z_\alpha,0),(x,t))\\
& = &  B_\alpha(\0,(Z_\alpha,0)) + L((Z_\alpha,0),(x,t)).
\end{eqnarray*}
The last equality follows from the fact that $(x,t)$ and $(Z_\alpha,0)$ are lying on an $\alpha$-ray. This means that it is enough to prove that for all $z\leq x$,
\begin{equation}\label{eq:zZalpha}
B_\alpha(\0,(z,0)) + L((z,0),(x,t))\leq B_\alpha(\0,(Z_\alpha,0)) + L((Z_\alpha,0),(x,t))\,.
\end{equation}

Suppose $\pp$ is a coalescence point of the $\alpha$-rays starting at $\0$, $(x,t)$ and $(z,0)$. Then
\[ B_\alpha(\0,(z,0)) = L(\pp,(z,0)) - L(\pp,\0)\]
and
\[  B_\alpha(\0,(Z_\alpha,0))=L(\pp,(Z_\alpha,0)) - L(\pp,\0)\,.\]
Furthermore, since $\pp, (Z_\alpha,0)$ and $(x,t)$ are elements of $\varpi_\alpha(x,t)$, we know that
\begin{eqnarray*}
L(\pp,(z,0)) + L((z,0),(x,t)) & \leq & L(\pp, (x,t))\\
& = & L(\pp,(Z_\alpha,0)) + L((Z_\alpha,0),(x,t)).
\end{eqnarray*}
From this, \eqref{eq:zZalpha} easily follows.

\hfill$\Box$\\

For $t\in\R$, define the positive measures $\nu^t_\alpha$ on $\R$, using Proposition \ref{prop:Busemann}, by
\begin{eqnarray}
\nonumber\nu^t_\alpha((x,y]) &:=& B_\alpha((x,t),(0,0))-B_\alpha((y,t),(0,0))\\
\label{eq:defnua}&=&B_\alpha((x,t),(y,t))\,,
\end{eqnarray}
where $x\leq y$. Proposition \ref{prop:BusMarkov} then shows that the process $(\nu^t_\alpha)_{t\in\R}$ is a Markov process in the space of positive and locally finite measures defined on $\R$: the future evolution of $\nu^t_\alpha$ depends on the Poisson process in the upper-half plane $\R\times (t,\infty)$ and on the present value of $\nu^t_\alpha$, not on the past of the process (which is of course independent of the Poisson process in $\R\times (t,\infty)$). By translation invariance (Proposition \ref{prop:Busemann}) the distribution of $\nu^t_\alpha$ does not depend on $t$:
\begin{equation}\label{eq:equinua}
\nu^t_\alpha\stackrel{\cal D}{=}\nu^0_\alpha\,.
\end{equation}
Thus, this distribution is an equilibrium (or time invariant) measure for the underlying  Markov process. In the next section we will describe the generator of this Markov process, which will be an extension of the classical Hammersley interacting particle process.

\newpage

\chapter{The Hammersley Interacting Fluid System}\label{ch:fluid}

\section{The Graphical Construction} 
It is well known that the classical Hammersley model, where all weights are 1, described in \cite{AlDi}, has a representation as an interacting particle system. The Hammersley Last Passage Percolation model with random weights has a similar description, although a better name might be an interacting fluid system.\\

We start by restricting the compound Poisson process $\{\omega_\pp\,:\,\pp\in\PP\}$ to $\R\times\R_+$. Then we choose a positive, locally finite measure $\nu$ defined on $\R$. Usually, these measures will be purely atomic, but this is not necessary. To each measure $\nu$ we associate a non-decreasing process $\nu(\cdot)$ defined by
\begin{equation}\label{eq:defnuprocess}
\nu(x)=\left\{\begin{array}{ll} \nu((0,x]) & \mbox{for } x\geq 0\\
-\nu((x,0]) & \mbox{for } x<0.\end{array}\right.
\end{equation}
Note that $\nu(\cdot)$ is a cadlag function. Although the details are a bit cumbersome, all the results we will show can be extended, mutatis mutandis, to the case where $\nu(x)=-\infty$ for $x<0$, which would correspond to a non-locally finite measure with an infinite fluid density to the left of $0$. This is a quite natural starting condition, but we will not use it explicitly.\\

The Hammersley interacting fluid system $(M_{\nu}^t)_{t\geq 0}$ is a stochastic process with values in the space of positive, locally finite measures on $\R$. Its evolution is defined as follows: if there is a Poisson point with weight $\omega$ at a point $(x_0,t)$, then 
$$M_\nu^{t}(\{x_0\}) = M_\nu^{t-}(\{x_0\}) + \omega\,$$
and for $x> x_0$,
\begin{equation}\label{eq:evolM}
M_\nu^t((x_0,x]) = (M_\nu^{t-}((x_0,x]) - \omega)_+\,.
\end{equation}

\begin{figure}[!ht]
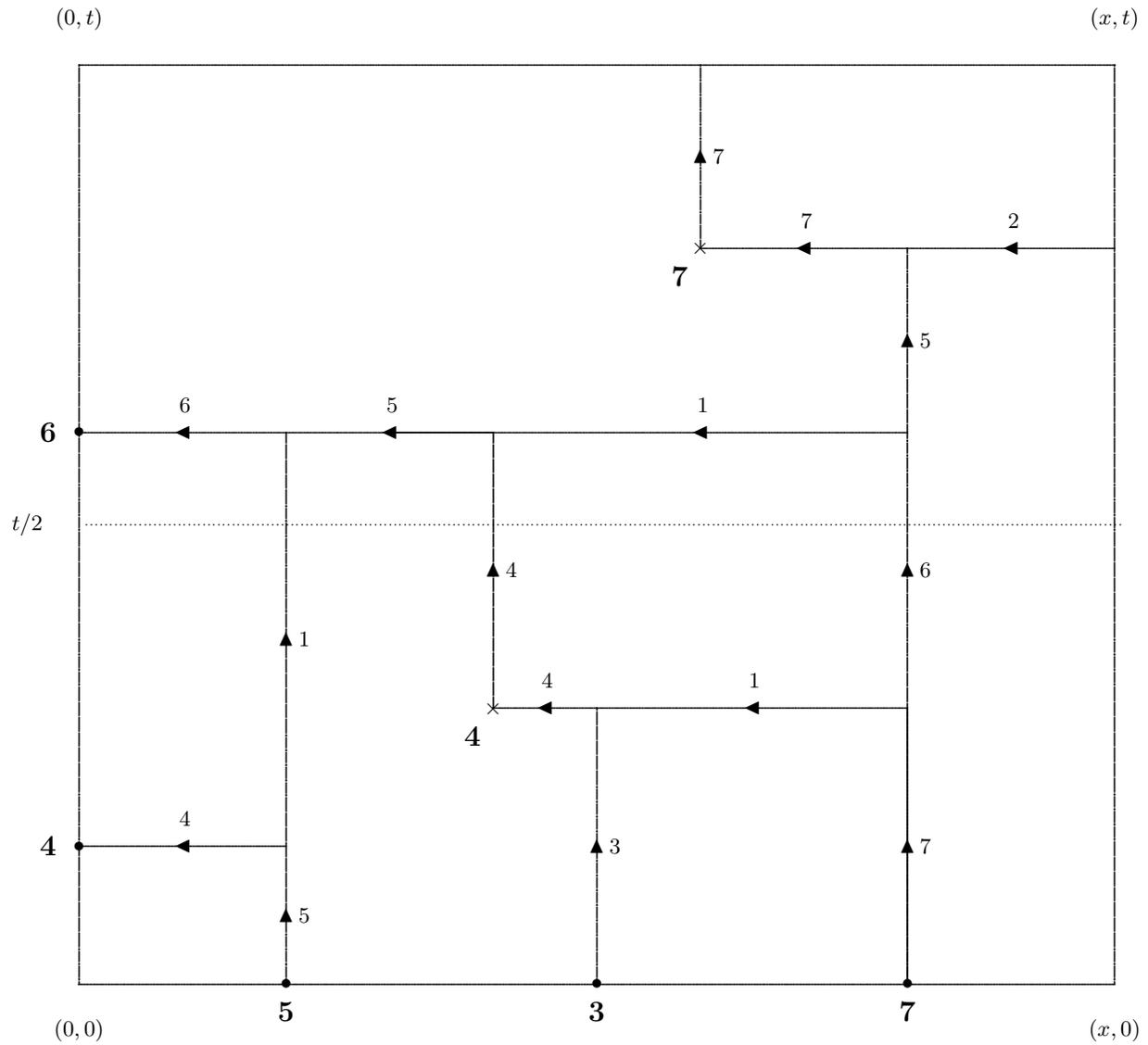

\begin{center}
\strut
$$
\beginpicture
\footnotesize
\setcoordinatesystem units <0.09\textwidth,0.08\textwidth>
  \setplotarea x from 0 to 10, y from 0 to 12
\multiput {\bf $\times$} at
4 3
6 8
/

\multiput {$\bullet$} at
2 0
5 0
8 0
0 1.5
0 6
/
\put {$(0,0)$} at 0.0 -0.5
\put {$(x,0)$} at 10 -0.5
\put {$(0,t)$} at 0.0 10.5
\put {$(x,t)$} at 10 10.5
\put {$t/2$} at -0.5 5

\setlinear
\plot 0 10  10 10  10 0 0 0 0 10
/
\plot 10 8 8 8
/
\plot 8 0 8 8 6 8 6 10
/
\plot 8 0 8 3 5 3
/
\plot 5 0 5 3 4 3 4 6 0 6
/
\plot 8 6 3 6
/
\plot 2 0 2 6
/
\plot 2 1.5 0 1.5
/
\multiput {$\blacktriangle$} at
2 0.75
5 1.5
8 1.5
8 4.5
8 7
6 9
4 4.5
2 3.75
/
\multiput {$\blacktriangleleft$} at
1 1.5
1 6
3 6
6 6
6.5 3
7 8
9 8
4.5 3
/

{\large
\put {\bf 5} at 2 -0.3
\put {\bf 3} at 5 -0.3
\put {\bf 7} at 8 -0.3
\put {\bf 4} at -0.3 1.5
\put {\bf 6} at -0.3 6
\put {\bf 4} at 3.8 2.7
\put {\bf 7} at 5.8 7.7
}

\put {$4$} at 4.1 4.5
\put {$4$} at 4.45 3.3
\put {$3$} at 5.1 1.5
\put {$1$} at 6.45 3.3
\put {$7$} at 8.1 1.5
\put {$6$} at 8.1 4.5
\put {$5$} at 8.1 7
\put {$2$} at 8.95 8.3
\put {$7$} at 6.95 8.3
\put {$1$} at 5.95 6.3
\put {$5$} at 2.95 6.3
\put {$6$} at 0.95 6.3
\put {$1$} at 2.1 3.75
\put {$4$} at 0.95 1.8
\put {$5$} at 2.1 .75
\put {$7$} at 6.1 9

\setdots<2pt>

\plot 0 5 10 5
/

\endpicture
$$
\caption{In this picture,  restricted to $[0,x]\times[0,t]$, the measure $\nu$ consists of three atoms of weight $5$, $3$ and $7$. The compound Poisson process consists of two points of weight $4$ and $7$. The measure $M_\nu^{t/2}$ consists of three atoms of weight $1$, $4$ and $6$, while at time $t$, it consists of one atom with weight $7$. A total weight of $4+6$ has left the box while a total weight of $2$ has entered.}\label{fig:hamprocess}
\end{center}
\end{figure}

Here, $M_\nu^{t-}$ is the ``mass distribution'' of the fluid at time $t$ if the Poisson point at $(x_0,t)$ would be removed. To the left of $x_0$ the measure does not change. In words, the Poisson point at $(x_0,t)$ moves a total mass $\omega$ to the left, to the point $x_0$, taking the mass from the first available fluid to the right of $x_0$. See Figure \ref{fig:hamprocess} for a visualization in case of atomic measures, of the process inside a space-time box.\\

It is not true that the evolution $M^t$ is well defined for all measures $\nu$ (e.g. if we start with a finite number of particles to the left of $0$, every particle would be pulled instantaneously to $-\infty$). We will follow the Aldous and Diaconis \cite{AlDi} graphical representation in the Last Passage Percolation model (compare to the result in the classical case, found in their paper):
\begin{prop}\label{prop:construction}
Let $\calN$ be the set of all positive, locally finite measures $\nu$ such that
\begin{equation}\label{eq:leftden}
\liminf_{y\to-\infty}\frac{\nu(y)}{y}>0\,.
\end{equation}
For each $\nu\in\calN$, the process defined by
\begin{equation}\label{eq:lastpassage}
L_{\nu}(x,t):=\sup_{z\leq x} \left\{ \nu(z) + L((z,0),(x,t))\right\}\ \ \ \ \ (x\in\R,\,\,t\geq 0)\,
\end{equation}
is well defined and the measure
\begin{equation}\label{eq:connection}
M^t_\nu((x,y]):=L_\nu(y,t)-L_\nu(x,t)\,,
\end{equation}
evolves according to the Hammersley interacting fluid system.
\end{prop}
\noindent{\bf Proof:} It follows immediately from the definition that $L_\nu(x,t)$ is increasing in $x$ and $t$, even if $L_\nu$ would not be finite everywhere. This implies that if we can prove that $L_\nu$ is finite on for example $\ZZ\times \ZZ_+$ with probability one, then almost surely, $L_\nu$ is finite everywhere. Therefore, we need only prove for any fixed point $(x,t)$ that $L_\nu(x,t)$ is finite with probability one. \\

Use \eqref{eq:sym} to see that for any $\eps>0$ and $r=\sqrt{|x-z|t}$ we have that
\[ \P\left(  L\left((z,0),(x,t)\right) > -\eps z\right) = \P\left( L\left((0,0),r,r)\right) > -\eps z\right).\]
Using Theorem \ref{thm:fluctuation}, we get
\begin{eqnarray*}
\P\left(  L\left((z,0),(x,t)\right) > -\eps z\right)&\leq &\P\left( L\left((0,0),(r,r)\right) > \eps (-x+ r^2/t)\right)\\
&\leq& ce^{-c'r}\,,
\end{eqnarray*}
for some positive constants $c$ and $c'$, and for large enough $r$. Borel-Cantelli then shows that
$$\frac{L\left((z,0),(x,t)\right)}{-z}\stackrel{\rm a.s.}\,,{\longrightarrow} 0\,,$$
as $z\to-\infty$. Now \eqref{eq:leftden} gives us \eqref{eq:lastpassage}.\\

We can show that $x\mapsto L_\nu((z,0),(x,t))$ is c\`adl\`ag using the fact that for $y\geq x$,
\[ L((z,0),(y,t)) \leq L((z,0),(x,t)) + L((x+,0),(y,t))\]
and
\[ \lim_{y\downarrow x} L((x+,0),(y,t))=0\,,\]
where $x+$ means that you are not allowed to use a possible Poisson point directly above $x$. Therefore, $M_\nu^t$ is indeed a locally finite measure on $\R$.\\

To see that $M_\nu^t$ follows the Hammersley interacting fluid dynamics that we have just defined, suppose that there is a Poisson point in $(x_0,t)$ with weight $\omega$. For $x<x_0$, this Poisson point has no effect, so $M_\nu^t = M_\nu^{t-}$ on $(-\infty,x_0)$. Clearly,
\[ M_\nu^{t}(\{x_0\}) = M_\nu^{t-}(\{x_0\}) + \omega\,.\]
If $x> x_0$, then the longest path to $(x,t)$ that attains the supremum in \eqref{eq:lastpassage} can either use the weight in $(x_0,t)$, which would give 
$$L_\nu(x,t)=L_\nu(x_0,t)=L_\nu(x_0,t-)+\omega\,,$$
or it could ignore the weight in $(x_0,t)$, which would give 
$$L_\nu(x,t)=L_\nu(x,t-)\,.$$ 
This proves that
\begin{eqnarray*}
M_\nu^{t}((x_0,x]) & = & L_\nu(x,t) - L_\nu(x_0,t)\\
& = & \max\left\{L_\nu(x_0,t-)+\omega\,,L_\nu(x,t-)\right\} - L_\nu(x_0,t-) - \omega\\
& = & \left(M^{t-}_\nu((x_0,x])-\omega\right)_+.
\end{eqnarray*}
\hfill$\Box$\\

\section{The Generator}
Recall that we have defined the Hammersley interacting fluid system $M^t$ as a Markov process on the space $\cal N$ (Proposition \ref{prop:construction}).
To compute the generator $G$ of $M^t$, we endow the space $\cal N$ with the weak topology of convergence of continuous functions with compact support, so $\nu_n\to \nu$ iff for all $\phi\in C_c(\R)$, we have
\[ \int_{\R} \phi(x)\nu_n(dx) \to \int_{\R} \phi(x)\nu(dx).\]

Now assume that we have a continuous bounded function $f:\calN\mapsto\R$ that only depend on``cylinders'', so on the configuration of the fluid in a compact interval. Consider the transition operator
$$P_t f(\nu):=\E\left(f(M^t)\mid M^0=\nu\right)\,.$$
\begin{prop}\label{prop:generator} Let $\cal R_{z,\omega}$ represent the  operator acting on $\cal N$ that moves a total mass $\omega$ to the left, to point $z$, taking mass from the first available fluid to the right of $z$. Let $f:\calN\mapsto\R$ be a continuous bounded cylinder function. Then
$$\lim_{t\to 0}\frac{P_t f(\nu)-f(\nu)}{t}=Gf(\nu)\,,$$
where
$$Gf(\nu) = \int_{-\infty}^\infty  \int_0^\infty  \left[f(R_{z,\omega}(\nu)) -
f(\nu) \right]\,dF(\omega) \,dz\,.$$
\end{prop}

\noindent{\bf Proof:} Suppose $f$ only depends on the configuration $\nu$ restricted to the interval $[a,b]$. We wish to determine $P_tf(\nu)$ for small values of $t$. Choose $N$ large with $N\geq -a$, and define $X_1(t)\geq X_2(t)\geq\ldots$ as the $x$-coordinates of the rightmost Poisson point in the strip $(-\infty,-N]\times [0,t]$, the second rightmost point, etc. Define $\omega_i(t)$ as the weight of the Poisson point with $x$-coordinate $X_i(t)$. The configuration of $M^t$ at $[a,b]$ will not depend on the compound Poisson process in the strip $(-\infty,a-N]\times [0,t]$, if we have the event
\[ E_t=\left\{\sum_{j=1}^i\omega_j(t)\leq \nu([X_i(t),a))\,\,\,\forall i\geq 1\right\}.\]

Define the following random walk in $\R_+^2$:
\[ S_n(t) = \left[-X_i(t) + a - N, \sum_{i=1}^n \omega_i(t)\right].\]
The steps to the right are iid with distribution ${\rm Exp}(t)$, and the steps up are iid with the distribution $F$ of the weights. We define the corresponding step-function for $x\in [0,\infty)$:
\[ \chi^{(t)}(x) = \sum_{i=1}^\infty 1_{\{x\geq -X_i(t) + a - N\}} \omega_i.\]
Clearly, $\{\chi^{(t)}\}$ is a decreasing family (as $t\to 0$) of increasing functions. In fact, for fixed $t$, $\chi^{(t)}$ describes a compound Poisson process, with intensity $t$ and jump distribution $F$. Define
\[ \psi(x) = \nu([a-N-x,a)).\]
Since $\nu\in{\cal N}$, we know that $\nu$ satisfies \eqref{eq:leftden}, so there exists $\eps>0$ and $N_0\geq 1$ such that for all $x\geq 0$ and all $N\geq N_0$,
\[ \psi(x)\geq \eps x + \eps N.\]

The event $E_t$ occurs when $\chi^{(t)}(x)\leq \psi(x)$ for all $x\geq 0$. Since the distribution of the weights $F$ has an exponential moment, we know that there exist constants $c_1,c_2>0$ (depending only on $\eps$ and $N_0$) such that for all $t$ small enough and $N\geq N_0$,
\[ \P(\forall\ x\geq 0\ :\ \chi^{(t)}(x) < \eps x + \eps N) \geq 1 - c_1te^{-c_2N}.\]
This follows from the classical result where the weights (i.e., the jumps up) have an exponential distribution. We conclude that
\[ \P(E_t^c) \leq c_1te^{-c_2N}.\]

Now choose $N_t=-\log(t)$. Also, define the event
\[ A_t = \{ \mbox{the strip } [a-N_t,b]\times [0,t] \mbox{ contains }0\mbox{ or }1\mbox{ Poisson point}\}.\]
Then, given that $M^0=\nu$, we get
\begin{eqnarray*}
\E(f(M^t) - f(\nu)) & = & \E(1_{E_t\cap A_t}
(f(M^t)-f(\nu)))\\
&&\qquad\qquad+\, \E(1_{E^c_t\cup A_t^c}(f(M^t)-f(\nu)))\\
& = & \E(f(M^t)-f(\nu)\,|\,E_t\cap A_t)(1-\P(E_t^c\cup A_t^c)) \\
&&\qquad\qquad\qquad+\,\, \E(1_{E^c_t\cup A_t^c}(f(M^t)-f(\nu)))\\
& = & \frac{t(b-a+N_t)}{1+t(b-a+N_t)}\\
&\times&\int_0^\infty\,\frac{1}{b-a+N_t}\, \int_{a-N_t}^b \big[f(R_{z,\omega}(\nu))\\
&& \qquad\qquad\qquad\qquad\qquad- \,f(\nu)\big]\,dz\,dF(\omega)\\
&&\qquad\qquad\qquad\qquad\qquad\qquad\qquad\qquad+\,\, o(t)\,.
\end{eqnarray*}

The $o(t)$-term follows from the fact that $\P(E_t^c)=o(t)$, $\P(A_t^c)=o(t)$ and $f$ is bounded. The factor in front of the integrals equals the probability of having 1 Poisson point in the strip, given that the number of points equals $0$ or $1$. Of course the $x$-coordinate of this one point will have a uniform distribution on $[a-N_t,b]$. The proposition now follows when we divide by $t$, take the limit $t\to 0$, and note that Poisson points to the right of $b$ have no influence on $f(M^{(t)})$, so
\[ \int_{a-N_t}^b f(R_{z,\omega}(\nu)) - f(\nu)\,dz = \int_{a-N_t}^\infty f(R_{z,\omega}(\nu)) - f(\nu)\,dz.\]
\hfill$\Box$\\

Assume that we have a probability measure $\Psi$ defined on $\calN$ and consider $\nu\in \calN$ as a realization of this probability measure. We say that $\nu$ is time invariant for the Hammersley interacting fluid process (in law) if
$$M_\nu^t\stackrel{\cal D}{=}M_\nu^0=\nu\,\,\,\mbox{ for all }\,\,\,t\geq 0\,.$$
In this case, we also say that the underlying probability measure on $\calN$ is an equilibrium measure. In particular, a probability measure $\Psi$ on $\calN$ is an equilibrium measure if it satisfies
\begin{equation}\label{eq:generator}
\int Gf(\nu)d\Psi (\nu) = 0\,.
\end{equation}

It was this interplay between the longest path description and the equilibrium interacting particle system that proved very fruitful in the results for the classical Hammersley process in \cite{CaGr1,CaGr2}. We will attempt the same in the interacting fluid system, but since the equilibrium solution to \eqref{eq:generator} is not explicitly known (in most cases), we needed to develop new tools and ideas, which in fact also had interesting applications for the classical case.

\section{Exit Points}
We have two remarks concerning the interactive fluid process. The first remark concerns the flux of the fluid system through the $t$-axis. We define the flux measure $\nu^*$ on $[0,\infty)$ such that for $t\geq s\geq 0$,
\begin{equation}\label{eq:defflux}
\nu^*((s,t]) = L_\nu(0,t) - L_\nu(0,s).
\end{equation}
In Figure \ref{fig:hamprocess} we can see that $\nu^*((0,t]) = 10$, whereas $\nu^*((0,t/2]) = 4$. For $x,t\geq 0$, we also have the equality
\begin{eqnarray}
\nonumber L_\nu(x,t) &=& \max\Big\{ \sup_{0\leq z\leq x}\{\nu(z) + L((z,0),(x,t))\}\,,\\ 
\nonumber&&\qquad\qquad\qquad\qquad\sup_{0\leq z\leq t}\{\nu^*(z) + L((0,z),(x,t))\}\Big\}\,.\\
\label{eq:fluxformula}
\end{eqnarray}

The way to prove \eqref{eq:fluxformula} leads us to our second remark: the supremum in the definition of $L_\nu(x,t)$ is actually attained (it is therefore a maximum). This is the statement of the proposition below.
\begin{prop}\label{prop:Znuexists}
With probability 1, the set $\{ z\leq x\ :\ L_\nu(x,t) = \nu(z) + L((z,0),(x,t))\}$ is non-empty, for all $(x,t)\in \R\times [0,\infty)$. Furthermore, if we define
\begin{equation}\label{eq:defZnu}
Z_\nu(x,t) = \sup \{ z\leq x\ :\ L_\nu(x,t) = \nu(z) + L((z,0),(x,t))\},
\end{equation}
then
\[ L_\nu(x,t) = \nu(Z_\nu(x,t)) + L((Z_\nu(x,t),0),(x,t)).\]
Also, the following local comparison property holds for $x<y$ and $t\geq 0$: if  $Z_\nu(y,t)\leq 0$ then
\[ L_\nu(y,t)-L_\nu(x,t)\leq L(\0,(y,t)) - L(\0,(x,t))\,;\]
If $Z_\nu(x,t)\geq 0$ then
\[ L(\0,(y,t)) - L(\0,(x,t))\leq L_\nu(y,t)-L_\nu(x,t)\,.\]
\end{prop}
\noindent{\bf Proof:} The first part of this proposition relies on the fact that on a compact interval, the sum of a non-decreasing right-continuous function ($z\mapsto\nu(z)$) and a non-increasing left-continuous function ($z\mapsto L((z,0),(x,t))$) attains its maximum, and the set of maxima is closed. \\

Now we have to show that with probability 1, for all $(y,s)\in \R\times [0,\infty)$, the supremum over $z\leq y$ can actually be restricted to a compact set. For a single $(x,t)$ (and therefore for a countable set of $(x,t)$'s), this has been done in the proof of Proposition \ref{prop:construction}, using \eqref{eq:leftden} and Theorem \ref{thm:shape} for $L$. Then we can conclude the desired compactness property for all $(y,s)$ by using the inequality
\[ \nu(z) + L((z,0),(y,s))\leq \nu(Z_\nu(x,t)) + L((Z_\nu(x,t),0),(y,s))\] 
for all $z\leq Z_\nu(x,t)$, $y\geq x$ and $s\leq t$.

This last inequality can be seen as follows: suppose there exists $z<Z_\nu(x,t)$ such that
\[ \nu(z) + L((z,0),(y,s))> \nu(Z_\nu(x,t)) + L((Z_\nu(x,t),0),(y,s)).\]
Then the geodesic $\varpi((Z_\nu(x,t),0),(x,t))$ must intersect the geodesic $\varpi((z,0),(y,s))$ in some point $\cc$. It follows that 
$$L((z,0),(y,s)) = L((z,0),\cc) + L(\cc,(y,s))\,,$$ 
and 
$$L((Z_\nu(x,t),0),(x,t)) = L((Z_\nu(x,t),0),\cc) + L(\cc,(x,t))\,.$$ 
Furthermore, we have that
\[ \nu(z) + L((z,0),(y,s))> \nu(Z_\nu(x,t)) + L((Z_\nu(x,t),0),\cc) + L(\cc,(y,s)).\]
Combining this gives
\begin{eqnarray*}
\nu(z) + L((z,0),(x,t))& \geq &  \nu(z) + L((z,0),\cc) + L(\cc,(x,t)) \\
& > & \nu(Z_\nu(x,t)) + L((Z_\nu(x,t),0),\cc)\\
&& \qquad\qquad\qquad\qquad\qquad +\,L(\cc,(x,t))\\
 & = & \nu(Z_\nu(x,t)) + L((Z_\nu(x,t),0),(x,t))\,,
\end{eqnarray*}
which contradicts the definition of $Z_\nu(x,t)$.\\

Now we prove the second statement. If $Z_\nu(y,t)\leq 0$ then there exists $z_y\leq 0$ such that
\[ L_\nu(y,t) = L((z_y,0),(y,t)) + \nu(z_y).\]
Define $\zz$ as the (an) intersection between the geodesic from $(z_y,0)$ to $(y,t)$ and the geodesic from $\0$ to $(x,t)$. Then we can see that
\begin{eqnarray*}
L_\nu(x,t) &\geq & L((z_y,0),(x,t)) + \nu(z_y)\\
& \geq & L((z_y,0),\zz) + L(\zz,(x,t)) + \nu(z_y).
\end{eqnarray*}
This implies that
\begin{eqnarray*}
L_\nu(y,t) - L_\nu(x,t) & \leq & L((z_y,0),(y,t)) - L((z_y,0),\zz) - L(\zz,(x,t))\\
& = & L(\zz,(y,t)) - L(\zz,(x,t))\\
& = & L(\zz,(y,t)) - L(\0,(x,t)) + L(\0,\zz)\\
& \leq & L(\0,(y,t)) - L(\0,(x,t)).
\end{eqnarray*}
The proof of the third statement follows the same line. 

\hfill$\Box$\\

Note that we have also proved the following statement:
\begin{equation}\label{eq:orderingZnu}
\forall\ y\geq x\ \forall\ 0\leq s\leq t\ :\ Z_\nu(y,s)\geq Z_\nu(x,t).
\end{equation}
Equation \eqref{eq:fluxformula} now can be seen as follows: define $\nu_-$ as the restriction of $\nu$ to $(-\infty,0]$ (so $\nu_-$ has no mass on the positive $x$-axis). Then for $x,t\geq 0$, there exists $Z\leq 0$, such that
\[ L_{\nu_-}(x,t) = \nu(Z) + L((Z,0),(x,t))\,.\]
The finite geodesic $\varpi((Z,0),(x,t))$ crosses the positive $t$-axis in some point, call it $(0,s)$. It is not hard to see that 
$$\nu^*(s) = \nu(Z) + L((Z,0),(0,s))\,$$ 
and
\[ L_{\nu_-}(x,t) = \nu^*(s) + L((0,s),(x,t))\geq \nu^*(\tilde s) + L((0,\tilde s),(x,t))\]
for all $0\leq \tilde s\leq t$. This follows from similar arguments as at the end of the proof of Proposition \ref{prop:Znuexists}. \\

Finally, we remark that
\begin{eqnarray*}
L_\nu(x,t) &=& \max\Big\{ \sup_{z\leq 0} \{\nu_-(z) + L((z,0),(x,t))\}\,, \\
&&\qquad\qquad\qquad\qquad\sup_{z>0}\{ \nu(z) + L((z,0),(x,t))\}\Big\}\,,
\end{eqnarray*}
from which Equation \eqref{eq:fluxformula} follows.

\section{Law of Large Numbers for Exit Points}
Define the following crossing points:
\begin{equation}\label{def:Zalpha}
(Z_\alpha(x,t),0) \mbox{ is the crossing point of }\varpi_\alpha(x,t)\ \mbox{and } \R\times\{0\}.
\end{equation}
This was already used in the proof of Proposition \ref{prop:BusMarkov}. Analogously,
\begin{equation}\label{def:Zalpha*}
(0,Z^*_\alpha(x,t)) \mbox{ is the crossing point of }\varpi_\alpha(x,t)\ \mbox{and } \{0\}\times\R.
\end{equation}
Recall \eqref{eq:defnua} and let $\nu_\alpha:=\nu_\alpha^0$. We have the following proposition.
\begin{prop}\label{prop:exitpoint}
For all $x\in \R$ and $t\geq 0$, we have
\[ Z_\alpha(x,t) = Z_{\nu_\alpha}(x,t).\]
\end{prop}
\noindent{\bf Proof:} From \eqref{eq:zZalpha} we immediately get that $Z_{\nu_\alpha}(x,t)\geq Z_\alpha(x,t)$. Now suppose $Z_{\nu_\alpha}(x,t)> Z_\alpha(x,t)$. Since
\begin{eqnarray*}
B_\alpha(\0,(x,t)) & =& L_{\nu_\alpha}(x,t) \\
& = & \nu_\alpha(Z_{\nu_\alpha}(x,t)) + L((Z_{\nu_\alpha}(x,t),0),(x,t))\\
& = & B_\alpha(\0,(Z_{\nu_\alpha}(x,t),0)) + L((Z_{\nu_\alpha}(x,t),0),(x,t))\,,
\end{eqnarray*}
we conclude that
\[  L((Z_{\nu_\alpha}(x,t),0),(x,t)) = B_\alpha((Z_{\nu_\alpha}(x,t),0),(x,t)).\]

Now define $\cc$ as the coalescing point of the $\alpha$-ray starting at $(Z_{\nu_\alpha}(x,t),0)$ and the $\alpha$-ray starting at $(x,t)$ (which goes through $(Z_\alpha(x,t),0)$). Since the two points are on one $\alpha$-ray, then
\begin{eqnarray*}
L(\cc,(x,t)) & = & B_\alpha(\cc,(x,t))\\
& = & B_\alpha(\cc,(Z_{\nu_\alpha}(x,t),0)) + B_\alpha((Z_{\nu_\alpha}(x,t),0),(x,t))\\
& = & L(\cc,(Z_{\nu_\alpha}(x,t),0)) + L((Z_{\nu_\alpha}(x,t),0),(x,t))\,,
\end{eqnarray*}
which would imply that there exists a longest path from $\cc$ to $(x,t)$ which is strictly below the original $\alpha$-ray, contradicting the uniqueness of the lowest finite geodesic.

\hfill $\Box$\\

Since the $\alpha$-ray $\varpi_\alpha(\0)$ has a.s. the asymptotic angle $\alpha$, by translation invariance, we can easily see that
\begin{equation}\label{eq:vanishexitpoint}
\lim_{t\to\infty}\P\big(|Z_{\alpha}(t,(\tan\alpha) t)|\geq\epsilon t\big)=\lim_{t\to\infty}\P\big(|Z^*_{\alpha}(t,(\tan\alpha) t)|\geq\epsilon t\big)=0\,.
\end{equation}

However, we can also control $Z_\nu$ for more general $\nu$. Compare the following theorem to Lemma 3.3 in \cite{FeMaPi}.
\begin{thm}\label{thm:genexitcontrol}
Suppose $\nu\in \calN$. Assume that
\begin{equation}\label{eq:genexitcontrol}
\liminf_{z\to-\infty}\frac{\nu((z,0])}{-z}\geq\frac{\gamma}{2}\sqrt{\tan\alpha}\,\,\mbox{ and }\,\,\limsup_{z\to\infty}\frac{\nu((0,z])}{z}\leq\frac{\gamma}{2}\sqrt{\tan\alpha}\,.
\end{equation}
Let $h\in \R$. Then, with probability one,
$$\lim_{t\to\infty}\frac{Z_{\nu}\big(t + h, (\tan\alpha)t\big)}{t}=0\,.$$
Furthermore, define
\begin{equation}\label{eq:deftildeZ}
\tilde{Z}_{\nu,h}\big(t,(\tan\alpha)t\big) = \arg\sup_{z\leq t}\left\{ \nu(z) + L(-\bft_\alpha+(z,0),(h,0))\right\}.
\end{equation}
where $\bft_\alpha:=(t,(\tan\alpha)t)$. Here we take the right-most maximum. The idea is that we move the origin to $-\bft_\alpha$ and look at the exit point for $(h,0)$. Then, with probability one,
$$\lim_{t\to\infty}\frac{\tilde{Z}_{\nu,h}\big(t,(\tan\alpha) t\big)}{t}=0\,.$$
\end{thm}
\noindent {\bf Proof:} We start with the proof for $Z_\nu$. Using the transformation \eqref{eq:sym}, we can assume, without loss of generality, that $\alpha=5\pi/4$.\\

The proof of this lemma is based on the following elementary estimate for the shape function $f(x,t)=\gamma\sqrt{xt}$: for $s\geq 0$
\begin{equation}\label{eq:curvature}
f(t+s,t)-f(t,t)\leq \left\{ \begin{array}{cl}
\frac{\gamma}{2}s - \frac{\gamma}{32}\frac{s^2}{t} & \mbox{if } s\leq 8t,\\
\vspace{0.2cm}
\frac{\gamma}{\sqrt{8}} s & \mbox{if } s\geq 8t.
\end{array}\right.
\end{equation}
Fix $\eps>0$ and suppose that $Z_{\nu,h}(t):=Z_{\nu}(t+h,t)\leq -\eps t-1$. Define $k=-\lceil Z_\nu(t)\rceil$ and $n=\lfloor t+h\rfloor$. It follows that
$$\nu(-k+1)+L\big((-k,0),(n+1,n+1)\big)-L\big((0,0),(n,n)\big)\geq 0\,.$$

Choose $\eta>0$ small enough (we will see how small). For $t$ big enough, we know from \eqref{eq:genexitcontrol} that
\[ \nu(-k+1)\leq -\frac12\gamma k + \eta k.\]
Here we use that $k\geq \eps t$. This implies that for $t$ big enough,
\begin{equation}\label{eq:bendL}
L\big((-k,0),(n+1,n+1)\big)-L\big((0,0),(n,n)\big)\geq \frac12\gamma k - \eta k\,.
\end{equation}

A straightforward application of Theorem \ref{thm:shape} and Borel-Cantelli shows that for $t$ big enough (and therefore for $n$ big enough), we will have that
\[ |L\big((0,0),(n,n)\big) - f(n,n)|\leq \eta n.\]
Below we will work out a more complicated application of Theorem \ref{thm:shape}. We now have for $t$ big enough
$$L\big((-k,0),(n+1,n+1)\big)-f(n,n)\geq \frac12\gamma k - \eta k -\eta n\,.$$
For $t$ big enough, we will have that $k\geq \frac12 \eps n$ (since $h$ is fixed). Now we subtract $f(n+k,n)-f(n,n)$ from both sides of the inequality and use \eqref{eq:curvature} to see that
$$L\big((-k,0),(n+1,n+1)\big)-f(n+k,n)$$
is greater or equal to
$${\left\{ \begin{array}{ll}
\frac{\gamma}{32}\frac{k^2}{n} - \eta k - \eta n & \mbox{if } \frac12 \eps n\leq k\leq 8n\,,\\
\\
\frac{\gamma}{\sqrt{8}} k -\eta k - \eta n & \mbox{if } k\geq 8n\,,
\end{array}\right.}\\
$$
which is greater or equal to
$$
\left\{ \begin{array}{ll}
\frac{\eps^2\gamma}{128}n - 9\eta n & \mbox{if } k\leq 8n\,,\\
\\
\frac{\gamma}{\sqrt{8}} k -2\eta k  & \mbox{if } k\geq 8n\,.
\end{array}\right.
$$

Therefore, if the set $\{t\geq 0\,:\, Z_{\nu}(t,t)\leq -\eps t-1\}$ is unbounded, it follows that for some small $\eta>0$, the event that
\[|L\big((-k,0),(n+1,n+1)\big)-f(n+k,n)|\geq \left\{ \begin{array}{ll}
\eta n & \mbox{if } k\leq 8n\,,\\
\\
\eta k  & \mbox{if } k\geq 8n\,.
\end{array}\right.\]
happen infinitely often for $n\geq 1$ and $k\geq 1$. Using Borel-Cantelli, this will have zero probability if for all $\eta>0$,
\begin{equation}\label{eq:borcant1}
\sum_{n=1}^\infty\, \sum_{k=0}^{8n}\, \P\left(|L\big((-k,0),(n+1,n+1)\big)-f(n+k,n)|> \eta n\right) < +\infty
\end{equation}
and
\begin{equation}\label{eq:borcant2}
\sum_{n=1}^\infty\, \sum_{k=8n}^{\infty}\, \P\left(|L\big((-k,0),(n+1,n+1)\big)-f(n+k,n)|> \eta k\right) < +\infty.
\end{equation}

Theorem \ref{thm:shape} gives us some control on the fluctuations of $L$ about its asymptotic shape. Note that for $0\leq k\leq 8n$, we have
\[\begin{array}{ll}
\P\big(|L\big((-k,0),(n+1,n+1)\big)\hspace{-0.2cm}&-\ f(n+k,n)|>\eta n\big) = \vspace{0.3cm}\\
& \P\left(|L\big(\0,r(1,1)\big) - f(n+k,n)|>\eta n\right)\,.
\end{array}\]
for $r=\sqrt{(n+1+k)(n+1)}$.\\

Let $u=n^{2/3}$. If we choose $n$ large enough, we can make sure that $u<\eta n/2$ and
\[|f(n+k,n) - \gamma r|<\eta n/2.\]
Also, for $n$ large enough, we have that $u\in \left[\,c_3 \sqrt{r}\log^2 r\,,\,c_4 r^{3/2}\log r\,\right]$ (see Theorem \ref{thm:shape}). This implies, using Theorem \ref{thm:shape}, that there exist $c_i>0$ such that for $n$ large enough and $0\leq k\leq 8n$,
\begin{eqnarray*}
 \P\left(|L\big(\0,r(1,1)\big) - f(n+k,n)|>\eta n\right)& \leq & \P\big(|L\big(\0,(r,r)\big) - \gamma r|>u\big)\\
 &\leq & c_1\exp\left(-c_2n^{1/6}/\log n\right).
\end{eqnarray*}
This clearly proves \eqref{eq:borcant1}. \\

For $k\geq 8n$, we define $u=k^{2/3}$, and in a similar way we find that
\[  \P\left(|L\big(\0,r(1,1)\big) - f(n+k,n)|>\eta n\right) \leq c_1\exp\left(-c_2k^{1/6}/\log k\right)\,.\]
This proves \eqref{eq:borcant2}.\\

 The proof that the set $\{t\geq 0\,:\, Z_{\nu}(t+h,t)\geq \eps t + 1\}$ is bounded with probability 1 is actually easier, since we will have that $0\leq k\leq n$. Therefore, we only need to use the following bound on the shape function $f(x,t)=\gamma\sqrt{xt}$: for $0\leq s\leq t$
\[ f(t-s,t) \leq \frac12\gamma s - \frac{\gamma}{8}\,\frac{s^2}{t}.\]
The remainder of the argument is similar to the previous case.\\

Now we consider $\tilde{Z}_{\nu,h}$. Since for all $t$ we have 
$$\tilde{Z}_{\nu,h}(t,t)\stackrel{\cal D}{=} Z_\nu(t+h,t)\,$$ 
(but not as processes!), the Borel-Cantelli type arguments for $Z_{\nu,h}$ hold in this case as well, mutatis mutandis.

\hfill $\Box$\\

\newpage
\chapter{Busemann Functions and Equilibrium Measures}\label{ch:Busemannequi}


\section{Busemann-Equilibrium Measures}
Fix $\alpha \in (\pi,3\pi/2)$, recall \eqref{eq:defnua}, and define the measure
$$\nu_\alpha((x,y]):=\nu_\alpha^0((x,y])=B_\alpha\left((x,0),(y,0)\right)\mbox{ for }\,\,\, x\leq y\in\R\,.$$
Of course, as an immediate consequence of Proposition \ref{prop:BusMarkov}, for all $(x,t)\in \R\times[0,\infty)$,
\begin{equation}\label{eq:LastBuse}
L_{\nu_\alpha}(x,t) = B_\alpha((0,0),(x,t))\,.
\end{equation}
Now we are able to make the first step to characterize the equilibrium measures for the Hammersley Interacting Fluid process:
\begin{prop}\label{prop:equilibrium}
Fix $\alpha\in(\pi,3\pi/2)$ and denote by $\Psi_\alpha$ the probability measure on $\calN$ induced by the random measure $\nu_\alpha$. Then  $\Psi_\alpha$ is an equilibrium measure for the Hammersley Interacting Fluid Process.
\end{prop}
\noindent{\bf Proof:} By \eqref{eq:equinua} and \eqref{eq:LastBuse}, we have the following:
\begin{eqnarray}
\nonumber M^t_{\nu_\alpha}((x,y])&=& L_{\nu_\alpha}(y,t)-L_{\nu_\alpha}(x,t)\\
\nonumber &=& B_\alpha((0,0),(y,t)) -B_\alpha((0,0),(x,t))\\
\nonumber&=&\nu^t_\alpha((x,y])\\
\nonumber&\stackrel{\cal D}{=}&\nu_\alpha((x,y])\,,
\end{eqnarray}
which shows that $\Psi_\alpha$ is indeed an equilibrium measure for the Hammersley Interacting Fluid process $M^t$.

\hfill$\Box$\\

We will show that the measure $\nu_\alpha$ has the following mixing property, usually called strong mixing in dynamical systems. We consider the $\sigma$-algebra $\calF = \sigma\{ \nu_\alpha((a,b])\ :\ a\leq b \in \R\}$ on the sample space $\Omega$, defined by the compound Poisson process restricted to $\R\times\R_-$. We can define the translation $\tau_t$ as an $\calF$-measurable map from $\Omega$ to $\Omega$, simply by translating all Poisson points by the vector $(t,0)$.
\begin{prop}\label{prop:mixing}
The random measure $\nu_\alpha$ has stationary and integrable increments. Furthermore, it satisfies
\begin{equation}\label{eq:defmixing}
 \forall\ A,B\in \calF:\ \lim_{t\to \infty}\ \P(A\cap \tau^{-1}_t(B)) = \P(A)\P(B).
 \end{equation}
In particular, this implies that $\nu_\alpha$ is ergodic.
\end{prop}
\noindent{\bf Proof:} Stationarity and integrability follow from Proposition \ref{prop:Busemann} and Proposition \ref{prop:int}. From translation invariance and a standard approximation of sets in $\calF$, it is enough to prove \eqref{eq:defmixing} for all $A,B\in \calF_h:=\sigma\{ \nu_\alpha((a,b])\ :\ a\leq b \in [0,h]\}$. \\

Consider the geodesics 
$$\varpi(-\bft_\alpha,\0)\mbox{ and }\varpi(-\bft_\alpha,(h,0))\,.$$ 
where $\bft_\alpha:=(-t,-t(\tan\alpha))$. Almost surely, these paths will converge to $\varpi_\alpha(\0)$ and $\varpi_\alpha((h,0))$, respectively, on any finite box. This means, that if we define for $a,b\in [0,h]$
\[ \nu^{(t)}_\alpha((a,b]) = L(-\bft_\alpha,(b,0)) - L(-\bft_\alpha(a,0))\,,\]
then for $t$ big enough, we have $\nu^{(t)}_\alpha = \nu_\alpha|_{[0,h]}$.\\

Clearly, $\tau^{-1}_{t+h}(B)$ is independent of $\nu_\alpha^{(t)}$, since they depend on the Poisson process to the left respectively to the right of the line $\{-t\}\times \R$. Define the event
\[ C_t = \{ \forall\ s\geq t:\ \nu^{(s)}_\alpha = \nu_\alpha|_{[0,h]}\}\]
and denote $A^{(t)}$ the counterpart of the event $A$ in 
$$\calF_h^{(t)}:=\sigma\{ \nu^{(t)}_\alpha((a,b])\ :\ a\leq b \in [0,h]\}\,.$$

To see what is meant by $A^{(t)}$, we define the index-set 
$$I=\{ \nu_\alpha((a,b])\ :\ a\leq b \in [0,h]\}\,,$$
and 
$$I^{(t)}=\{ \nu^{(t)}_\alpha((a,b])\ :\ a\leq b \in [0,h]\}\,.$$
There is a canonical bijection $i:I\to I^{(t)}$. Define $\cal B$ as the product $\sigma$-algebra on $\R^I$, and likewise ${\cal B}^{(t)}$. Extend the canonical map $i$ such that $i:\R^I\to\R^{I^{(t)}}$. Define the map $\phi: \Omega \to \R^I$ by 
\[ \omega \mapsto \phi(\omega):=\{\nu_\alpha((a,b])(\omega)\ :\ a\leq b \in [0,h]\},\]
and likewise $\phi_t:\Omega \to \R^{I^{(t)}}$. We know that 
$$\calF_h = \phi^{-1}(\cal B)\mbox{ and }\calF_h^{(t)}=\phi_t^{-1}(\cal B^{(t)})\,.$$
This means that there exists $U\in \cal B$, such that $A=\phi^{-1}(U)$. We define 
$$A^{(t)} = \phi_t^{-1}(i(U))\,.$$

Thus
\begin{eqnarray*}
|\P(A\cap \tau_{t+h}^{-1}(B)) - \P(A)\P(B)| & \leq & |\P(A\cap \tau_{t+h}^{-1}(B)\cap C_t)\\
&&\qquad\qquad\qquad -\, \P(A)\P(B)| \\
&&\qquad\qquad\qquad\qquad+\, \P(C^c_t)\\
& = & |\P(A^{(t)}\cap \tau_{t+h}^{-1}(B)\cap C_t)\\
&&\qquad\qquad\qquad-\, \P(A)\P(B)|\\
&&\qquad\qquad\qquad\qquad+\, \P(C^c_t)\\
& \leq & |\P(A^{(t)}\cap \tau_{t+h}^{-1}(B))\\
&&\qquad\qquad\qquad-\, \P(A)\P(B)|\\
&&\qquad\qquad\qquad\qquad  +\, 2\P(C^c_t)\\
& = & \P(B)|\P(A^{(t)}) - \P(A)|\\
 &&\qquad\qquad\qquad\qquad+\, 2\P(C^c_t)\\
& \leq & 4\P(C^c_t).
\end{eqnarray*}
The proposition now follows from the fact that\footnote{We note that, in the classical model, we have independent increments even if the probability of the event $C_t^c$ does not decay to $0$ very fast. This indicates that, to show mixing by using these events may not be the best strategy.} $\P(C^c_t)\to 0$.

\hfill $\Box$\\

\section{Uniqueness of the Equilibrium Measure}

Now we state the most important result of this chapter.
\begin{thm}\label{thm:invariant}
If we start the Hammersley Interacting Fluid Process with $\nu_\alpha$ then
\[M_t^{\nu_\alpha}\stackrel{\cal D}{=} \nu_\alpha\,\,\,\mbox{ for all }\,\,\,t\geq 0\,.\]
The process $x\mapsto \nu_\alpha(x)$ is stationary and ergodic and its intensity is given by
\begin{equation}\label{eq:gamma}
\E \nu_\alpha(1)=\frac{\gamma(F)}{2}\sqrt{\tan\alpha}\,.
\end{equation}
Finally, consider a random $\nu\in \calN$, which is time invariant, and which defines a stationary and ergodic process on $\R$ with $\E\nu(1)\in(0,\infty)$. Define $\alpha\in (\pi,3\pi/2)$ by
\[ \alpha =\arctan \left(\frac{2}{\gamma(F)}\E \nu(1)\right)^2.\]
Then $\nu \stackrel{\cal D}{=}  \nu_\alpha$.
\end{thm}

\noindent {\bf Proof:} The first statement is an immediate consequence of Proposition \ref{prop:equilibrium}. The fact that $x\mapsto \nu_\alpha(x)$ is stationary, ergodic and that 
$\E\nu_\alpha(1)\in(0,\infty)$ follows directly from Proposition \ref{prop:mixing}. \\
 
 Fix $\alpha\in(\pi,3\pi/2)$ and set $\rho=\rho(\alpha):=\sqrt{\tan\alpha}$. Since the model is invariant under the map $(x,t)\to(\rho x,t/\rho)$,
$$\E \nu_\alpha(1)= \E \nu _{5\pi/4}(\rho) =\E \nu _{5\pi/4}(1)\sqrt{\tan\alpha}\,.$$

By Proposition \ref{prop:exitpoint}, 
$$Z_{5\pi/4}(t):=Z_{5\pi/4}(t,t)=Z_{\nu_{5\pi/4}}(t,t)\,.$$
Now, for all $t\geq 0$,
$$L_{\nu_{5\pi/4}}(t,t)=\nu_{5\pi/4}\left(Z_{5\pi/4}(t)\right)+L\big((Z_{5\pi/4}(t),0),(t,t)\big)\,.$$
If $Z_{5\pi/4}(t)\geq 0$ then
$$0\leq L_{\nu_{5\pi/4}}(t,t)-L\big(\0,(t,t)\big)\leq \nu_{5\pi/4}\left(Z_{5\pi/4}(t)\right)\,.$$
On the other hand, if $Z_{5\pi/4}(t)< 0$ then 
$$Z^*_{5\pi/4}(t):=Z^*_{\nu_{5\pi/4}}(t,t)=Z^*_{{5\pi/4}}(t,t)\geq 0\,.$$
 It follows from \eqref{eq:defflux} and \eqref{eq:LastBuse} that
\begin{equation}\label{eq:defnu*}
\nu^*_\alpha(x) = B_\alpha(\0,(0,x)).
\end{equation}
From the additivity of the Busemann function, we know that
$$L_{\nu_{5\pi/4}}(t,t)=\nu^*_{5\pi/4}\left(Z^*_{5\pi/4}(t)\right)+L\big((\0,Z^*_{{5\pi/4}}(t)),(t,t)\big)\,.$$
which finally implies that,
\begin{eqnarray*}
0&\leq& L_{\nu_{5\pi/4}}(t,t)-L\big(\0,(t,t)\big)\\
&&\qquad\qquad\qquad\leq\, \max\left\{\nu_{5\pi/4}(Z_{5\pi/4}(t)),\nu^*_{5\pi/4}(Z^*_{5\pi/4}(t))\right\}\,.
\end{eqnarray*}

Clearly, from \eqref{prop:Busesym}  (Proposition \ref{prop:Busemann}) it follows that 
$$\nu^*_{5\pi/4} \stackrel{\cal D}{=} \nu_{5\pi/4}\,.$$ 
Let $\eps > 0$. The ergodicity of $\nu_{5\pi/4}$ (and of $\nu^*_{5\pi/4}$) and the fact that  
$\E \nu_\alpha(1)\in(0,\infty)$ imply that there exists $\eta>0$ such that for all $t$ big enough,
\[ \P\left(\max\left\{\nu_{5\pi/4}(\eta t)\,,\, \nu^*_{5\pi/4}(\eta t)\right\}\geq \eps t\right)\leq \eps/2\,.\]
By \eqref{eq:vanishexitpoint}, for any $\eta>0$, 
$$\P\left(\max\left\{Z_{5\pi/4}(t),Z^*_{5\pi/4}(t)\right\}\geq \eta t\right)\leq \eps/2\,,$$ 
 for large enough $t$. Combining all this gives, for $t$ big enough,
\[ \P(L_{\nu_{5\pi/4}}(t,t)-L\big(\0,(t,t)\big)\geq \eps t)\leq \eps\,.\]
Therefore,
\begin{equation}\label{eq:int2}
\frac{L_{\nu_{5\pi/4}}(t,t)-L\big(\0,(t,t)\big)}{t} \stackrel{\cal D}{\longrightarrow} 0 \,.
\end{equation}

By \eqref{eq:LastBuse}, and using Proposition \ref{prop:Busemann}, we get
\begin{eqnarray*}
\E L_{\nu_{5\pi/4}}(t,t) & = & \E B_{5\pi/4}(\0,(t,t))\\
& = & \E B_{5\pi/4}(\0,(t,0)) + \E B_{5\pi/4}((t,0),(t,t))\\
& = & 2t\E(\nu_{5\pi/4}(1)).
\end{eqnarray*}
The ergodic theorem applied to $B_{5\pi/4}(\0,(t,t))$ implies that\footnote{The same method to prove Proposition \ref{prop:mixing} can be used to prove ergodicity of $B_{5\pi/4}$.}, with probability one,
$$\lim_{t\to\infty}\frac{B_{5\pi/4}(\0,(t,t))}{t}= 2\E \nu_{5\pi/4}(1)\,.$$
Combining this with Theorem \ref{thm:shape} and \eqref{eq:int2}, one gets \eqref{eq:gamma}.\\

Now we need to address the uniqueness of $\nu_\alpha$. Suppose $\nu\in\calN$ is ergodic and time invariant. Define $Z(t)=Z_\nu(t,(t\tan\alpha))$ and $Z_h(t)=Z_\nu(t+h,(\tan\alpha)t)$. Now define, as in Theorem \ref{thm:genexitcontrol},
\[\tilde{Z}(t) = \arg\sup_{z\leq t}\left\{ \nu(z) + L(-\bft_\alpha+(z,0),\0)\right\}\]
and
\[ \tilde{Z}_h(t) = \arg\sup_{z\leq t+h}\left\{ \nu(z) + L(-\bft_\alpha+(z,0),(h,0))\right\}.\]
Here, we take the right-most location of the maximum. The intuition for $\tilde{Z}(t)$ and $\tilde{Z}_h(t)$ is that we place again the origin at $-\bft_\alpha$, and look at the exit-point for the path that starts at $-\bft_\alpha$, picks up mass from $\nu$ and then goes to $\0$, resp. $(h,0)$.\\

Clearly, we have
\[ (Z(t), Z_h(t))\stackrel{\cal D}{=} (\tilde{Z}(t), \tilde{Z}_h(t)).\]
Since $\nu$ is ergodic, and by our choice of $\alpha$, $\nu$ satisfies \eqref{eq:genexitcontrol}. Theorem \ref{thm:genexitcontrol} then tells us that
\[ t^{-1}(\tilde{Z}(t),\tilde{Z}_h(t)) \stackrel{\rm a.s.}{\longrightarrow} \0\,.\]
This means that the two geodesics
$$\varpi(-\bft_\alpha+(\tilde{Z}(t),0),\0)\,\mbox{ and }\,\varpi(-\bft_\alpha+(\tilde{Z}(t),0),(h,0))\,$$
will converge, in any , to the $\alpha$-rays 
$$\varpi_\alpha(\0)\,\mbox{ and }\,\varpi_\alpha((h,0))$$
respectively (this follows from Theorem \ref{thm:unicoa}). However, these two $\alpha$-rays will coalesce, which means that with probability 1, there exists $t_0>0$ such that for all $t\geq t_0$, the two converging paths coalesce, which in turn implies that $\tilde{Z}(t)=\tilde{Z}_h(t)$ (because they are both the right-most point where the maximum takes place and, as soon as they coalesce, they get the same exit point). \\

Now define
\[ \tilde{L}(t) = \sup_{z\leq t}\left\{ \nu(z) + L((-t+z,-(\tan\alpha)t),\0)\right\}\]
and
\[ \tilde{L}_h(t) = \sup_{z\leq t+h}\left\{\nu(z) + L((-t+z,-(\tan\alpha)t),(h,0))\right\}.\]
We also have that
\[ (\tilde{L}(t), \tilde{L}_h(t)) \stackrel{\cal D}{=} (L_\nu(t,(\tan\alpha)t),L_\nu(t+h,(\tan\alpha)t)).\]
Furthermore, if $t\geq t_0$, then
\[ \tilde{L}_h(t) - \tilde{L}(t) = B_\alpha(\0,(h,0))=\nu_\alpha((0,h]).\]
This proves that
\begin{eqnarray*}
M^\nu_{t(\tan\alpha)}((t,t+h]) & = & L_\nu(t+h,(\tan\alpha)t) - L_\nu(t,(\tan\alpha)t)\\
 & \stackrel{\cal D}{\longrightarrow} & \nu_\alpha((0,h]).
\end{eqnarray*}
Since $\nu$ is time invariant and ergodic, we see that
\[ \nu((0,h]) \stackrel{\cal D}{=} M^\nu_{(\tan\alpha)t}((t,t+h]) \stackrel{\cal D}{=} \nu_\alpha((0,h]).\]

In principle, we need to show convergence for a finite number of $h$'s simultaneously, but it is not hard to see that the ideas we used can be extended to that case, at the cost of some notational burden. \\
 
 Note that we have proved that for any deterministic $\nu$ satisfying (\ref{eq:genexitcontrol}), 
$$M^\nu_{(\tan\alpha)t}([t,t+h]) \stackrel{\cal D}{\to}\nu_{\alpha}(h)\,$$
as a process in $h$. This shows that in a rarefaction fan, the fluid process converges locally to the correct equilibrium process (local equilibrium).

\hfill $\Box$\\

\begin{coro}\label{cor:nu*}
\[ \E(\nu_\alpha^*(1)) = \frac{\gamma(F)}{2\sqrt{\tan\alpha}}.\]
In particular, for all $\alpha\in (\pi,3\pi/2)$, we have
\[ \E(\nu_\alpha(1))\cdot \E(\nu_\alpha^*(1)) = \frac{\gamma(F)^2}{4}.\]
\end{coro}
\noindent {\bf Proof:} Remember that 
$$\nu^*_\alpha(x) = B_\alpha(\0,(0,x))\,.$$
For $\alpha = 5\pi/4$, the result follows from Proposition \ref{prop:sym} and Theorem \ref{thm:invariant}. Now use the map $(x,t)\mapsto (\rho x,t/\rho)$ to see that
\[ \nu^*_\alpha(x) \stackrel{\cal D}{=} \nu^*_{5\pi/4}\left(x/\sqrt{\tan\alpha}\right).\]

\hfill $\Box$\\

\subsection{The Classical Hammersley Process}
 
\begin{thm}\label{thm:Poisson}
Consider the classical Hammersley model and let $\Psi_\rho$ denote a probability measure on $\calN$ induced by a one-dimensional homogeneous Poisson process of intensity $\rho>0$. Then $\Psi_\rho$ is an equilibrium measure.
\end{thm}

The proof of this theorem will be an immediate consequence of a version of Burke's Theorem (Theorem \ref{thm:Burke} below). To formulate this theorem, we start by considering $M^t$ on a fixed interval $[0,R]$. Also, we choose $\nu^*$ as an independent Poisson process of intensity $1/\rho$. The evolution of $\nu$ restricted to $[0,R]$ only depends on $\nu^*$, on the Poisson process in the strip $[0,R]\times [0,\infty)$ and on the starting configuration $\nu$ on $[0,R]$. This description is called the Hammersley process with sources ($\nu|_{[0,R]}$) and sinks ($\nu^*$). In fact, $(M^t_\nu|_{[0,R]})_{t\geq 0}$ is still a Markov process, since $\nu^*$ is a Poisson process.\\

When we look at the paths of the Hammersley particles induced by $M^t_\nu$, which in fact correspond to the level sets of $L_\nu(x,t)$, they enter the box $[0,R]\times [0,T]$ either at the bottom (as sources, i.e. atoms of $\nu$), or at the east, i.e. the right-side of the box. We will call this process of ``entries'' $E_R$. Furthermore, the Hammersley paths exit the box on the left-side of the box as sinks (atoms of $\nu^*$) or at the top (atoms of $M^T_\nu$). Finally, the paths have bottom-left corners at the Poisson points of $\PP$, and we call the process of upper-right corners ${\PP}^*$.
\begin{thm}\label{thm:Burke}
Using the notations introduced above, the process $E_R$ is a Poisson process of intensity $1/\rho$, the process $M^T_\nu$ is a Poisson process of intensity $\rho$ and the process $\PP^*$ is a (two-dimensional) Poisson process of intensity 1. Furthermore, $E_R$, $M^T$ and $\PP^*$ are all independent.
\end{thm}

The idea of the proof is to show reversibility of the process $M^t_\nu$. We introduce the space $E$ as the state-space of $M^t_\nu$, so
\[ E = \sqcup_{n=0}^\infty E_n,\]
where $E_0=\{\emptyset\}$ and
\[E_n = \{(x_1,x_2,\ldots,x_n)\ :\ 0\leq x_1\leq x_2\leq \ldots \leq x_n\leq R\}.\]
We endow $E$ with the usual topology, which makes it into a locally compact space. We can view the measure $\nu$ as a measure on $E$. One of the goals of the following calculations will be to prove that $\nu$ is an equilibrium measure.\\

We can also write down the generator of the Markov process $M^t_\nu|_{[0,R]}$: if $f\in C_0(E)$, then for $x\in E$
\[ Gf(x) = \int _0^{R} f({\cal R}_tx)dt + \frac1{\rho}f({\cal L}x)-\left(\frac1{\rho}+R\right)f(x)\]
where ${\cal L}$ corresponds to an exit to the left and ${\cal R}_t$ corresponds to
an insertion of a new Poisson point at $t$. So ${\cal L}:E\to E$ where
\[ {\cal L}x = \left\{ \begin{array}{ll}
(x_2,\ldots ,x_n) & {\rm if}\ x\in E_n \ (n\geq 2)\,,\\
\emptyset & {\rm if}\ x\in E_0\sqcup E_1\,.
\end{array}
\right. \]
 And for $0<t<R$, ${\cal R}_t:E\to E$ where
\[ {\cal R}_tx=\left\{ \begin{array}{ll}
(x_1,\ldots ,x_{i-1},t,x_{i+1}, \ldots ,x_n) & {\rm if}\ x_{i-1}<t\leq x_i\ (x\in E_n)\,,\\
(x_1,\ldots ,x_n,t) & {\rm if}\ x_n<t\ (x\in E_n)\,.
\end{array}
\right. \]
Here we use the convention that $x_0=0$. \\

Note that there is a difference with Proposition \ref{prop:generator}, because of the sinks. Viewing $G$ as an operator from $L^1(\mu)$ to $L^\infty(\mu)$, we can calculate the dual operator $G^*$. This will be the generator of the time-reversed Markov process. We define ${\cal {\cal R}}$ as an exit to the right and ${\cal L}_s$ as a new point at $s$ such that the point directly to the left of $s$ moves to the right, that is ${\cal {\cal R}}:E\to E$ where 
\[ {\cal R}x = \left\{ \begin{array}{ll}
(x_1,\ldots ,x_{n-1}) & {\rm if}\ x\in E_n \ (n\geq 2)\,,\\
\emptyset & {\rm if}\ x\in E_0\sqcup E_1\,,
\end{array}
\right. \]
and, for $0<s<R$, ${\cal L}_s:E\to E$ where
\[{\cal L}_sx=\left\{ \begin{array}{ll}
(x_1,\ldots ,x_{i-1},s,x_{i+1}, \ldots ,x_n) & {\rm if}\ x_i\leq s< x_{i+1}\ (x\in E_n),\\
(s,x_1,\ldots ,x_n) & {\rm if}\ s<x_1\ (x\in E_n).
\end{array}
\right. \]
Define
\[ G^*g(y) = \int _0^{R} g({\cal L}_sy)ds + \frac1{\rho}g({\cal R}y)-\left(\frac1{\rho}+R\right)g(y)\,.\]
\begin{lem}\label{lem:dual}
For all $f,g\in C_c(\R)$,
$$\int_E Gf(x)g(x)d\Psi_\rho(x)=\int_E f(y)G^*g(y)d\Psi_\rho(y)\,.$$
\end{lem}
\noindent {\bf Proof:} See the Appendix in \cite{CaGr1}.

\hfill $\Box$\\

\noindent {\bf Proof of Theorem \ref{thm:Burke}} Lemma \ref{lem:dual} shows that $G^*$ is the generator of the time reversed process. Notice that $G^*$ exactly corresponds to the generator of a Hammersley process with sinks and sources that moves to the right: put the sinks on the right-hand side as a Poisson process of intensity $1/\rho$, and let the particles jump to the right to a Poisson process of intensity $1$ in the plane.\\
 
 We now have two conclusion. First of all, since $G^*1=0$, $\Psi_\rho$ is an equilibrium measure  for the classical Hammersley process. Secondly, since the time reversed process has the same generator as a Hammersley process moving to the right, which we denote by $\tilde{M}$, we conclude that the upper-left corners of the Hammersley paths $\PP^*$ must be a Poisson process of intensity $1$, since they correspond to the points the $\tilde{M}$ particles jump to.\\ 
 
The process of entry-points $E_R$ is slightly more subtle, since not all entries can be retraced by looking at the process $M$; after all, if $M^t_\nu=\emptyset$, there might be entries caused directly by sinks. However, this holds in exactly the same way for the right moving Hammersley process $\tilde{M}$, and therefore we can conclude that the process $E_R$, which are the entry-points for the Hammersley process, corresponds to the exit-points of the time-reversed Hammersley process, which in turn correspond to the sinks of the process $\tilde{M}$, which by construction are a Poisson process of intensity $1/\rho$. Furthermore, we can conclude that all three processes $M^T_\nu$, $\PP^*$ and $E_R$ are independent, since this is true by construction for the process $\tilde{M}$.

\hfill $\Box$\\

\begin{coro}\label{thm:classicalPoisson}
Consider the classical Hammersley model. Then  $\gamma(\delta_1)=2$ and $\nu_\alpha$ is a one-dimensional homogeneous Poisson process of intensity $\rho=\sqrt{\tan\alpha}$.
\end{coro}
\noindent {\bf Proof:} Together with Theorem \ref{thm:invariant}, Theorem \ref{thm:Poisson} implies that 
$$\nu_\rho\stackrel{\cal D}{=}\nu_{\alpha}\,\,\mbox{ with }\,\,\alpha=\alpha(\rho)=\arctan\left(\frac{2}{\gamma(\delta_1)}\rho\right)^2\,.$$
We have seen that, for $\rho=1$, $\nu^*_1(1)=\nu_1(1)=1$. Therefore, by Corollary \ref{cor:nu*},
$$\frac{\gamma(\delta_1)^2}{4}=1\,,$$
and hence $\gamma(\delta_1)=2$. In particular, $\rho(\alpha)=\sqrt{\tan\alpha}$.

\hfill $\Box$\\

\section{The Multi-Class Process}

For two positive measures $\nu$ and $\bar\nu$ on $\R$, we say that $\bar\nu$ dominates $\nu$, notation $\bar\nu\geq\nu$, whenever $\bar\nu(I)\geq \nu(I)$ for all measurable $I\subseteq \R$.
\begin{prop}\label{prop:coupling}
Suppose we have two measures $\nu, \bar\nu\in {\cal N}$ such that $\bar\nu\geq \nu$. Define the corresponding interacting fluid system as $M^\nu_t$ and $M^{\bar\nu}_t$, using the same weighted Poisson process (basic coupling). Then $M^{\bar\nu}_t\geq M^\nu_t$ (as measures). If $\bar\nu|_{(-\infty,0)}=\nu|_{(-\infty,0)}$, then  $M^{\bar\nu}_t([0,x])-M^\nu_t([0,x])$ is non-increasing in $t$ for all $x\geq 0$.
\end{prop}
\noindent{\bf Proof:} Fix an interval $[-K,K]$ and a time $t$. There exists (a random) $M>0$ such that $M^\nu_t$ and $M^{\bar\nu}_t$ restricted to $[-K,K]$ only depend on Poisson points in $[-M,K]\times [0,t]$ and on $\nu$ and $\bar\nu$ restricted to $[-M,K]$ (it is not hard to see that we can take $M=Z_\nu(-K,t)$, see \eqref{eq:orderingZnu}). This means that we are only dealing with a finite number of Poisson points, so if we can prove that the premise ``$M^{\bar\nu}_s\geq M^\nu_s$ for all $s<t$'' implies that $M^{\bar\nu}_t\geq M^\nu_t$, we will have proved the first statement, since it is obviously true for $t=0$.\\

Suppose there exists a Poisson point at $(x_0,t)$ with weight $\omega$ for some $x_0\in [-M,K]$, since otherwise the implication is immediate. We then know, using Proposition \ref{prop:construction} and \eqref{eq:evolM}, that if $x_0< x\leq y$,
\begin{eqnarray*}
M^{\nu}_t((x,y]) & = & (M^{\nu}_{t-}((x_0,y]) - \omega)_+  - (M^{\nu}_{t-}((x_0,x]) - \omega)_+ \\
& \leq & (M^{\bar\nu}_{t-}((x_0,y]) - \omega)_+  - (M^{\bar\nu}_{t-}((x_0,x]) - \omega)_+ \\
& = & M^{\bar\nu}_t((x,y]).
\end{eqnarray*}
The inequality follows from the fact that if $A\geq B$ and $\tilde{A}\geq \tilde{B}\geq 0$, then $(A-\omega)_+-(B-\omega)_+\leq (A+\tilde{A}-\omega)_+ - (B+\tilde{B}-\omega)_+$. If $x\leq x_0< y$ or $x\leq y\leq x_0$, the implication is straightforward, following a similar split up.\\

The second statement follows from a similar reasoning: suppose there is a Poisson point at $(x_0,t)$ with weight $\omega$. If $x> x_0\geq 0$,
\begin{eqnarray*}
&&M^{\bar\nu}_t((x_0,x])-M^\nu_t((x_0,x])  =  (M^{\bar\nu}_{t-}((x_0,x]) - \omega)_+\\
&&\qquad\qquad\qquad\qquad\qquad\qquad\qquad\qquad -\, (M^{\nu}_{t-}((x_0,x]) - \omega)_+\\
&&\qquad\qquad\qquad\qquad\qquad\qquad\leq\, M^{\bar\nu}_{t-}((x_0,x])-M^{\nu}_{t-}((x_0,x])\,.
\end{eqnarray*}
The inequality follows from the fact that 
$$(A-c)_+-(B-c)_+\leq A_+-B_+$$ 
whenever $c\geq 0$ and $A\geq B$. Since 
$$M^{\bar\nu}_t([0,x_0])-M^\nu_t([0,x_0]) =  M^{\bar\nu}_{t-}([0,x_0])-M^{\nu}_{t-}([0,x_0])\,,$$ 
this shows that
\[ M^{\bar\nu}_t([0,x])-M^\nu_t([0,x])\leq M^{\bar\nu}_{t-}([0,x])-M^{\nu}_{t-}([0,x]).\]

Now suppose $x_0<0$ and $x\geq 0$. Note that under the condition on $\bar \nu$, we have that for all $s\geq 0$ and all $\eps>0$, $L_{\bar\nu}(-\eps,s)=L_{\nu}(-\eps,s)$, so
\[ M^{\bar\nu}_s([0,x])-M^\nu_s([0,x]) = L_{\bar\nu}(x,s)-L_\nu(x,s).\]
When $L_{\bar\nu}(x,t)$ does not use the weight at $(x_0,t)$, we know that $L_{\bar\nu}(x,t)=L_{\bar\nu}(x,t-)$ and that $L_{\nu}(x,t)\geq L_{\nu}(x,t-)$, which implies the desired result. If $L_{\bar\nu}(x,t)$ does use the weight at $(x_0,t)$, then it is not hard to see that $L_\nu(x,t)$ will also use the weight at $(x_0,t)$ (the longest path corresponding to $\bar \nu$ is always to the right of the path corresponding to $\nu$), which means that only the mass on the $x$-axis strictly to the left of $0$ is used, and therefore 
$$M^{\bar\nu}_t([0,x])=M^\nu_t([0,x])\,.$$

Finally, when $x_0=x$, we get that 
$$M^{\bar\nu}_{t}([0,x])=M^{\bar\nu}_{t-}([0,x])+\omega$$ 
and 
$$M^{\nu}_t([0,x])=M^\nu_{t-}([0,x])+\omega\,,$$ 
and when $x_0>x$, we have 
$$M^{\bar\nu}_{t}([0,x])=M^{\bar\nu}_{t-}([0,x])\,$$
and 
$$M^{\nu}_t([0,x])=M^\nu_{t-}([0,x])\,.$$

\hfill $\Box$\\

In other words, Proposition \ref{prop:coupling} tells us that the Hammersley Interacting Fluid System is monotone: if one starts the fluid process with the same Poisson weights (basic coupling) and with ordered initial configurations, then the order is preserved for all $t\geq 0$. This coupled process is called the Multi-Class Fluid System. It is just a convention to describe a coupled process with ordered initial configurations \cite{FeMa}.\\

For any countable $D\subseteq(\pi,3\pi/2)$, one can a.s. construct simultaneously a collection of equilibrium processes $\{\Psi_\alpha:\alpha\in D\}$ by using the same Poisson weights on $\R\times\R_-$ and the Busemann functions $B_{\alpha}$. It turns out that this collection respects the order induced by the angles $\alpha\in D$. More precisely:
\begin{thm}\label{thm:multi-class}
If $\bar\alpha>\alpha$ then $\nu_{\bar\alpha}\geq\nu_{\alpha}$. In particular, for any countable subset  $\{\alpha_{i}\,:\,i\in\ZZ\}\subseteq D$, if one runs simultaneously (basic coupling) the Interacting Fluid Processes on $\R\times\R_+$ with initial measures $(\nu_{\alpha_i}\,:\,i\in\ZZ)$ 
then, whenever $\alpha_i>\alpha_j$, $M^{\nu_{\alpha_i}}_t\geq M^{\nu_{\alpha_j}}_t$ for all $t\geq 0$, and
$$(M^{\nu_{\alpha_i}}_t\,:\,i\in \ZZ)\stackrel{\cal D}{=}(\nu_{\alpha_i}\,:\,i\in \ZZ)\,.$$
\end{thm}
\noindent{\bf Proof:}
Almost sure coalescence of $\alpha$-rays for fixed $\alpha\in(\pi,3\pi/2)$ implies almost sure coalescence for a given countable $D\subseteq(\pi,3\pi/2)$. This allows us to construct simultaneously the Busemann functions $\{B_\alpha\,:\,\alpha\in D\}$, as a function of the underlying compound Poisson process. Thus, time invariance of the Busemann multi-class measure $(\nu_{\alpha_i}\,:\,i\in\ZZ)$ follows from translation invariance of the compound Poisson process (as in the proof of Proposition \ref{prop:Busemann}).\\

To see that it is indeed a multi-class measure, let $\mm$ be the crossing point between $\varpi_{\bar\alpha}((z,0))$ and $\varpi_\alpha((z',0))$. Furthermore, denote $\cc$ as the coalescence point of the two $\alpha$-rays $\varpi_\alpha((z,0))$ and $\varpi_\alpha((z',0))$, and denote $\bar\cc$ as the coalescence point of the two $\bar\alpha$-rays $\varpi_{\bar\alpha}((z,0))$ and $\varpi_{\bar\alpha}((z',0))$. Then
\begin{eqnarray}
\nonumber \nu_{\bar\alpha}\big([z,z']\big)- \nu_{\alpha}\big([z,z']\big)&=& \Big\{L(\bar\cc,(z',0))-L(\bar\cc,(z,0))\Big\}\\
\nonumber&-&\Big\{ L(\cc,(z',0))-L(\cc,(z,0)) \Big\}\\
\nonumber&=& L(\bar\cc,(z',0))\\
\nonumber&&\qquad\qquad-\big\{L(\bar\cc,\mm)+L(\mm,(z',0))\big\}\\
\nonumber&+& L(\cc,(z,0))\\
\nonumber&&\qquad\qquad-\big\{L(\cc,\mm)+L(\mm,(z,0))\big\}\\
\nonumber&\geq& \,\,0\,.
\end{eqnarray}

\hfill $\Box$\\

This result is also new in the classical Hammersley Interacting System, where a different and explicit description of the multi-class equilibrium measure with a finite number of classes is given by Ferrari and Martin \cite{FeMa}.

\newpage
\chapter{Second-Class Particles}\label{ch:2ndclass}

\section{Second-Class Particles and Exit Points}
Proposition \ref{prop:coupling} can be used to define the notion of second-class particles. In the interacting fluid system we can define it analogously to the interacting particle case, with a slight adaptation due to the continuous weights. \\

We start by changing $\nu$ into $\bar\nu$, by putting an extra weight $\eps>0$ in $0$, so
\[\bar\nu([0,x])=\nu([0,x])+\eps\,\,\,\mbox{ for }\,\,\,x\geq 0\,.\]
With this new process, and using the same Poisson weights, we define $M^{\bar\nu}_t$. Clearly,
\[ M^{\bar\nu}_t([0,x])\leq M^\nu_t([0,x]) + \eps.\]
Now define the location of the second class particle $X_\nu(t)$ as
\[ X_\nu(t) = \inf\{x\geq 0: M^{\bar\nu}_t([0,x])=M^\nu_t([0,x])+\eps\}\,.\]

By Proposition \ref{prop:coupling}, $X_\nu(t)$ is a non-decreasing function of $t$, meaning that the second class particle moves to the right. In fact, the extra mass $\eps$ will spread out, and our definition coincides with the rightmost point of this spread-out mass. This is a natural choice, since we will show that it does not depend on the total mass $\eps$, while for example the leftmost point does depend on $\eps$.\\

There is the following important connection between the longest path description and the second class particle. Let $\nu^+$ be the process defined by $\nu^+(x)=\nu(x)$ for $x\geq0$, and by $\nu^+(x)=-\infty$ for $x<0$. We also define the process $\nu^-$ by setting $\nu^-(x)=0$ for $x\geq0$, and $\nu^-(x)=\nu(x)$ for $x<0$. Then
\[ L_{\nu^+}(x,t) = \left\{ \begin{array}{ll}
\sup\{L((z,0),(x,t)) + \nu(z)\ :\ 0\leq z\leq x\} & \mbox{if}\ x\geq 0\\
-\infty & \mbox{if}\ x<0\,,
\end{array}\right.\]
and
\[ L_{\nu^-}(x,t) = \sup\{L((z,0),(x,t)) + \nu(z)\ :\ z<0\ \mbox{and}\ z\leq x\}\,.\]
Clearly,
$$L_\nu(x,t)=\max\big\{L_{\nu^+}(x,t),L_{\nu^-}(x,t)\big\}\,.$$

Now suppose $x\geq 0$. If $L_{\nu^+}(x,t)\geq L_{\nu^-}(x,t)$, there exists a longest path that does not use any weight of $\nu$ on $(-\infty,0)$. This means that if we add a weight $\eps>0$ in the origin, 
$$L_{\bar\nu}(x,t)=L_{\bar\nu^+}(x,t)=L_\nu(x,t)+\eps\,.$$ 
Using Proposition \ref{prop:construction}, we see that this means that 
$$M^{\bar\nu}_t(x)=M^\nu_t(x)+\eps\,$$
so $X_\nu(t)\leq x$.\\

If on the other hand we start with $X_\nu(t)\leq x$, we conclude that 
$$M^{\bar\nu}_t(x)=M^\nu_t(x)+\eps\,$$
using Proposition \ref{prop:coupling} and the fact that $M^\nu_t$ and $M^{\bar\nu}_t$ are right-continuous. This in turn means that 
$$L_{\bar\nu}(x,t)=L_{\nu}(x,t)+\eps\,,$$ 
which is only possible if $L_{\nu^+}(x,t)\geq L_{\nu^-}(x,t)$.\\

We have shown that
\begin{equation}\label{eq:secondclass}
\{ X_\nu(t)\leq x\} = \{L_{\nu^+}(x,t)\geq L_{\nu^-}(x,t)\}\,.
\end{equation}
Note that this can be rewritten as
\begin{equation}\label{eq:secondclassZ}
\{ X_\nu(t)\leq x\} = \{ Z_\nu(x,t) \geq 0\}.
\end{equation}
This means that the path of the second class particle corresponds to a competition interface, a fact well known for the totally asymmetric exclusion process \cite{FePi}. This allows us to show that the second class particle satisfies a strong law whenever $\nu^+$ and $\nu^-$ have asymptotic intensities. The proof of this does not use a coupling of two invariant versions of the fluid process, as is usual in the interacting particle case, but it uses the longest path description in a direct way. \\

We would like to point out that in our general set-up, with random weights on the Poisson points, we do not have an equivalent of Burke's Theorem. This means that the time-reversed process is not a Hammersley interacting fluid system. Therefore, the path of a second class particle in general does not coincide in law with a longest path in the interacting fluid system, in contrast to the classical case, where the statement is true \cite{CaGr2}. However, we do have the following connection.
\begin{prop}\label{prop:second-class}
Assume that the distribution of $\nu$ is translation invariant. Then, for any $t\geq 0$, we have that
\[ X_{\nu}(t) - x \stackrel{\cal D}{=} -Z_\nu(x,t).\]
\end{prop}
\noindent{\bf Proof:} This follows almost immediately from \eqref{eq:secondclassZ}, since that equality can be rewritten as
\[ \{ X_{\nu}(t) - x  \leq h\} = \{ Z_\nu(x + h,t) \geq 0\}.\]
Now use translation invariance to see that
\[ Z_\nu(x + h,t) \stackrel{\cal D}{=} Z_\nu(x,t) + h.\]

\hfill $\Box$\\

When we consider all $\alpha$-rays starting at the line $\R\times\{t\}$ and we move from left to right, $(X_{\nu_\alpha}(t),t)$ is the first point where the $\alpha$-ray passes the origin. It is tempting to think that the $\alpha$-ray starting at $(X_{\nu_\alpha}(t),t)$ actually passes through the origin, but this is false in general. In fact, after time $t$, most $\alpha$-rays will have coalesced with other rays, and the crossings with the $x$-axis will be quite far apart; we would conjecture they are order $t^{2/3}$ apart.

\section{Law of Large Numbers for Second-Class Particles}
Proposition \ref{prop:second-class} allows us to use Theorem \ref{thm:shape} and Theorem \ref{thm:genexitcontrol} to prove a strong law for the second class particle in the case of $\nu_\alpha$.  However, we are able to prove a strong law even for deterministic initial conditions that satisfy a density property given in the following lemma, whose proof is very similar to Theorem \ref{thm:genexitcontrol}.
\begin{lem}\label{lem:genexitcontrol2}
Suppose $\nu\in \calN$. Assume that
\begin{equation}\label{eq:genexitcontrol2a}
\liminf_{d\to\infty} \inf_{z\in [-4d,4d]} \frac{\nu((z-d,z+d])}{2d}\geq\frac{\gamma}{2}\sqrt{\tan\alpha}\,.
\end{equation}
Let $\eps>0$. Then, with probability one,
\[ \liminf_{t\to \infty} \frac{Z_\nu((1+\eps)t,(\tan\alpha)t)}{t}>0.\]
Now assume that
\begin{equation}\label{eq:genexitcontrol2b}
\limsup_{d\to\infty} \sup_{z\in [-4d,4d]} \frac{\nu((z-d,z+d])}{2d}\leq\frac{\gamma}{2}\sqrt{\tan\alpha}\,.
\end{equation}
Then, with probability one,
\[ \limsup_{t\to \infty} \frac{Z_\nu((1-\eps)t,(\tan\alpha)t)}{t}<0.\]
\end{lem}

\noindent{\bf Remark:} Note that in our density condition, we do not allow the midpoint of the interval to be much larger than $d$. The reason for this might be more clear if we think of a Poisson process: if we fix $d$, we can always find some $z\in \R$ such that the interval $[z-d,z+d]$ is empty!\\

\noindent{\bf Proof of Lemma \ref{lem:genexitcontrol2}}: The proof of this lemma relies on the proof of Theorem \ref{thm:genexitcontrol}. As usual, we will assume without loss of generality that $\alpha=5\pi/4$. We will start with the first statement.\\

We need to prove that for all $\eta>0$ small enough, the set 
$$\{t\geq 0\ :\ Z_\nu((1+2\eps)t,t)\leq \eta t\}\,$$
is bounded with probability one. Choose $\eta<\eps/2$. Define the translated measure $\nu_{2\eps t}$ such that for all $b\geq a$,
\[ \nu_{2\eps t}((a,b]) = \nu((a+2\eps t,b+2\eps t])\,.\]
Using the definition in \eqref{eq:defnuprocess}, we see that for $z\in \R$,
\[ \nu_{2\eps t}(z) = \nu(z+2\eps t) - \nu(2\eps t)\,.\]
Define $k=-\lceil Z_\nu((1+2\eps)t,t)-2\eps t\rceil$ and $n=\lfloor t\rfloor$. If $t$ is such that 
$$Z_\nu((1+2\eps)t,t)\leq \eta t\,,$$ 
then
\begin{eqnarray*} 
&&\nu(-k+1+2\eps t) + L((-k+2\eps t,0),((1+2\eps)t,t))\\
&&\qquad\qquad\qquad\qquad - \big[ \nu(2\eps t) + L((2\eps t,0), ((1+2\eps)t,t))\big] \geq 0\,.
\end{eqnarray*}
This event has the same probability as the event
\[ \nu_{2\eps t}(-k+1) + L((-k,0),(t,t)) - L((0,0),(t,t)) \geq 0\,.\]

We know that $k\geq \eps t$ (since $\eta <\eps$). We can use \eqref{eq:genexitcontrol2a} to see that for $t$ big enough,
\[ -\nu_{2\eps t}(-k+1) = \nu((-k+1+2\eps t, 2\eps t])\geq \frac12 \gamma k - \eta k\,.\]
To see this, define half the length of the interval by $d=(k-1)/2$, and the midpoint by $z=\eps t -(k-1)/2$. We see then see that $|z| \leq 4d$. This gives us, using the fact that $n\leq t\leq n+1$,
\[ L((-k,0),(n+1,n+1)) - L((0,0),(n,n)) \geq \frac12 \gamma k - \eta k\,.\]
This is exactly \eqref{eq:bendL} in the proof of Theorem \ref{thm:genexitcontrol}, and we can follow that proof from this point. The second statement of the lemma follows similarly.

\hfill $\Box$\\

With this lemma we can proof the following result:
\begin{thm}\label{thm:second-class}
Assume that
\begin{equation}\label{eq:second-int}
\lim_{x\to\infty}\frac{\nu\big(x)}{x}=\lim_{x\to-\infty}\frac{\nu(x)}{x}=\frac{\gamma}{2}\sqrt{\tan\alpha}\,.
\end{equation}
Then, with probability one,
$$\lim_{t\to\infty}\frac{X_\nu(t)}{t}=\frac{1}{\tan\alpha}\,.$$
\end{thm}
\noindent {\bf Proof:} As usual, we will assume without loss of generality that $\alpha=5\pi/4$. Suppose $\eps>0$ and $X_\nu(t)\leq t - 2\eps t$. Define $n=\lfloor t\rfloor$. Then for $t$ large enough, we have
\[X_\nu(n)\leq (1-2\eps)n + 1 + 2\eps\leq (1-\eps)n\,.\]
By  \eqref{eq:secondclassZ}, this implies that $Z_\nu((1-\eps)n,n)\geq 0$.\\

If we can show that $\nu$ satisfies \eqref{eq:genexitcontrol2b}, then we can use Lemma \ref{lem:genexitcontrol2} to see that $Z_\nu((1-\eps)n,n)\geq 0$ can happen only for finitely many $n\geq 1$, which gives
\[ \liminf_{t\to \infty} \frac{X_\nu(t)}{t} \geq 1\,.\]
Bounding the limit from above can be done using the analogous argument.\\

To see that \eqref{eq:genexitcontrol2b} indeed holds, remark that for all $1>\eta>0$, there exists $R>1$ such that for all $z\geq R$,
\[ \left|\frac{\nu((0,z])}{z} - \frac12\gamma\right| < \eta \ \ \mbox{and}\ \ \left|\frac{\nu((-z,0])}{z} - \frac12\gamma\right| < \eta.\]
Choose $M>1$ such that $\nu([-R,R])<M$. Now choose $d>MR/\eta$. If $z>d+R$ or $z<-d-R$, we get that
\begin{eqnarray*}  
\nu((z-d,z+d])&\leq& d\gamma + (2|z|+2d)\eta\\
& \leq& d\gamma + 10d\eta\,.
\end{eqnarray*}
If $|z+d|<R$,
\begin{eqnarray*}
 \nu((z-d,z+d])&\leq& \frac12(2d+R)\gamma + M + |z-d|\eta \\
 &\leq& d\gamma + (5d+2+\frac12\gamma)\eta\,,
 \end{eqnarray*}
and a similar bound holds when $|z-d|<R$. This proves that $\nu$ satisfies \eqref{eq:genexitcontrol2b}.

\hfill $\Box$\\

\section{The Rarefaction Fan}
Assume now that we have a non-homogeneous initial measure $\nu$ that satisfies the following rarefaction assumption:
\begin{equation}\label{eq:assumpint}
a_\nu:=\limsup_{x\to \infty} \frac{\nu^+(x)}{x} < \liminf_{y\to -\infty} \frac{\nu^-(y)}{y}=:b_\nu\,.
\end{equation}
In this regime the behavior of the second class particle is different from the previous one \cite{CaPi2}.
\begin{thm}\label{thm:2ndclasspart}
For each $v\in\left[\frac{\gamma^2}{4b_\nu^{2}},\frac{\gamma^2}{4a_\nu^{2}}\right]$ let
$$\alpha(v):=\arctan(1/v)\in(0,\pi/2)\,.$$
If \eqref{eq:assumpint} holds then
$$\frac{X_\nu(t)}{t}\stackrel{\cal D}{\to}V_\nu\mbox{ as }t\to\infty\,,$$
where $V_\nu$ is a random variable with support $\left[\frac{\gamma^2}{4b_\nu^{2}},\frac{\gamma^2}{4a_\nu^{2}}\right]$ and
\begin{equation}\label{2class-dist}
\P\left(V_\nu \leq v\right) = \P\left(\sup_{z\geq 0}  \{\nu(z) -\nu_{\alpha(v)}(z)\} \geq \sup_{z<0} \{\nu(z) -\nu_{\alpha(v)}(z)\}\right)\,.
\end{equation}
for $v\in \left(\frac{\gamma^2}{4b_\nu^{2}},\frac{\gamma^2}{4a_\nu^{2}}\right)$, where $\nu_\alpha$ and $\nu$ are independent random measures.
\end{thm}

 \noindent{\bf Proof:} Condition \eqref{eq:assumpint} ensures that with probability one, at least one of the two suprema is finite, and hence the right hand side probability is well defined for every $v\in \left(\frac{\gamma^2}{4b_\nu^{2}},\frac{\gamma^2}{4a_\nu^{2}}\right)$. (Notice that at the boundary, \eqref{2class-dist} may fail!) \\
 
 For each $v,t>0$  let 
$$\alpha=\alpha_v=\arctan(1/v)\in(0,\pi/2)\,,$$
$$x=x_t:=t/ \tan\alpha\,\,\mbox{ and }\,\,\xx_\alpha := (x,x\tan\alpha)\,.$$ 
\eqref{eq:secondclass} can be rephrased to
$$\{ X_\nu(t)\leq x\} = \{L_{\nu^+}(x,t)-L(\0,(x,t))\geq L_{\nu^-}(x,t)-L(\0,(x,t))\}\,.$$
Thus,
\begin{eqnarray}
\nonumber \{ X_\nu(t)\leq x\} & = & \{ L_{\nu^+}(\xx_\alpha)\geq L_{\nu^-}(\xx_\alpha)\}\\
\nonumber& = & \Big\{ \sup_{0\leq z\leq x} \left[ \nu(z) + L((z,0),\xx_\alpha)\right] \\
\nonumber&&\qquad\qquad\geq \sup_{ z< 0} \left[ \nu(z) + L((0,z),\xx_\alpha)\right]\Big\}\\
\nonumber& = & \Big\{ \sup_{0\leq z\leq x} \left[ \nu(z) + L((z,0),\xx_\alpha) - L(\0,\xx_\alpha)\right] \\
\nonumber& & \ \ \ \ \ \ \ \ \geq \sup_{z<0} \left[ \nu(z) + L((0,z),\xx_\alpha)- L(\0,\xx_\alpha)\right]\Big\}\,.\\
\label{2classBuse}
\end{eqnarray}

Notice that, for any compact set $K\subset \R^2$,
\begin{equation}\label{eq:busconv}
\left(L(\zz,\xx_\alpha) - L(\0,\xx_\alpha)\,;\, \zz\in K\right) \stackrel{\cal D}{\longrightarrow} \left(B_\alpha(\0,\zz)\,;\, \zz\in K\right)\,\,\mbox{ as }\,\,x\to \infty\,.
\end{equation}
This remark is the core of the proof. It follows from Theorem \ref{thm:unicoa}: we can take $n\geq 1$ big enough, so that $K\subset [-n,n]\times [-n,n]$. The $\alpha$-rays starting at $(-n,n)$ and $(n,-n)$ will coalesce at some point $\cc\in \R^2$. Furthermore, the longest paths to $\xx_\alpha$, starting at $(-n,n)$ and $(n,-n)$, will converge to the respective $\alpha$-rays in a bigger bounded square, containing $\cc$. From that time on, for all $\zz\in K$, we will have
\[   L(\zz,\xx_\alpha) - L(\0,\xx_\alpha) = B_\alpha(\0,\zz).\]

Notice that, since $\nu_\alpha$ is a function of the underlying two dimensional compound Poisson random set $\PP$ restricted to the upper half plane, it is independent of $\nu$. Together with \eqref{2classBuse} and Lemma \ref{2classcontrol} below, this proves the theorem for $\frac{\gamma}{2}\sqrt{\tan\alpha}\in(a_\nu,b_\nu)$.\\
  
 Now, if $b_\nu <\frac{\gamma}{2}\sqrt{\tan\alpha}$ then a.s. 
 $$\sup_{y< 0} \nu^-(y)-\nu^{-}_{\alpha}(y) = \infty\,.$$ 
 By the same reasoning as in the  proof of Lemma \ref{2classcontrol},
$$\sup_{z\leq 0} \left\{ \nu(z) + L((z,0),\xx_{\alpha}) - L(\0,\xx_{\alpha})\right\}\stackrel{\cal D}{\longrightarrow} \infty\,\,,\mbox { as }\,\,x\to\infty\,.$$
Hence, both functions in \eqref{2class-dist} are zero. The case where $a_\nu>\frac{\gamma}{2}\sqrt{\tan\alpha}$ follows from a similar argument.

\hfill $\Box$\\

\begin{lem}\label{2classcontrol}
If $\frac{\gamma}{2}\sqrt{\tan\alpha}\in\left(a_\nu,b_\nu\right)$ then, as $x\to\infty$,
\[ \sup_{0\leq z\leq x} \left\{ \nu(z) + L((z,0),\xx_{\alpha}) - L(\0,\xx_{\alpha})\right\}\stackrel{\cal D}{\longrightarrow} \sup_{z\geq 0} \left\{\nu(z) - \nu_{\alpha}(z)\right\}\,,\]
and
\[ \sup_{z< 0} \left\{\nu(z) + L((z,0),\xx_{\alpha}) - L(\0,\xx_{\alpha})\right\}\stackrel{\cal D}{\longrightarrow} \sup_{z< 0} \left\{ \nu(z) - \nu_{\alpha}(z)\right\}.\]
\end{lem}
\noindent{\bf Proof:}  We will prove that $a_\nu < \frac{\gamma}{2}\sqrt{\tan\alpha}$ implies the first statement. The proof of the second statement follows exactly the same reasoning. \\

Pick $\alpha'\in(0,\pi/2)$ and $a'>0$ such that 
$$a_\nu<a'=\frac{\gamma}{2}\sqrt{\tan\alpha'}<\frac{\gamma}{2}\sqrt{\tan\alpha}\,.$$ 
Let $\tilde Z_{\alpha'}(x)$ be the crossing point between $\varpi_{\alpha'}(\0)$ and the line $\{(y,t)\,:\,y\in\R\}$ (recall that $t=x/\tan\alpha$). Let $E_x$ be the event that $\tilde Z_{\alpha'}(x)$ is to the right of $\xx_\alpha$. By Proposition \ref{prop:Znuexists} (local comparison), we have that, under the event $E_x$,
\[ L(\0,\xx_\alpha) - L((0,z),\xx_\alpha) \geq \nu_{\alpha'}(z)\,\]
(translate the origin to $\xx_\alpha$ and reflect both coordinates), and hence
\[ \nu(z) - [L(\0,\xx_\alpha) - L((0,z),\xx_\alpha)] \leq\nu(z) - \nu_{\alpha'}(z)\,.\]
From Theorem \ref{thm:genexitcontrol} (recall that $\alpha'<\alpha$), we conclude that
$$\lim_{x\to\infty}\P\left(E_x\right)=1\,.$$
By our choice of $\alpha'$, $a'>a_\nu$, and thus there exists $M>0$ (random) such that
$$0\leq\sup_{z\geq0}\{\nu(z) - \nu_{\alpha'}(z)\}=\sup_{0\leq z\leq M}\{\nu(z) - \nu_{\alpha'}(z)\}<\infty\,.$$
 This means that for each $\eps>0$, we can find $R>0$ such that for all $x$ large enough, the probability of the event that 
\begin{eqnarray*}
&& \sup_{0\leq z\leq x} \left\{ \nu(z) + L((z,0),\xx_\alpha) - L(\0,\xx_\alpha)\right\}\\
&&\qquad\qquad\qquad\qquad  >  \sup_{0\leq z\leq R} \left\{ \nu(z) + L((z,0),\xx_\alpha) - L(\0,\xx_\alpha)\right\}
\end{eqnarray*}
is smaller than $\eps$. Now we can use \eqref{eq:busconv} to finish the proof of the first statement of Lemma \ref{2classcontrol}.

\hfill $\Box$\\

In the classical Hammersley process $\nu_\alpha$ is a Poisson process of intensity
$$\rho=\sqrt{\tan\alpha}=\frac{1}{\sqrt v}\,,$$
which allows us to compute the law of the asymptotic speed in some cases.

\subsection{Poisson Initial Configuration} Choose $\lambda>0$ and $\mu>\lambda$, and assume that $(\nu^+(y)\,;\,y\geq 0)$ and $(\nu^-(y)\,;\,y<0)$ are independent Poisson counting processes of intensity $\lambda$ and $\mu$, respectively. We will only consider the case $\rho\in(\lambda,\mu)$, the other cases are trivial.\\

Define two asymmetric simple random walks $W^+$ and $W^-$, with
\[\P(W^+(n+1)-W^+(n) = +1)= p^+ := \frac{\lambda}{\lambda + \rho}\]
and
\[\P(W^-(n+1)-W^-(n) = +1)= p^- := \frac{\rho}{\mu + \rho}.\]
Since $\nu^+$, $\nu^+_{\rho}$, $\nu^-$ and $\nu^-_{\rho}$ are independent Poisson counting process, it is not hard to see that
$$\sup_{z\geq 0}\{ \nu^+(z)-\bar\nu^{+}_{\rho}(z)\} \stackrel{\cal D}{=} \sup_{n\geq 0} W^+(n)\,,$$
and that 
$$\sup_{z< 0} \{\nu^-(z)-\bar\nu^{-}_{\rho}(z)\} \stackrel{\cal D}{=}  \sup_{n\geq 0} W^-(n)\,.$$
Furthermore, it is well known (and easy to see) that 
$$S^i:=\sup_{n\geq 0} W^i(n)\sim {\rm Geo}(r^i)\,,$$ 
where 
$$r^i=\frac{p^i}{1-p^i}\,,$$
since $p^i<0.5$ (for $i=+,-$). We find
\begin{eqnarray*}
\P(S^+\geq S^-) & = & \sum_{k=0}^\infty \P(S^+\geq k)\P(S^-=k)\\
& = & \sum_{k=0}^\infty (r^+)^k(r^-)^{k}(1-r_-) \\
&=&\frac{1-r^-}{1-r^+r^-} =\frac{\mu-\rho}{\mu-\lambda}\,.
\end{eqnarray*}
Therefore, by Theorem \ref{thm:2ndclasspart} (and taking $\rho_v=1/\sqrt{v}$), this proves that
\[ \P(V^\nu \geq v) = \left\{\begin{array}{ll}
0 & \mbox{if }\,\, v\geq 1/ \lambda^2\\
\frac{\frac{1}{\sqrt{v}}-\lambda}{\mu-\lambda} & \mbox{if }\,\, 1/\mu^2< v< 1/\lambda^2\\
1 & \mbox{if }\,\,v\leq 1/\mu^2.
\end{array}
\right. \]
This agrees with the results found in \cite{CaDo,CoPi}.

\subsection{Periodic Initial Configuration} Choose $\lambda>0$ and $\mu>\lambda$, and assume that $(\nu^+(y)\,;\,y\geq 0)$ and $(\nu^-(y)\,;\,y<0)$ are deterministic periodic configurations of intensity $\lambda$ and $\mu$, i.e., $\nu^+$ is concentrated on $\{k/\lambda\,:\,k\geq 1\}$ while $\nu^-$ is concentrated on $\{k/\mu\,:\,k\leq -1\}$. Suppose $\lambda<\rho<\mu$. Theorem \ref{thm:2ndclasspart} then tells us that we should consider
\[ \P\left(\sup_{z\geq 0}  \{\nu(z) -\nu_{\rho}(z)\} \geq \sup_{z<0} \{\nu(z) -\nu_{\rho}(z)\}\right)\,,\]
where $\nu_\rho$ is a Poisson process on $\R$ of intensity $\rho$. \\

Define
\[ S^+ = \sup_{z\geq 0}  \{\nu(z) -\nu_{\rho}(z)\}\ \ \ \mbox{and}\ \ \ S^- = \sup_{z<0} \{\nu(z) -\nu_{\rho}(z)\}.\]
Clearly, for $z\geq 0$, we have $\nu(z) = \lfloor\lambda z\rfloor$, where $ \lfloor x\rfloor$ denotes the integer part of $x$. Therefore,
\begin{eqnarray*}
S^+ & = & \sup_{z\geq 0} \{ -\nu_\rho(z) + \lfloor\lambda z\rfloor\}\\
& = & \lfloor\sup_{z\geq 0} \{ -\nu_\rho(z) + \lambda z\}\rfloor\\
& = & -\lceil\inf_{z\geq 0} \{ \nu_\rho(z) - \lambda z\}\rceil.
\end{eqnarray*}

The infimum of a Poisson process minus a linear function is studied in \cite{Py}. Theorem 3 in \cite{Py} entails that for $k\geq 0$
\begin{eqnarray*}
\P(S^+\geq k) & = & \P(\lceil\inf_{z\geq 0} \{ \nu_\rho(z) - \lambda z\}\rceil \leq -k)\\
& = & \P(\inf_{z\geq 0} \{ \nu_\rho(z) - \lambda z\}\leq -k)\\
& = & (1-p_+)^k,
\end{eqnarray*}
where $p_+=p_+(\lambda,\rho)\in (0,1)$ is the positive solution of
\[ p_+ = 1-e^{-p_+\rho/\lambda}.\]
For $z<0$ we have $\nu(z)=-\lfloor\mu |z|\rfloor$. So when we switch to positive $z$, we get
\begin{eqnarray*}
S^- & \stackrel{\cal D}{=} & \sup_{z\geq 0} \{\nu_\rho(z) - \lfloor\mu z\rfloor\}\\
& = & \lceil\sup_{z\geq 0} \{\nu_\rho(z) - \mu z\}\rceil.
\end{eqnarray*}
Now we can use Equation (7) of \cite{Py}: for $k\geq 0$
\begin{eqnarray*}
\P(S^-\leq k) & = & \P(\sup_{z\geq 0} \{\nu_\rho(z) - \mu z\}\leq k)\\
& = & (1-\rho/\mu)\sum_{i=0}^k (-\rho/\mu)^i\frac{(k-i)^i}{i!}e^{\rho(k-i)/\mu}.
\end{eqnarray*}
In particular, $\P(S^- = 0) = 1-\rho/\mu$. Therefore,
\begin{eqnarray*}
\P(S^+\geq S^-) & = & \sum_{k=0}^\infty \P(S^-\leq k)p_+(1-p_+)^k\\
& = & p_+(1-\rho/\mu)\\
&&\qquad\times\,\sum_{k=0}^\infty\sum_{i=0}^k (-\rho/\mu)^i\frac{(k-i)^i}{i!}e^{\rho(k-i)/\mu}(1-p_+)^k\,.
\end{eqnarray*}

\begin{figure}[t]
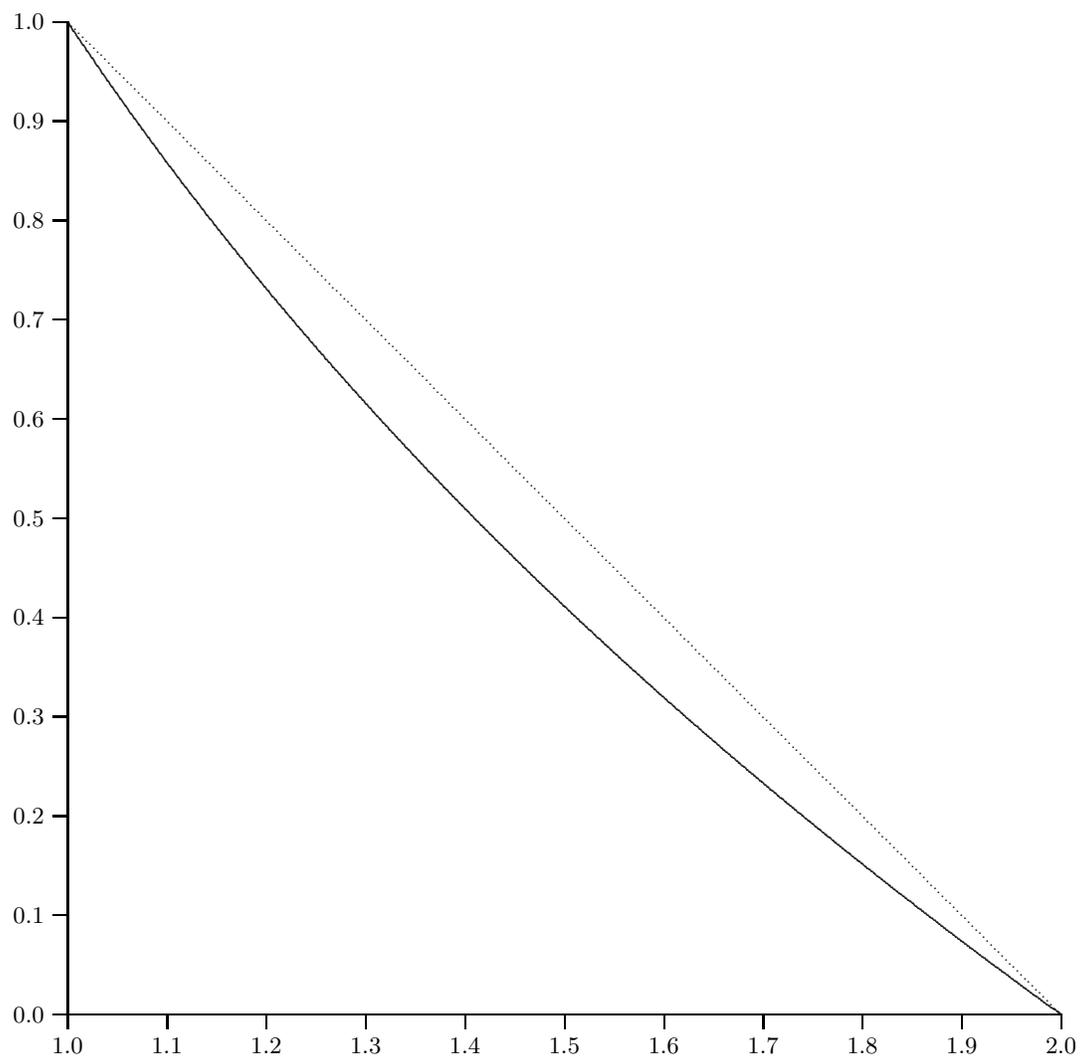

\[
\beginpicture
\footnotesize
\setcoordinatesystem units <0.8\textwidth,0.8\textwidth>
  \setplotarea x from 1 to 2, y from 0 to 1
  \axis bottom
    shiftedto y=0.0
    ticks numbered from 1 to 2 by 0.1
    /
  \axis left
    shiftedto x=1.0
    ticks numbered from 0 to 1 by 0.1
    /
\setquadratic
\plot
1.0000    1.0000     1.0050    0.9925     1.0100    0.9850    1.0150    0.9777    1.0200    0.9704    1.0250    0.9631    1.0300    0.9558    1.0350    0.9486    1.0400    0.9414    1.0450    0.9343    1.0500    0.9272    1.0550    0.9201    1.0600    0.9131    1.0650    0.9062    1.0700    0.8992    1.0750    0.8923    1.0800    0.8855    1.0850    0.8787    1.0900    0.8719    1.0950    0.8652    1.1000    0.8585    1.1050    0.8518    1.1100    0.8451    1.1150    0.8385
    1.1200    0.8320    1.1250    0.8255    1.1300    0.8190    1.1350    0.8125    1.1400    0.8061    1.1450    0.7997
    1.1500    0.7933    1.1550    0.7870    1.1600    0.7807    1.1650    0.7744    1.1700    0.7682    1.1750    0.7620
    1.1800    0.7558    1.1850    0.7496    1.1900    0.7435    1.1950    0.7374    1.2000    0.7314    1.2050    0.7253
    1.2100    0.7193    1.2150    0.7134    1.2200    0.7074    1.2250    0.7015    1.2300    0.6956    1.2350    0.6898
    1.2400    0.6839    1.2450    0.6781    1.2500    0.6723    1.2550    0.6666    1.2600    0.6608    1.2650    0.6551
    1.2700    0.6494    1.2750    0.6438    1.2800    0.6382    1.2850    0.6325    1.2900    0.6270    1.2950    0.6214
    1.3000    0.6159    1.3050    0.6104    1.3100    0.6049    1.3150    0.5994    1.3200    0.5940    1.3250    0.5885
    1.3300    0.5831    1.3350    0.5778    1.3400    0.5724    1.3450    0.5671    1.3500    0.5618    1.3550    0.5565
    1.3600    0.5512    1.3650    0.5460    1.3700    0.5407    1.3750    0.5355    1.3800    0.5303    1.3850    0.5252
    1.3900    0.5200    1.3950    0.5149    1.4000    0.5098    1.4050    0.5047    1.4100    0.4996    1.4150    0.4946
    1.4200    0.4895    1.4250    0.4845    1.4300    0.4795    1.4350    0.4746    1.4400    0.4696    1.4450    0.4647
    1.4500    0.4597    1.4550    0.4548    1.4600    0.4500    1.4650    0.4451    1.4700    0.4402    1.4750    0.4354
    1.4800    0.4306    1.4850    0.4258    1.4900    0.4210    1.4950    0.4162    1.5000    0.4115    1.5050    0.4067
    1.5100    0.4020    1.5150    0.3973    1.5200    0.3926    1.5250    0.3880    1.5300    0.3833    1.5350    0.3787
    1.5400    0.3740    1.5450    0.3694    1.5500    0.3648    1.5550    0.3602    1.5600    0.3557    1.5650    0.3511
    1.5700    0.3466    1.5750    0.3421    1.5800    0.3376    1.5850    0.3331    1.5900    0.3286    1.5950    0.3241
    1.6000    0.3197    1.6050    0.3152    1.6100    0.3108    1.6150    0.3064    1.6200    0.3020    1.6250    0.2976
    1.6300    0.2932    1.6350    0.2889    1.6400    0.2845    1.6450    0.2802    1.6500    0.2759    1.6550    0.2716
    1.6600    0.2673    1.6650    0.2630    1.6700    0.2587    1.6750    0.2545    1.6800    0.2502    1.6850    0.2460
    1.6900    0.2418    1.6950    0.2376    1.7000    0.2334    1.7050    0.2292    1.7100    0.2250    1.7150    0.2208
    1.7200    0.2167    1.7250    0.2125    1.7300    0.2084    1.7350    0.2043    1.7400    0.2002    1.7450    0.1961
    1.7500    0.1920    1.7550    0.1879    1.7600    0.1839    1.7650    0.1798    1.7700    0.1758    1.7750    0.1717
    1.7800    0.1677    1.7850    0.1637    1.7900    0.1597    1.7950    0.1557    1.8000    0.1517    1.8050    0.1478
    1.8100    0.1438    1.8150    0.1398    1.8200    0.1359    1.8250    0.1320    1.8300    0.1280    1.8350    0.1241
    1.8400    0.1202    1.8450    0.1163    1.8500    0.1125    1.8550    0.1086    1.8600    0.1047    1.8650    0.1009
    1.8700    0.0970    1.8750    0.0932    1.8800    0.0893    1.8850    0.0855    1.8900    0.0817    1.8950    0.0779
    1.9000    0.0741    1.9050    0.0703    1.9100    0.0666    1.9150    0.0628    1.9200    0.0590    1.9250    0.0553
    1.9300    0.0515    1.9350    0.0478    1.9400    0.0441    1.9450    0.0404    1.9500    0.0367    1.9550    0.0330
    1.9600    0.0293    1.9650    0.0256    1.9700    0.0219    1.9750    0.0182    1.9800    0.0146    1.9850    0.0109
    1.9900    0.0073    1.9950    0.0040    2.0000         0
/
\setdots <2pt>
\setlinear
\plot
1 1 2 0
/
\endpicture
\]
\caption{Picture of $P(S^+\geq S^-)$, together with uniform distribution.}\label{fig:detlamu}
\end{figure}

We were not able to simplify this formula significantly. For $\lambda=1$ and $\mu=2$, Figure \ref{fig:detlamu} gives the graph of $P(S^+\geq S^-)$ as a function of $\rho\in (1,2)$.
We compare it with the probability 
$$\frac{\mu-\rho}{\mu-\lambda}\,$$ 
which we would get if we would take a Poisson process of intensity $\lambda$ for $x\geq 0$, and of intensity $\mu$ for $x\leq 0$. \\

An interesting limit is $\lambda\to\infty$. This means that $S^+=0$, and we get
\[\P(S^+\geq S_-) = 1-\rho/\mu.\]
This is exactly the same probability as when we would have a Poisson process of intensity $\mu$ to the left of $0$, and no particles to the right of $0$!

\newpage

\chapter{Fluctuation Results}\label{ch:CLT}

\section{Some Elementary Identities}
We have seen that there is only one family of ergodic equilibrium measures for the Hammersley interacting fluid system. Let us denote it by $\{\nu_\lambda\,:\,\lambda>0\}$, where
\begin{equation}\label{eq:intensity}
\lambda:=\E\nu_\lambda(1)=\frac{\gamma}{2}\sqrt{\tan\alpha}\,
\end{equation}
(recall \eqref{eq:shape}). Let
\begin{equation}\label{eq:direction}
\phi:=\left(\frac{\gamma}{2\lambda}\right)^2=\frac{1}{\tan\alpha}\,\,\,\,\mbox{ and }\,\,\,\, \psi:=\frac{\gamma^2}{2\lambda}\,.
\end{equation}
Thus, the speed $\phi$ corresponds to the characteristic angle $\alpha$ associated to the equilibrium measure $\nu_\lambda$. For simple notation, put $L_\lambda:=L_{\nu_\lambda}$,
$$\nu_\lambda(x):=B_\alpha(\0,(x,0))\,,\,\,\,\nu^t_\lambda(x):=B_\alpha((0,t),(x,t))\,,$$
$$\nu^*_\lambda(t):=B_\alpha(\0,(0,t))\,,\,\,\, \nu^{*,x}_\lambda(t):=B_\alpha((x,0),(x,t))\,.$$

\begin{prop}\label{prop:VarCov}
For given  $\lambda>0$ and $x,t\geq0$
\begin{equation}\label{eq:VarCov1}
\Var\, L_\lambda(x,t)=\Var \,\nu^*_\lambda(t) -\Var \,\nu_\lambda(x)+2\Cov \left(\nu_\lambda(x),L_\lambda(x,t)\right)\,.
\end{equation}
In particular,
\begin{equation}\label{eq:VarCov}
\Var\, L_\lambda(\phi t,t)=2\Cov \left(\nu_\lambda(\phi t),L_\lambda(\phi t,t)\right)\,.
\end{equation}
\end{prop}
\noindent{\bf Proof:} By Proposition \ref{prop:Busemann} (additivity),
$$L_\lambda(x,t)=\nu_\lambda(x)+\nu_\lambda^{*,x}(t)\,.$$
Thus,
\begin{eqnarray}
\nonumber \Var \,L_\lambda(x,t) &=& \Var \,\nu_\lambda(x)+\Var\, \nu_\lambda^{*,x}(t)+2\Cov\left(\nu_\lambda(x),\nu_\lambda^{*,x}(t)\right)\\
\nonumber &=& \Var \,\nu_\lambda(x)+\Var\, \nu_\lambda^{*,x}(t)\\
\nonumber&&\qquad\qquad\qquad+2\Cov\left(\nu_\lambda(x),L_\lambda(x,t)-\nu_\lambda(x)\right)\\
\nonumber &=& \Var\, \nu_\lambda^{*,x}(t)- \Var \,\nu_\lambda(x)+2\Cov\left(\nu_\lambda(x),L_\lambda(x,t)\right)\\
\nonumber &=& \Var\, \nu_\lambda^{*}(t)-\Var \,\nu_\lambda(x)+2\Cov\left(\nu_\lambda(x),L_\lambda(x,t)\right)\,.
\end{eqnarray}
(In the last step we also used translation invariance.) \\

By Proposition \ref{prop:sym} (rescaling), with 
$$\tilde\alpha:=5\pi/4\,\mbox{ and }\,r:=\phi^{-1/2}=\sqrt{\tan\alpha}\,,$$
we have that
$$\nu_\lambda(\phi t)=B_\alpha(\0,(\phi t,0))\stackrel{\cal D}{=}B_{5\pi/4}(\0,(\phi^{1/2}t,0))\,$$
and
$$\nu^{*}_\lambda(t)=B_\alpha(\0,(0,t))\stackrel{\cal D}{=}B_{5\pi/4}(\0,\phi^{1/2}t)\,.$$

Again, by Proposition \ref{prop:sym} (reflection),
$$B_{5\pi/4}(\0,(\phi^{1/2}t,0))\stackrel{\cal D}{=}B_{5\pi/4}(\0,\phi^{1/2}t)$$
which shows that
$$\nu_\lambda(\phi t)\stackrel{\cal D}{=}\nu^{*}_\lambda(t)\,.$$
In particular,
$$\Var \,\nu_\lambda(\phi t)= \Var\, \nu_\lambda^{*}(t)$$
and hence,
\begin{equation*}
\Var\, L_\lambda(\phi t,t)=2\Cov\left(\nu_\lambda(\phi t),L_\lambda(\phi t,t)\right)\,.
\end{equation*}

\hfill $\Box$\\

Recalling that
$$L_{\lambda}(x,t):=\sup_{z\leq x} \left\{ \nu_\lambda(z) + L((z,0),(x,t))\right\}\,,$$
one may guess that the correlation between $L_\lambda$ and $\nu_\lambda$ is proportional to the amount collected by $L_\lambda$ in the positive $x$-axis. Which is to say that
\begin{equation}\label{eq:exitconj}
\Cov\left(\nu_\lambda(x),L_\lambda(x,t)\right)\sim \E Z^+_\lambda(x,t)
\end{equation}
Since $\phi=(\tan\alpha)^{-1}$, one may expect that
$$\E Z^+_\lambda(\phi t,t)\sim  t^{\xi}\,,$$
where $\xi\in(0,1)$ is the critical exponent for transversal fluctuations of $\alpha$-rays. This indicates that the variance of $L_\lambda$ along the characteristic speed $\phi$ has sub-linear growth.\\

Another remarkable identity is given below:
\begin{prop}\label{prop:formula}
For $x\in\R$ and $t\geq 0$ let
$$b_{x,t}:=x-\phi t\,.$$
Then
\begin{equation}\label{eq:formula}
\E\big(\left\{L_\lambda(x,t)-\left[\nu_\lambda(b_{x,t})+\psi t\right]\right\}^2\big)=\Var\big(L_\lambda\left(\phi t,t\right)\big)\,.
\end{equation}
\end{prop}
\noindent{\bf Proof:} By Proposition \ref{prop:Busemann} (additivity),
$$L_\lambda(x,t)=\nu_\lambda(b_{x,t})+ B_\alpha((b_{x,t},0),(x,t))\,.$$
On the other hand (translation invariance),
$$B_\alpha((b_{x,t},0),(x,t))\stackrel{\cal D}{=}B_\alpha(\0,(x-b_{x,t},t))=B_\alpha(\0,(\phi t ,t))=L_\lambda(\phi t,t)\,.$$
Since
$$\E L_\lambda(\phi t, t)=\psi t\,,$$
we get that
$$L_\lambda(x,t)-\left[\nu_\lambda(b_{x,t})+\psi t\right]\stackrel{\cal D}{=}L_\lambda(\phi t ,t))-\psi t\,,$$
which shows \eqref{eq:formula}.

\hfill $\Box$\\

To illustrate the importance of \eqref{eq:formula}, let us assume that \eqref{eq:exitconj} holds and that $z_t$ is a path with speed $a\neq\phi$, and let $b_t=z_t-\phi t$. By \eqref{eq:formula}, we get that
$$\lim_{t\to\infty}\frac{L_\lambda(z_t,t)-\left[\nu_\lambda(b_{t})+\psi t\right]}{\sqrt{t}}=0\,,$$
in the $\LL^2$ sense. Therefore, if one has Gaussian fluctuations for the equilibrium measures (which is a plausible assumption) then $L_\lambda$ will have the same behavior along any  speed $a\neq\phi$. This will be made precise in the classical Hammersley model.

\section{The Exit Point Formula}

\begin{thm}\label{thm:exitpoint}
Consider the classical Hammersley model. Then
$$\Cov\left(\nu_\lambda(x),L_\lambda(x,t)\right)=\lambda \E \left(Z^+_\lambda(x,t)\right)\,.$$
In particular,
$$\Var\, L_\lambda(x,t)=\frac{t}{\lambda} - \lambda x+2\lambda \E \left(Z^+_\lambda(x,t)\right) \,.$$
\end{thm}
\noindent{\bf Proof:} For each $\epsilon>0$ we couple $\nu^+_\lambda$ with $\nu_{\lambda+\epsilon}$ and let the negative part stay the same. Denote this new measure by $\nu_{\lambda,\epsilon}$. Let
$$f(\epsilon):=\E \left(L_{\nu_{\lambda,\epsilon}}(x,t)\right)\,,$$
and, for each $n\geq 0$ let (disregarding the dependence on $x$ and $t$ for a while)
$$a_n:=\E (L_{\nu_{\lambda,\epsilon}}(x,t)\mid \nu_{\lambda,\epsilon}(x)=n)\,.$$
Note that $a_n$ does, in fact, not depend on $\epsilon$ since we are condition a Poisson process  on the interval $[0,x]$. Then
$$f(\epsilon)=\sum_{n=0}^{\infty}\frac{x(\lambda+\epsilon)^n}{n!}e^{-x(\lambda+\epsilon)} a_n$$
and
\begin{eqnarray}
\nonumber\frac{d}{d\epsilon} f(0)&=&\frac{1}{\lambda} \sum_{n=0}^{\infty}\frac{(x\lambda)^n}{n!}e^{-x\lambda} a_n - x\sum_{n=0}^{\infty}\frac{(x\lambda)^n}{n!}e^{-x\lambda} a_n\\
\nonumber&=& \frac{1}{\lambda}\E \left(L_{\nu_\lambda} \nu_\lambda\right)-x\E \left(L_{\nu_\lambda}\right)\,,
\end{eqnarray}
which shows that
\begin{eqnarray}
\nonumber\lambda\frac{d}{d\epsilon} f(0)&=&\E \left(L_{\nu_\lambda} \nu_\lambda\right)-\lambda x\E \left(L_{\nu_\lambda}\right)\\
\nonumber&=& \E \left(L_{\nu_\lambda} \nu_\lambda\right)- \E\left(\nu_\lambda\right)\E \left(L_{\nu_\lambda}\right)\\
\nonumber&=&\Cov(\nu_\lambda,L_\lambda)\,.
\end{eqnarray}

Now we are going to calculate this derivative in a different way. Note that the probability of finding more than one Poisson point of $\nu_{\lambda,\epsilon}-\nu_\lambda$ lying within $[0,x]$ is $o(\epsilon)$. If we have exactly one Poisson point $z\in[0,x]$, then this point will contribute to $L_{\nu_{\lambda,\epsilon}}$ if and only if $z\in[0,Z_\lambda(x,t)]$. So
\begin{eqnarray}
\nonumber f(\epsilon)&=&\E\left(L_{\nu_{\lambda,\epsilon}}(x,t)\right)\\
\nonumber&=&\E \left(L_\lambda(x,t)\right)+\epsilon\int_0^x\P\left(Z_\lambda(x,t)>z\right)dz + o(\epsilon)\\
\nonumber&=&f(0)+\epsilon\int_0^x\P\left(Z_\lambda(x,t)>z\right)dz + o(\epsilon)\,.
\end{eqnarray}
Thus,
$$\frac{d}{d\epsilon} f(0)=\int_0^x\P\left(Z_\lambda(x,t)>z\right)dz=\E \left(Z^+_\lambda(x,t)\right)\,,$$
which shows the first part of the proposition. The second part follows by combining the first part together with Proposition \ref{prop:VarCov}.

\hfill $\Box$\\

\begin{coro}\label{coro:crossformula}
Let $\alpha\in(\pi,3\pi/2)$ and for $x,t\geq 0$ define
$$B_\alpha(x,t):=B_\alpha(\0,(x,t))\,.$$
Then
$$\Var B_\alpha(x,t)=\frac{t}{\sqrt{\tan\alpha}}-(\sqrt{\tan\alpha}) x + 2(\sqrt{\tan\alpha})\E Z_\alpha^+(x,t)\,.$$
In particular,
$$\Var B_\alpha\left(t,(\tan\alpha)t\right)= 2(\sqrt{\tan\alpha})\E Z_\alpha^+(t,(\tan\alpha)t)\,.$$
\end{coro}

\noindent{\bf Proof:} It follows from Theorem \ref{thm:exitpoint} and Corollary \ref{thm:classicalPoisson}.

\hfill $\Box$\\

\section{Gaussian Fluctuations}
To illustrate the importance of \eqref{eq:formula}, let us restrict ourselves to the classical Hammersley model. Thus,
$$\phi:=\frac{1}{\lambda^2}\,\,\mbox{ and }\,\, \psi:=\frac{2}{\lambda}\,.$$
\begin{thm}\label{thm:dependence}
Consider the classical Hammersley model. Given a deterministic path $(z_t)_{t\geq0}$, let 
$$b_{t}=z_t-\lambda^{-2} t\,.$$ 
Then
$$\lim_{t\to\infty}\frac{\E\big(\left\{L_\lambda(z_t,t)-\left[\nu_\lambda(b_{t})+2\lambda^{-1} t\right]\right\}^2\big)}{t}=0\,.$$
\end{thm}

\noindent{\bf Proof:} Notice that
$$0\leq t^{-1}Z^+_\lambda(\lambda^{-2}t,t)\leq \lambda^{-2}\,.$$
Together with the Dominated Convergence Theorem, Theorem \ref{thm:genexitcontrol} implies that
$$\lim_{t\to\infty}\frac{\E \left(Z^+_\lambda(\lambda^{-2} t,t)\right)}{t}=0\,.$$
By Proposition  \ref{prop:formula} and Theorem \ref{thm:exitpoint}, this proves Theorem \ref{thm:dependence}. 

\hfill $\Box$\\

\begin{coro}\label{coro:clt}
Consider the classical Hammersley model. Assume that $(z_t)_{t\geq0}$ is a deterministic path such that
$$\lim_{t\to\infty}\frac{z_t}{t}=a\,.$$
Then
\begin{equation}\label{eq:var}
\lim_{t\to\infty}\frac{\Var \big(L_\lambda(z_t,t)\big)}{t}=\sigma^2:=|a\lambda-\frac{1}{\lambda}|\,.
\end{equation}
Furthermore, if $a\neq \lambda^{-2}$ then
\begin{equation}\label{eq:clt}
\lim_{t\to\infty}\P\left(\frac{L_\lambda(z_t,t)-[ \lambda z_t+\lambda^{-1}t]}{\sigma t^{1/2}} \leq u\right)=\P(N\leq u)\,,
\end{equation}
where $N$ is a standard Gaussian random variable.
\end{coro}

\noindent{\bf Proof:} Theorem \ref{thm:dependence} shows that
$$\lim_{t\to\infty}\frac{L_\lambda(z_t,t)-\left[\nu_\lambda(b_{t})+2\lambda^{-1}t\right]}{t^{1/2}}=0\,,$$
in the $\LL^2$ sense. Since $\nu_\lambda$ is a one-dimensional Poisson process of intensity $\lambda$, this implies \eqref{eq:var} and \eqref{eq:clt}.

\hfill$\Box$\\

It is also interesting to remark that Proposition \ref{prop:formula} implies that
\begin{equation}\label{eq:errorstat}
 \frac{L_1(z_t,t)-[\nu_1(b_t)+2t]}{t^{1/3}}\,\,\stackrel{\cal D}{=}\,\,\frac{L_1(t,t)-2t}{t^{1/3}}\,.
\end{equation}
Baik and Rains \cite{BaRa} proved that the right hand side of \eqref{eq:errorstat} converges in distribution to $F_0$, the zero-mean Tracy-Widon distribution. Therefore, the left hand side of  \eqref{eq:errorstat} has the same limit behavior.\\

These results naturally lead us to a central limit theorem for the Busemann function. Without lost of generality, let take $\alpha=5\pi/4$ and let
$$B(\beta,t):=B_\alpha(\0,\vec\beta)\,,$$
where $\vec\beta:=(\cos\beta,\sin\beta)$. Notice that
$$\E B(\beta,t)=(\cos\beta+\sin\beta)t\,.$$
If $\beta\in[0,\pi/2]$ then
\begin{equation}\label{eq:cltBuse1}
B(\beta,t)=L_1(t\vec\beta)\,
\end{equation}
(as process in $t$). This is just a consequence of \eqref{eq:LastBuse}. On the other hand, if $\beta\in[\pi/2,\pi]$ then
\begin{equation}\label{eq:cltBuse2}
B(\beta,t)=Y(t)-X(t)\,,
\end{equation}
where $Y(t)$ and $X(t)$ are two independent Poisson process with intensity $1$. \\

\begin{figure}[!h]
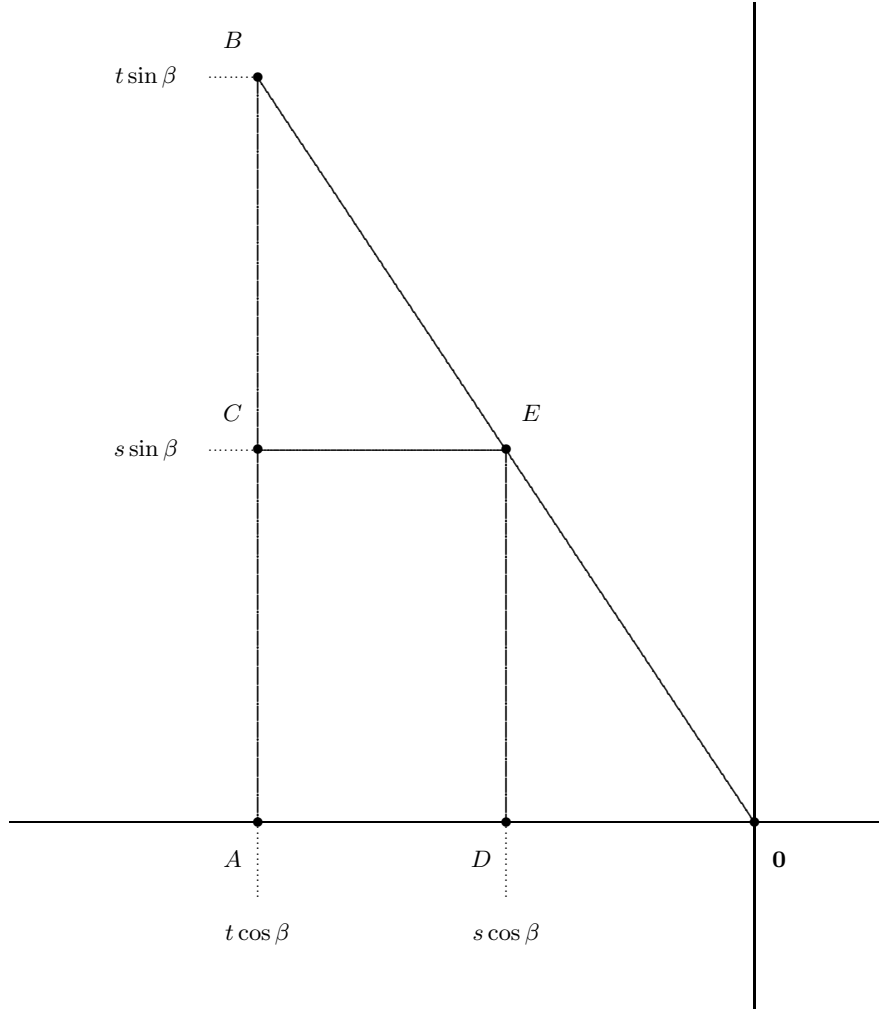

\[
\beginpicture
\footnotesize
\setcoordinatesystem units <0.2\textwidth,0.3\textwidth>
  \setplotarea x from -3 to .5, y from -0.5 to 2.2
  \axis bottom
    shiftedto y=0.0
    /
  \axis left
    shiftedto x=0.0
    /
\multiput{$\bullet$} at
0 0
-1 0
-2 0
-2 1
-2 2
-1 1
/
\put {$\0$} at 0.1 -0.1
\put {$A$} at -2.1 -0.1
\put {$B$} at -2.1 2.1
\put {$C$} at -2.1 1.1
\put {$D$} at -1.1 -0.1
\put {$E$} at -0.9 1.1

\put {$s\cos\beta$} at -1 -0.3
\put {$t\cos\beta$} at -2 -0.3
\put {$s\sin\beta$} at -2.45 1
\put {$t\sin\beta$} at -2.45 2

\setlinear
\plot
0 0 -2 0 -2 2 0 0
/
\plot
-2 1 -1 1 -1 0
/
\setdots <2pt>
\plot
-1 0 -1 -0.2
/
\plot -2 0 -2 -0.2
/
\plot -2 1 -2.2 1
/
\plot -2 2 -2.2 2
/

\endpicture
\]
\caption{Independent increments for $L_{1}$.}\label{fig:ABCDE}
\end{figure}

To obtain \eqref{eq:cltBuse2}, note that, by additivity and anti-symmetry (Proposition \ref{prop:Busemann}),
\begin{eqnarray*}
B(\vec\beta,t)&=&B_\alpha\big(\0,(t\cos\beta,0)\big)+B_\alpha\big((t\cos\beta,0),t\vec\beta\big)\\
&=&B_\alpha\big((t\cos\beta,0),t\vec\beta\big)-B_\alpha\big((t\cos\beta,0),\0\big)\\
& = & L_{1}((t\cos\beta,0), t\vec\beta) - L_{1}((t\cos\beta,0),\0).
\end{eqnarray*}
By Theorem \ref{thm:Burke} (Burke's Theorem), both $t\mapsto L_{1}(\0,(t,0))$ and $t\mapsto L_{1}(\0,(0,t))$ are Poisson processes with intensity $1$, and $t\mapsto L_{\nu_\alpha}(\0,t\vec\beta)$ has independent increments. This last statement can be seen in the following way: consider the increments $L_{1}(\0,s\vec\beta)$ and $L_{1}(s\vec\beta,t\vec\beta)$, for $0\leq s\leq t$; see Figure \ref{fig:ABCDE}.
Then both increments depend on the independent Poisson processes $L_{1}$ restricted to (the line segments) $\0 A$ and $AB$, and on the Poisson process in the triangle $\0 AB$. Burke's Theorem for $L_{1}$ now implies that $L_{1}$ restricted to $CE$ is independent of $L_{1}$ restricted to $ED$. Since $L_{1}$ restricted to $BE$ depends on $L_{1}$ on $CE$, $CB$ and the Poisson process in $ECB$, whereas $L_{1}$ restricted to $\0 E$ depends on $L_{1}$ on $\0 D$, $DE$ and the Poisson process in $\0 DE$, and since these six processes are all independent, we proved that the increments are independent.

\begin{coro}\label{coro:TCLBus}
For $\beta\in[0,\pi]$
$$\lim_{t\to\infty}\frac{\Var B\big(\beta,t\big)}{t}=|\cos\beta-\sin\beta|\,.$$
Further, for $\beta\in[0,\pi/4)$
$$\lim_{t\to\infty}\P\left(\frac{B\big(\beta,t\big)-(\cos\beta+\sin\beta)t }{(\sqrt{\cos\beta-\sin\beta})t^{1/2}}\leq u\right)= \P(N\leq u)\,,$$
and for $\beta\in(\pi/4,\pi]$
$$\lim_{t\to\infty}\P\left(\frac{B\big(\beta,t\big)-(\cos\beta+\sin\beta)t }{(\sqrt{\sin\beta-\cos\beta})t^{1/2}}\leq u\right)= \P(N\leq u)\,,$$
where $N$ denotes a standard Gaussian random variable.
\end{coro}
\noindent{\bf Proof:} It follows by combining \eqref{eq:cltBuse1} and \eqref{eq:cltBuse2} together with Corollary \ref{coro:clt}.

\hfill $\Box$\\

\section{Cube-Root Asymptotics for $L_\la$}

In this section we will show how, in the classical Hammersley process, we can use the results we have obtained in the previous sections to prove cube-root fluctuations of $L_\lambda$ and $L$.\\

 In \cite{BaDeJo} results on the limiting distribution of $L$ were derived, and building on these methods, Baik and Rains \cite{BaRa} derived limiting distributions for $L_\lambda$, using exact formulas for this length, and heavy analytic methods to analyze asymptotics of these exact formulas. We were not able to obtain these limiting distributions using our methods, but we were able to find the correct scaling. Furthermore, we were able to link these results to transversal fluctuations of the longest path. Also, our methods have since been extended to other models, where the analytic methods were not (yet) available.

\subsection{Upper Bounds}
We start by bounding the fluctuations of $L_1(t,t)$. We wish to use Chebyshev's inequality, so the following lemma is useful to us:
\begin{lem}\label{lem:Varla}
For all $\la>0$, we have
\[ \Var\big(L_\la(t,t)\big) \leq t|\la - \frac{1}{\la}| + \Var\big(L_1(t,t)\big).\]
\end{lem}
{\bf Proof:} Using reflection in the diagonal, we can see that
\[ L_\la(t,t) \stackrel{\cal D}{=} L_{1/\la}(t,t).\]
This means that it is enough to prove the lemma for $0<\la<1$. \\

Now couple $\nu_\la$ and $\nu_1$, so $\nu_1$ contains more points. Suppose for $x\in\R$ and $t>0$ that $Z_\la(x,t)> Z_1(x,t)$. Then
\begin{eqnarray*}
&&L((Z_1(x,t),0),(x,t)) - L((Z_\la(x,t),0),(x,t))\\
&&\qquad \qquad \qquad\qquad\qquad\qquad>\nu_1(Z_\la(x,t)) - \nu_1(Z_1(x,t))\\
&&\qquad \qquad \qquad\qquad\qquad\qquad\geq \nu_\la(Z_\la(x,t)) - \nu_\la(Z_1(x,t))\,.
\end{eqnarray*}
This would imply that
\[ L((Z_1(x,t),0),(x,t)) + \nu_\la(Z_1(x,t))> L_\la(x,t),\]
which is absurd. Therefore, for all $x\in\R$ and $t>0$, we have that $Z_\la(x,t)\leq Z_1(x,t)$. Using this and Theorem \ref{thm:exitpoint}, we get
\begin{eqnarray*}
\Var\big(L_\la(t,t)\big) & = & t(\frac1{\la} - \la) + {2}{\la}\E\big(Z^+_{\la}(t,t)\big)\\
& \leq & t(\frac1{\la} - \la) + 2\E\big(Z^+_{1}(t,t)\big)\\
& = & t(\frac1{\la} - \la) + \Var\big(L_1(t,t)\big)
\end{eqnarray*}

\hfill$\Box$\\

\begin{thm}\label{thm:main}
There exists a constant $c_1>0$ such that for all $t>0$ and all $r>0$, we have
\[ \P \big(Z_1(t,t)>r\,\E Z^+_1(t,t)\big) \leq \frac{c_1t^2}{(\E Z^+_1(t,t))^3}\left(\frac1{r^3}+\frac1{r^4}\right).\]
\end{thm}

\noindent{\bf Proof:} Note that for any $\la \geq 1$ and $0<u<t$,
\begin{eqnarray*}
\P \big(Z_1(t,t)>u\big) & = & \P \Big(\exists\ z>u : \nu_1(z) + L((z,0),(t,t)) = L_1(t,t)\Big)\\
& \leq & \P \Big(\exists\ z>u: \nu_1(z) + L_\la(t,t)-L_\la(z,0)\\
&& \qquad \qquad \qquad\qquad \qquad\qquad \qquad\geq L_1(t,t)\Big)\\
& = & \P \Big(\exists\ z>u: \nu_1(z) - \nu_\la (z) \\
&& \qquad\qquad \qquad\qquad \qquad\geq L_1(t,t) - L_\la (t,t)\Big).
\end{eqnarray*}

For the inequality we used the fact that for any $\la>0$ and any $\0\leq \xx\leq \yy \in \R^2$, we have
$$L_\la(\yy) \geq L_\la(\xx)+L(\xx,\yy)\,,$$
and hence,
$$L_\la(\yy) - L_\la(\xx)\geq L(\xx,\yy)\,.$$

Note that this is true for any choice of $\nu_\la$! Since we took $\la> 1$, we can choose $\nu_\la$ as a thickening of $\nu_1$, so that 
$$\tilde{\nu}_{\la-1}(z):=\nu_\la (z)-\nu_1(z)\,$$
is in itself a Poisson process with intensity $\la -1$. This means that
\[ \P \big(Z_1(t,t)>u\big) \leq \P \big(\tilde{\nu}_{\la -1}(u) \leq L_\la(t,t)-L_1(t,t)\big).\]

To have a useful bound for all $0\leq u\leq \frac34 t$, we choose $\la$ such that
\[ \E\tilde{\nu}_{\la -1}(u) - \E\big(L_\la(t,t)-L_1(t,t)\big)=(\la -1)u -t\left(\la + \frac1{\la}
-2\right)\]
is maximal. This means that we choose
\[ \la_u = (1-u/t)^{-1/2}. \]

Some useful elementary inequalities, that hold for all $0<u \leq \frac34 t$, are
\begin{equation}\label{eq:Ebig}
\begin{array}{l}
\la_u  \leq  2   \\
\vspace{-0.2cm}\\
\E\tilde{\nu}_{\la_u -1}(u) - \E\big(L_{\la_u}(t,t)-L_1(t,t)\big)  \geq  \frac{1}{4}u^2/t \\
\vspace{-0.2cm}\\
\la_u-1/{\la_u} \leq  2u/t.
\end{array}
\end{equation}
Using Theorem \ref{thm:exitpoint}, Lemma \ref{lem:Varla} and \eqref{eq:Ebig}, we get
\begin{eqnarray*}
\Var\big(L_{\la_u}(t,t)-L_1(t,t)\big) & \leq & 2\left\{\Var\big(L_{\la_u} (t,t)\big)+\Var\big(L_1(t,t)\big)\right\}\\
& \leq & 8\E Z^+_1(t,t) + 2t({\la_u} - 1/{\la_u})\\
& \leq & 8\E Z^+_1(t,t) + 4u.
\end{eqnarray*}

Now we can use Chebyshev's inequality:
\begin{eqnarray}\label{eq:Z>u}
\P \big(Z_1(t,t)>u\big) & \le & \P \big(\tilde{\nu}_{{\la_u} -1}(u) \leq L_{\la_u}(t,t)-L_1(t,t)\big) \nonumber \\
& \leq & \P\left(\tilde{\nu}_{{\la_u} -1}(u) \leq \E \tilde{\nu}_{{\la_u} -1}(u) -
u^2/(8t)\right)\nonumber\\
&+& \P\left(L_{\la_u}(t,t)-L_1(t,t) \geq \E \tilde{\nu}_{{\la_u} -1}(u) - u^2/(8t)\right) \nonumber \\
\nonumber& \leq & \frac{64t^2(\la_u - 1)u}{u^4} \\
&+& \P\Big(L_{\la_u}(t,t)-L_1(t,t)\nonumber\\
&&\qquad\qquad\qquad \geq \E\left(L_{\la_u}(t,t)-L_1(t,t)\right) + \frac{u^2}{8t}\Big) \nonumber \\
& \leq & \frac{64t^2}{u^3} + \frac{64t^2(8\E Z^+_1(t,t) + 4u)}{u^4} \nonumber \\
\nonumber& \leq & 512\left(\frac{t^2}{u^3} + \frac{t^2\E Z^+_1(t,t)}{u^4}\right)\,.
\end{eqnarray}

If $t\ge u\geq \frac34 t$, we see that
\begin{equation*}
\P \big(Z_1(t,t)>u\big)  \leq  \P \big(Z_1(t,t)>\frac 34 t\big)
\leq \frac{c_1t^2}{u^3} + \frac{c_1t^2\E Z^+_1(t,t)_+}{u^4},
\end{equation*}
where we use \eqref{eq:Z>u} and choose $c_1=(4/3)^4\times 512$. This means that \eqref{eq:Z>u} is true for all $0< u \leq t$, if we replace $512$ by $c_1$. Since $Z_1(t,t)\leq t$ with probability one, we have that \eqref{eq:Z>u} is true for all $u>0$. The theorem now follows from choosing $u=r\, \E Z^+_1(t,t)$.

\hfill $\Box$\\

\begin{coro}
\label{cor:var_expectation}
\[\limsup_{t\to \infty} \frac{\Var
(L_1(t,t))}{t^{2/3}}=2\limsup_{t\to \infty} \frac{\E Z^+_1(t,t)}{t^{2/3}} < +\infty\,.\]
\end{coro}

\noindent{\bf Proof:}
By Theorem \ref{thm:main},
\begin{eqnarray*}
1&=&\int_0^\infty\P \big(Z^+_1(t,t)>r\, \E Z^+_1(t,t)\big)dr\\
&\leq&\int_0^{1/2}\P \big(Z^+_1(t,t)>r\, \E Z^+_1(t,t)\big)dr\\
&&\qquad\qquad\qquad\qquad+\,\int_{1/2}^\infty\P \big(Z^+_1(t,t)>r\, \E Z^+_1(t,t)\big)dr\\
&\leq&\frac12+ c_1\left\{ \int_{1/2}^{\infty}\left(\frac1{r^3}+\frac1{r^4}\right)dr\right\}\frac{t^2}{(\E Z^+_1(t,t))^3}\,,
\end{eqnarray*}
and hence,
$$\frac{\E Z^+_1(t,t)}{t^{2/3}}\leq \left(2c_1\left\{ \int_{1/2}^{\infty}\left(\frac1{r^3}+\frac1{r^4}\right)dr\right\} \right)^{1/3} \,.$$ 
By Theorem \ref{thm:exitpoint}, this shows Corollary \ref{cor:var_expectation}.

\hfill $\Box$\\

As a consequence, we get:
\begin{coro}\label{coro:exitail}
\label{cor:main}
There exists a constant $c_2>0$ such that for all $r>0$,
\[ \P \big(Z_1(t,t)>r\,t^{2/3}\big)\leq\frac{c_2}{r^3}.\]
\end{coro}

\noindent{\bf Proof:} This follows from Theorem \ref{thm:main} and the previous corollary.

\hfill $\Box$\\

\subsection{Lower Bounds}

Using the upper bound results from the previous section, we can also get a lower bound result of the same order. We will use the sources and sinks representation of the Last Passage Percolation model: for $x,t\geq 0$ let $\bar{\nu}_\la$ be defined on $[-t,x]$ by
\[\bar\nu_\la(z)=\left\{\begin{array}{ll} \nu_\la(x) & \mbox{for } z\in[0,x]\\
\nu^*_\la(t) & \mbox{for } z\in(0,t]\,.\end{array}\right.\]
Remember that $\nu^*_\la$ is the Poisson process of ``sinks'' on the $t$-axis, with intensity $1/\la$. Thus $\bar\nu_\la$, restricted to $[0,x]$, is a Poisson process with intensity $\la$ (sources) while $\bar\nu_\la$ restricted to $(0,t]$ is a Poisson process with intensity $1/\la$ (sinks).\\

 Denote $\zz:=(z,0)$ for $z\in[0,x]$ and $\zz:=(0,z)$ for $z\in[-t,0)$. By \eqref{eq:fluxformula} we have that
$$L_\la(x,t)=\sup_{z\in[-t,x]}\left\{\bar\nu_\la(z)+L(\zz,(x,t))\right\}\,.$$
Define
\begin{equation}\label{eq:exitsup}
\bar Z_\la(x,t):=\sup\left\{z\in[-t,x]\,:\,\bar\nu_\la(z)+L(\zz,(x,t)=L_\la(x,t)\right\}\,,
\end{equation}
and
\begin{equation}\label{eq:exitinf}
\bar Z'_\la(x,t):=\inf\left\{z\in[-t,x]\,:\,\bar\nu_\la(z)+L(\zz,(x,t)=L_\la(x,t)\right\}\,.
\end{equation}
Thus, $\bar Z_\la$ and $\bar Z'_\la$ are exit points for $\bar\nu_\la$. \\

By definition, and using the invariance under reflection, we can see that
\begin{equation}\label{eq:exitsym}
\bar Z'_\la(x,t)\leq \bar Z_\la(x,t)\,\,\mbox{ and }\,\,\bar Z_\la(x,t)\stackrel{\cal D}{=}-\bar Z'_{1/\la}(t,x)\,.
\end{equation}
It also clear that 
$$\bar Z^+_\la(x,t)= Z^+_\la(x,t)\,.$$ 

By using the same reasoning as in the proof of  Proposition \ref{prop:Znuexists}, the following local comparison property holds. For any $x_2\geq x_1\geq x$, if $ \bar Z_\la(x,t)\geq 0$ then
\begin{equation*}
L(\0,(x_2,t)) - L(\0,(x_1,t)) \leq L_\la(x_2,t) - L_\la(x_1,t)\,.
\end{equation*}
Furthermore, for any $x_1\leq x_2\leq x$ , if $\bar Z'_\la(x,t)\leq 0$ then
\begin{equation*}
L(\0,(x_2,t)) - L(\0,(x_1,t)) \geq L_\la(x_2,t) - L_\la(x_1,t).
\end{equation*}
From now on we denote
$$\Delta(z,t):=L(\0,(t,t))-L(\0,(t-z,t))\,,$$
so that, if  $ \bar Z_\la(t,t)\geq 0$ then
\begin{equation}\label{eq:upboundLLla}
\Delta(z,t)) \leq L_\la(t,t) - L_\la(t-z,t)\,,
\end{equation}
and  if $\bar Z'_\la(x,t)\leq 0$ then
\begin{equation}\label{eq:lowboundLLla}
\Delta(z,t) \geq L_\la(t,t) - L_\la(t-z,t)\,.
\end{equation}

\begin{lem}\label{lem:lambda}
For $r>0$ let 
$$\la_+=1+ \frac{r}{t^{1/3}}>1\,\mbox{ and }\,\la_-=1-\frac{r}{t^{1/3}}<1\,.$$
Then there exist  constants $c_1,c_2>0$ such that
$$\P\left(\bar Z_{\la_+}(t,t)< 0\right)\leq \frac{c_1}{r^3}\,,$$
and
$$\P\left(\bar Z'_{\la_-}(t,t)>0\right)\leq \frac{c_2}{r^3}\,,$$
for all sufficiently large $t$ and $r\in[1,t^{1/3})$.
\end{lem}

\noindent{\bf Proof:}
By rescaling and translation invariance, we see that
\[ \bar Z_\la(t,t) \stackrel{\cal D}{=} \frac1{\la}\bar Z_1(\la t,t/\la)\,\,\mbox{ and }\,\, \bar Z_\la(x+h,t)\stackrel{\cal D}{=}\bar Z_\la(x,t)+h\,.\]
So, by Corollary \ref{cor:main},
\begin{eqnarray*}
\P\left(\bar Z_{\la_+}(t,t)< 0\right) & = & \P\left(\bar Z_1(\la_+ t,t/\la_+)< 0\right) \\
& \leq & \P\left(\bar Z'_1(\la_+ t,t/\la_+) < 0\right)\\
& = & \P\left(\bar Z_1(t/\la_+,\la_+ t)> 0\right)\\
& = & \P\left(\bar Z_1(\la_+ t,\la_+ t)> (\la_+ - 1/\la_+)t\right)\\
& = & \P\left(\bar Z_1(\la_+ t,\la_+ t)> 2rt^{2/3} + O(t^{1/3})\right)\\
& = & O(r^{-3}),
\end{eqnarray*}
for large enough $t$ and $r\in[1,t^{1/3})$. Analogously,
\begin{eqnarray*}
\P\left(\bar Z'_{\la_-}(t,t)>0\right) & = & \P\left(\bar Z'_1(\la_- t,t/\la_-)>0\right)\\
& \leq & \P\left(\bar Z_1(\la_- t,t/\la_-)>0\right)\\
& = & \P\left(\bar Z_1(t/\la_-,t/\la_-)>t(1/\la_- - \la_-)\right)\\
& = & O(r^{-3}),
\end{eqnarray*}
for large enough $t$ and $r\in[1,t^{1/3})$.

\hfill$\Box$\\

\begin{lem}\label{lem:lower2}
Fix $\delta>0$. Then
\begin{equation} \label{eq:lower2}
\limsup_{t\to\infty}\P\left(\sup_{z\in[0,\eps t^{2/3}]} \big\{\nu_1(z) -\Delta(z,t)\big\} > \delta t^{1/3}\right) = O(\epsilon^3)\,.
\end{equation}
\end{lem}

\noindent{\bf Proof:} Let $\la = 1 - rt^{-1/3}$, for some $r>0$. We introduce a Poisson process $\nu_\la$ of sources on the $x$-axis, independent of $\nu_1$, with intensity $\la$ and denote
\[ \nu_{\la}^t(z) = L_{\la}(t,t) - L_{\la}(t-z,t).\]

Clearly, $\nu_\la^t$ is also a Poisson process of intensity $\la$, independent of $\nu_1$. Using  \eqref{eq:lowboundLLla}, we get that the probability in the left hand side of \eqref{eq:lower2} is bounded by
$$\P\left(\sup_{0\leq z \leq \eps t^{2/3}} \big\{\nu_1(z) - \nu_{\la}^t(z)\big\} > \delta t^{1/3}\right) + \P\left(\bar Z'_{\la}(t,t)>0\right)\,,$$
which is equal to
\begin{equation}\label{eq:lower21}
\P\left(\sup_{0\leq z \leq \eps t^{2/3}} \big\{\nu_1(z) - \nu_{\la}(z)\big\} > \delta t^{1/3}\right) + \P\left(\bar Z'_{\la}(t,t)>0\right)\,.
\end{equation}

We will treat these two terms separately. We start with the second term. By Lemma \ref{lem:lambda}
$$\P\left(\bar Z'_{\la}(t,t)>0\right)  =O(r^{-3})$$
for large enough $t$ and $r\in[1,t^{1/3})$.

The first term in \eqref{eq:lower21} concerns the hitting time for the difference of two independent Poisson processes. The process 
$$\{z\mapsto t^{-1/3}\big[\nu_1(z) - \nu_{\la}(z)\big]\}$$ 
converges, as $t\to \infty$ (in the topology of uniform convergence on compacta), to the drifting Brownian motion process, namely
\[ W_{r}(z) := W(2z)+rz\,,\]
where $W$ is the standard Brownian motion on $\R^+$. We now get, for $r<\delta/\eps$,
 \begin{eqnarray*}
 \P\left(\sup_{z\in[0,\eps]} W_{r}(z) >\delta \right)&\leq&\P\left(\sup_{z\in[0,2\eps]} W(z) >\delta-r\epsilon \right)\\
 &=&\P\left(\sup_{z\in[0,1]} \frac{W}{(2\eps z)}{\sqrt{2\eps}} >\frac{\delta-r\epsilon}{\sqrt{2\eps}} \right)\\
 &=&\P\left(\sup_{z\in[0,1]} W( z) >\frac{\delta-r\epsilon}{\sqrt{2\eps}} \right)\\
&=&\sqrt{\frac{2}{\pi}}\int_{\frac{\delta-r\epsilon}{\sqrt{2\eps}}}^\infty e^{-u^2/2}du\,.
\end{eqnarray*}
Taking $r=\delta/(2\eps)$, we get
 \begin{eqnarray*}
 \sqrt{\frac{2}{\pi}}\int_{\frac{\delta-r\epsilon}{\sqrt{2\eps}}}^\infty e^{-u^2/2}du&=&2\int_{\frac{\delta}{\sqrt{8\eps}}}^\infty \frac{e^{-u^2/2}}{\sqrt{2\pi}}du\\
 &\sim&\frac{2\sqrt{8\eps}e^{-\delta^2/(16\eps)}}{\delta\sqrt{2\pi}} \,,
\end{eqnarray*}
using Mills' ratio approximation for the tail of a normal distribution in the last step. Using Donsker's Theorem, this means that, with this choice of $r$, our estimate for the second term in \eqref{eq:lower21} is dominant.

\hfill$\Box$\\

\begin{lem}\label{lem:lower1}
Let $\eta>0$. Then there exist $\delta=\delta(\eta)>0$ and $\epsilon=\epsilon(\delta)>0$ such that
\begin{equation} \label{eq:lower1}
\limsup_{t\to\infty}\P\left(\sup_{z>\eps t^{2/3}} \big\{\nu_1(z) -\Delta(z,t)\big\} < \delta t^{1/3}\right)< \eta\,.
\end{equation}
\end{lem}

\noindent{\bf Proof:} The proof of  this lemma is similar to the previous one, but now we have a different choice of $\la$ and $r$. As before, we introduce a Poisson process $\nu_\la$ of sources on the $x$-axis, independent of $\nu_1$, with intensity $\la=1 + rt^{-1/3}$ and denote
\[ \nu_\la^t(z) = L_\la(t,t) - L_\la(t-z,t)\,.\]

Clearly, $\nu_\la^t$ is also a Poisson process of intensity $\la$, independent of $\nu_1$. Using  \eqref{eq:upboundLLla}, we can bound the left hand side of \eqref{eq:lower1} by
$$\P\left(\sup_{z\in(\eps t^{2/3},rt^{2/3}]} \big\{\nu_1(z) - \nu_\la^t(z)\big\} < \delta t^{1/3}\right) + \P\left(\bar Z_\la(t,t)< 0\right)\,,$$
that is equal to
\begin{equation}\label{eq:lower11}
\P\left(\sup_{z\in(\eps t^{2/3},rt^{2/3}]} \big\{\nu_1(z) - \nu_\la(z)\big\} < \delta t^{1/3}\right) + \P\left(\bar Z_\la(t,t)< 0\right)\,.
\end{equation}
(Notice that we also have restricted $z$ to $(\eps t^{2/3},rt^{2/3}]$).\\

By Lemma \ref{lem:lambda}, we see that
$$\P\left(\bar Z_\la(t,t)< 0\right) = O(r^{-3})\,,$$
for large enough $t$ and $r\in[1,t^{1/3})$. This means that we can take $r$ big enough, such that the second term in \eqref{eq:lower11} is smaller than $\eta/2$.\\

Now consider the first term in \eqref{eq:lower11}. After rescaling, the term can be written as
\begin{eqnarray*}
&&\P\left(\sup_{z\in(\eps t^{2/3},rt^{2/3}]} \big\{\nu_1(z) - \nu_\la(z)\big\} < \delta t^{1/3}\right) \\
&&\qquad\qquad\qquad = \,\P\left(\sup_{z\in(\eps, r]}  \big\{t^{-1/3}[\nu_1(zt^{2/3}) - \nu_\la(zt^{2/3})]\big\} < \delta\right)\,.
\end{eqnarray*}

We know that the process 
$$\{z\mapsto t^{-1/3}[\nu_1(zt^{2/3}) - \nu_\la^t(zt^{2/3})]\}$$ 
converges, as $t\to \infty$, to a drifting Brownian motion process, namely
\[ W_r(z) := W(2z) - rz\  \ , z\geq 0\,,\]
where $W$ is  standard Brownian motion on $\R_+$. A standard application of Donsker's Theorem then shows that
\begin{eqnarray*}\label{eq:Donsker1}
 &&\lim_{t\to\infty} \P\left(\sup_{z\in(\eps t^{2/3},rt^{2/3}]} \big\{\nu_1(z) - \nu_\la^t(z)\big\} < \delta t^{1/3}\right)\\
  && \qquad\qquad=\,\P\left(\sup_{z\in(\eps,r]}\  W_r(z) \leq \delta\right)\\
 && \qquad\qquad\leq\,\P\left(\sup_{z\in[0,r]}\  W_r(z) \leq \delta\right)+ \P\left(\sup_{z\in[0,\eps]}\  W_r(z) > \delta\right)\,.
\end{eqnarray*}

Remember that we have fixed $r>0$ already. Now we choose $\delta=\delta(\eta,r)$ {small enough}, such that
\[ \P\left(\sup_{z\in[0,r]}\  W_r(z) \leq \delta\right) < \eta/4\,.\]
For fixed $r$ and $\delta(\eta,r)$, we can now choose $\eps>0$ small enough (as in the previous estimates), such that
\[ \P\left(\sup_{z\in[0,\eps]}\  W_r(z) > \delta\right) < \eta /4\,.\]

\hfill$\Box$\\

Now we are able to prove the main result of this section:
\begin{thm}\label{thm:lowboundZ}
\[\lim_{\eps\to 0} \limsup_{t\to\infty} \P\left(0\leq Z_1(t,t)\leq \eps t^{2/3}\right) = 0\,.\]
\end{thm}

\noindent{\bf Proof:}
Choose $\eta>0$. It is enough to find $\eps>0$ such that
\begin{equation}\label{eq:lower3}
\limsup_{t\to\infty} \P\left(0\leq Z_1(t,t)\leq \eps t^{2/3}\right) < \eta\,.
\end{equation}
The probability in the left hand side of \eqref{eq:lower3} is bounded by
\begin{eqnarray*}
&&\P\Big(\sup_{z>\eps t^{2/3}} \big\{\nu_1(z) + L((z,0),(t,t)))\big\}\\
&&\qquad\qquad\qquad\qquad< \sup_{0\leq z \leq \eps t^{2/3}} \big\{\nu_1(z) + L((z,0),(t,t))\big\}\Big)\,.
\end{eqnarray*}
For any given $\delta>0$, this is bounded by 
\begin{eqnarray*}
&& \P\left(\sup_{z>\eps t^{2/3}} \big\{\nu_1(z) - \Delta(z,t)\big\} < \delta t^{1/3}\right)\\
&&\qquad\qquad\qquad+\,\,\P\left(\sup_{0\leq z \leq \eps t^{2/3}} \big\{\nu_1(z) -\Delta(z,t)\big\} > \delta t^{1/3}\right)\,.
\end{eqnarray*}

We note that the random function $L(\cdot,\cdot)$ on $\R\times (0,\infty)$ is independent of $\nu_1$. Furthermore,
\[\{z\mapsto L(\0,(t,t))-L((z,0),(t,t))\} \stackrel{\cal D}{=}\{z\mapsto L(\0,(t,t))-L(\0,(t-z,t))\}\,\]
(as processes). Therefore, \eqref{eq:lower3} follows by combining Lemma \ref{lem:lower1} together with Lemma \ref{lem:lower2}.

\hfill$\Box$\\

\begin{coro}
\[ \liminf_{t\to\infty}\frac{\E Z^+_1(t,t)}{t^{2/3}} = \liminf_{t\to\infty}\frac{\Var(L_1(t,t))}{2t^{2/3}} > 0\,.\]
\end{coro}

\noindent{\bf Proof:} Since, with probability 1, either $\bar Z_1(t,t)\geq 0$ or $\bar Z_1(t,t)<0$, we have that
\begin{eqnarray*}
1&=&\P(\bar Z_1(t,t)\geq 0)+\P(\bar Z_1(t,t)<0)\\
&\leq &\P(\bar Z_1(t,t)\geq 0)+\P(\bar Z'_1(t,t)\leq0)\\
&=& 2\P(\bar Z_1(t,t)\geq 0)\,,
\end{eqnarray*}
and hence,
\begin{eqnarray*}
\frac12&\leq&\P(\bar Z_1(t,t)\geq 0)\\
&=&\P(Z_1(t,t)\geq 0)\\
&\leq & \P\left(0\leq Z_1(t,t)\leq \eps t^{2/3}\right)+\P\left(Z_1(t,t)> \eps t^{2/3}\right)\\
&\leq &\P\left(0\leq Z_1(t,t)\leq \eps t^{2/3}\right)+\frac{\E Z^+_1(t,t)}{\eps t^{2/3}}\,,
\end{eqnarray*}
which shows that
$$\epsilon\left[\frac12-  \P\left(0\leq Z_1(t,t)\leq \eps t^{2/3}\right)\right]\leq\frac{\E Z^+_1(t,t)}{ t^{2/3}}\,.$$
By Theorem \ref{thm:lowboundZ}, one can choose $\eps_0$ so that
$$\P\left(0\leq Z_1(t,t)\leq \eps_0 t^{2/3}\right)<\frac14\,,$$
for large enough $t$, which implies that
$$\frac{\eps_0}{4}<\frac{\E Z^+_1(t,t)}{t^{2/3}}\,,$$
for large enough $t$. Theorem \ref{thm:exitpoint} completes the proof of the corollary.

\hfill $\Box$\\

\section{Cube-Root Asymptotics for $L$}
In the preceding sections it was shown that, for the classical Hammersley model in equilibrium, the variance of of $L_1(t,t)$ is of order $t^{2/3}$. The key result to extend the cube-root asymptotics to $L$ is the following theorem:
\begin{thm}\label{thm:L_1L}
There exists a constant $c_1>0$ such that for all $r>0$
$$\P\left(|L_1(t,t)-L(t,t)|\geq r t^{1/3}\right)\leq \frac{c_1}{r^{5/4}}\,.$$
\end{thm}

Before proving Theorem \ref{thm:L_1L} we have:

\begin{coro}\label{coro:tight}
There exists a constant $c_2>0$ such that
$$\E|L(t,t)-2t|\leq c_2t^{1/3}\,.$$
\end{coro}

\noindent{\bf Proof:}
Indeed,
\begin{eqnarray}
\nonumber\E|L(t,t)-2t|&\leq& \E|L(t,t)-L_1(t,t)|+\E|L_1(t,t)-2t|\\
\nonumber &\leq&\left\{1+c_1\int_1^\infty r^{-5/4}dr\right\}t^{1/3}+\sqrt{\Var L_1(t,t)}\,.
\end{eqnarray}
By Corollary \ref{cor:var_expectation}, this proves Corollary \ref{coro:tight}.

\hfill $\Box$\\

To prove Theorem \ref{thm:L_1L} we require the following lemma, the proof of which uses again the local comparison argument.

\begin{lem}\label{lem:tightL1}
There exist constants $c_1,c_2>0$ such that for all $t>0$ and all $r\in[1,c_1t^{4/5}]$
\begin{equation}\label{eq:tightL1}
\P\Big(\sup_{|z|\leq r^{7/12}t^{2/3}}\{\nu_1(z)-\Delta(z,t)\}> rt^{1/3}\Big)\leq \frac{c_2}{r^{5/4}}\,.
\end{equation}
\end{lem}

\noindent{\bf Proof:}
First, let $K=r^{7/12}$. As in the proof of Lemma \ref{lem:lower2}, for $\la=1-ut^{-1/3}$ we get that the left hand side of \eqref{eq:tightL1} is bounded by
$$ \P\left(\sup_{|z| \leq K t^{2/3}} \big\{\nu_1(z) - \nu_{\la}(z)\big\} > r t^{1/3}\right) + \P\left(\bar Z'_{\la}(t,t)>0\right)\,,$$
and
$$\P\left(\bar Z'_{\la}(t,t)>0\right) = O(u^{-3})\,,$$
uniformly for $u\in[1,t^{1/3}/2]$. To deal with the first term, we first note that w.l.g. we can restrict $z$ to $[0,Kt^{2/3}]$ (for $z\in[-Kt^{2/3},0]$ the argument is similar). Also that
$$z\in[0,K]\mapsto M_{u}(z):= t^{-1/3}\big\{\nu_1(zt^{2/3}) - \nu_{\la}(zt^{2/3})\big\}-uz\,,$$
is a zero mean martingale, and that
$$\E(M_{u}(z)^2)\leq z+ (1-\lambda)z\leq 2z\,.$$
Hence,
 \begin{eqnarray*}
&&\P\left(\sup_{z \in[0, K t^{2/3}]} \big\{\nu_1(z) - \nu_{\la}(z)\big\} > r t^{1/3}\right)\\
&&\qquad\qquad\qquad\qquad=\, \P\left(\sup_{z \in[0, K]} \big\{M_{u}(z)+uz\big\} > r\right)\\
&&\qquad\qquad\qquad\qquad\leq \,\P\left(\sup_{z \in[0, K]} \big\{M_{u}(z)\big\} > r-uK\right)\,.
 \end{eqnarray*}
Taking $u=r/(2K)$ and using Doob's sub-martingale inequality, we get
 \begin{eqnarray*}
&&\P\left(\sup_{z \in[0, K t^{2/3}]} \big\{\nu_1(z) - \nu_{\la}(z)\big\} > r t^{1/3}\right)\\
&&\qquad\qquad\qquad\qquad\qquad\leq\, \P\left(\sup_{z \in[0, K]} \big\{M_{u}(z)\big\} > r/2\right)\\
 &&\qquad\qquad\qquad\qquad\qquad=\,\P\left(\sup_{z \in[0, K]} \big\{M_{u}(z)^2\big\} > r^2/4\right)\\
&&\qquad\qquad\qquad\qquad\qquad\qquad\leq\,\frac{4\E (M_{u}(K)^2)}{r^2}\leq \frac{8K}{r^2}<\frac{8}{r^{4/3}}\,.\\
 \end{eqnarray*}
This means that, for $u\sim rK^{-1}=r^{5/12}\in[1,t^{1/3}/2]$, our estimate for the second term is dominant, and the lemma is proved.

\hfill $\Box$\\

\noindent{\bf Proof of Theorem \ref{thm:L_1L}:} By Corollary \ref{coro:exitail}
\begin{eqnarray*}
&&\P\left(L_1(t,t)\neq\sup_{|z|\leq Kt^{2/3}}\{\nu_1(z)+L((0,z),(t,t))\}\right)\\
&&\qquad\qquad\qquad\qquad\qquad\leq\, 2\P\left(Z_1(t,t)>Kt^{2/3}\right)=\,O(K^{-3})\,.
\end{eqnarray*}
So taking $K=r^{7/12}$ and using Lemma \ref{lem:tightL1} (and translation invariance) we have, for all $r\in[1,c_1t^{4/5}]$, 
\begin{eqnarray*}\label{eq:tightL}
&&\P\left( |L_1(t,t)-L(\0,(t,t))|\geq rt^{1/3}\right)\\
&&\qquad\qquad\qquad=\,\P\left(L_1(t,t)-L(\0,(t,t))\geq rt^{1/3}\right)\\
&&\qquad\qquad\qquad=\,\P\left(\sup_{|z|\leq r^{7/12}t^{2/3}}\{\nu_1(z)-\Delta(z,t)\}\geq rt^{1/3}\right)\\
&&\qquad\qquad\qquad\qquad\qquad\qquad\qquad\qquad+\, 2\P\left(Z_1(t,t)>Kt^{2/3}\right)\\
&&\qquad\qquad\qquad\qquad\quad\qquad=\,O(r^{-5/4})+O(r^{-7/4})\,.
\end{eqnarray*}

If $r> c_1 t^{4/5}$, we note that
\begin{eqnarray}
\nonumber\P\left(L_1(t,t)-L(\0,(t,t))\geq rt^{1/3}\right)&\leq&\P\left(L_1(t,t)\geq rt^{1/3}\right)\\
\nonumber&\leq&2\P\left(P_t\geq \frac12rt^{1/3}\right)\\
\nonumber&\leq&2\P\left(P_t\geq \frac{c_1^{5/6}}{2}r^{1/6}t\right)\,,
\end{eqnarray}
where $P_t$ is a Poisson random variable with expectation $t$. This implies that
$$\P\left(L_1(t,t)-L(\0,(t,t))\geq rt^{1/3}\right)$$
 tends to zero faster than any negative power of $r$, if $r>t^{4/5}$, uniformly in all large $t$. Hence, we can conclude the proof of Theorem \ref{thm:L_1L}.

\hfill $\Box$\\

\cleardoublepage

\end{document}